%% file: ell_2_PTRF.tex
\documentclass[11pt, a4paper]{article}

\input{packages_notes}

\input{commandes_guillaume}
\usepackage{braket}
\usepackage{authblk}
\begin{document}
\title{A geometrical viewpoint on the benign overfitting property of the minimum $\ell_2$-norm interpolant estimator and its universality.}
\author[1]{Guillaume Lecu{\'e} and Zong Shang  \\ email: \href{mailto:lecue@essec.edu}{lecue@essec.edu}, email: \href{mailto:zong.shang@ensae.fr}{zong.shang@ensae.fr} \\ESSEC, business school, 3 avenue Bernard Hirsch, 95021 Cergy-Pontoise, France.\\ CREST Laboratory, CNRS, GENES, Ecole Polytechnique, Institut Polytechnique de Paris, 5, avenue Henry Le Chatelier 91120 Palaiseau, France.}

\maketitle

\begin{abstract}

In the linear regression model, the minimum $\ell_2$-norm interpolant estimator $\hat\bbeta$ has received much attention since it was proved to be consistent even though it fits noisy data perfectly under some condition on the covariance matrix $\Sigma$ of the input vector, known as \textit{benign overfitting}. Motivated by this phenomenon, we study the generalization property of this estimator from a geometrical viewpoint. Our main results extend and  improve the convergence rates as well as the deviation probability from \cite{TB21}.  Our proof differs from the classical bias/variance analysis and is based on the \textit{self-induced regularization} property introduced in  \cite{BMR}: $\hat\bbeta$ can be written as a sum of a ridge estimator $\hat\bbeta_{1:k}$ and an overfitting component $\hat\bbeta_{k+1:p}$ which follows a decomposition of the features space $\bR^p=V_{1:k}\oplus^\perp V_{k+1:p}$ into the space $V_{1:k}$ spanned by the top $k$ eigenvectors of $\Sigma$ and $V_{k+1:p}$ spanned by the $p-k$ last ones. We also prove a matching lower bound for the expected prediction risk thus obtain the sufficient and necessary conditions for benign overfitting of $\hat\bbeta$.   The two geometrical properties of random Gaussian matrices at the heart of our analysis are the Dvoretsky-Milman theorem and  isomorphic and restricted isomorphic properties. In particular, the Dvoretsky dimension appearing naturally in our geometrical viewpoint, coincides with the effective rank from \cite{MR4263288,TB21} and is the key tool for handling the behavior of the design matrix restricted to the sub-space $V_{k+1:p}$ where overfitting happens.

We extend these results to heavy-tailed scenarii proving the universality of this phenomenon beyond exponential moment assumptions. This phenomenon is unknown before and is widely believed to be a significant challenge. This follows from an anistropic version of the probabilistic Dvoretsky-Milman theorem that holds for heavy-tailed vectors which is of independent interest.
\end{abstract}

\section{Introduction} 
\label{sec:introduction}

One side of the 'deep learning experiment' that has received much attention from the statistical community over the last five years is about an unexpected behavior of some deep learning models that  fit (potentially noisy) data perfectly but can still generalize well \cite{MR3997901,DBLP:journals/cacm/ZhangBHRV21,DBLP:conf/iclr/ZhangBHRV17}. This is a bit unexpected because it is quite rare to find a book on statistical learning written before the 'era of deep learning', that does not advise any engineer or student to avoid an algorithm that would rely too much on training data. Overfitting on the training dataset was  something that any statistician  would try to avoid or to correct by (explicit) regularization when designing their learning models. This principle has been profoundly rethought over the last five years.  

This property is nowadays called the \textit{benign overfitting} (BO) phenomenon \cite{DBLP:journals/corr/abs-2105-14368,BMR} and has been the subject of much recent work in machine learning and the statistical communities. The motivation is to identify situations where benign overfitting holds, i.e. when an estimator with a perfect fit on the training data (called interpolant estimator) can still generalize well. Several results have been obtained on BO in various statistical frameworks such as regression \cite{DBLP:conf/nips/BelkinHM18,DBLP:conf/aistats/BelkinRT19,DBLP:journals/jsait/MuthukumarVSS20,DBLP:journals/simods/BelkinHX20,DBLP:journals/corr/abs-1903-08560,DBLP:journals/corr/abs-1906-03667,DBLP:conf/nips/XuH19,DBLP:journals/corr/abs-2012-00807,DBLP:journals/corr/abs-2112-04470,MR4124325,Mei_Montanari,Chinot_Lerasle,DBLP:conf/colt/ZouWBGK21,steph_BO,DBLP:journals/jsait/MuthukumarVSS20,DBLP:journals/corr/abs-2106-09276} and classification \cite{DBLP:conf/nips/BelkinHM18,MR1673273,DBLP:journals/corr/abs-2105-02083,DBLP:journals/corr/abs-2101-11815,DBLP:journals/jmlr/WynerOBM17}. The whole point is to find conditions under which interpolant estimators can generalize well (that is to predict well on out-of-sample data). i.e. to prove that their excess risk or estimation risk tend to zero -- in other words, when interpolant estimators are consistent. 

In the linear regression model, it has been understood that BO  principally occurs for anisotropic design vectors having some special decay properties on their spectrum \cite{MR4263288,TB21,DBLP:journals/corr/abs-2112-04470} except for the result from \cite{DBLP:journals/corr/abs-2111-05987} where BO was achieved by the minimum $\ell_1$-norm interpolant estimator in the isotropic case. 

\subsection{The anistropic Gaussian design linear regression model with Gaussian noise}\label{subsec:model_gauss}
Let $N\in\bN$ be the number of samples and $p\in\bN\cup\{\infty\}$ be the dimension of the features space. We consider $y=\bX \bbeta^*+\bxi$ where $\bX:\bR^p\to \R^N$ is a Gaussian matrix with i.i.d. $\cN(0,\Sigma)$ row vectors, $\bxi$ is an independent Gaussian noise $\cN(0,\sigma_\xi^2 I_N)$ and $\bbeta^*\in\bR^p$ is the unknown parameter of the model.   We write 
\begin{equation}\label{eq:design_matrix}
\bX = \begin{pmatrix}
(\Sigma^{1/2}G_1)^\top\\
\vdots\\
(\Sigma^{1/2}G_N)^\top
\end{pmatrix} = \bG^{(N\times p)}\Sigma^{1/2} 
\end{equation}where $G_1, \ldots, G_N$ are $N$ i.i.d. $\cN(0,I_p)$ and $\bG^{(N\times p)}$ is a $N\times p$ standard Gaussian matrix (with i.i.d. $\cN(0,1)$ entries). We consider this model as a benchmark model because it is likely not reflecting real world data but it is the one that is expected to be universal in the sense that, results obtained in other more realistic statistical models could be compared with or tend to that one obtained in this ideal benchmark Gaussian model. 

We will indeed prove the universality of the Gaussian case  in Section~\ref{sec:the_heavy_tailed_case} where heavy-tailed scenarios are considered. A taste of this section may be found in Theorem~\ref{theo:main_intro_weak_moment} where the very same rate of convergence as in the Gaussian case is obtained under the only existence of $\log N$ moments of the marginals of the design vector and of only $r>4$ moments on the noise.

The relevance of the approximation of large neural networks by linear models via the neural tangent kernel feature map  \cite{DBLP:conf/nips/JacotHG18,DBLP:conf/stoc/JacotGH21} in some regimes has been much discussed in the machine learning community for instance in \cite{DBLP:journals/corr/abs-2105-14368,DBLP:conf/nips/0001ZB20,MR4263288} and is also a motivation for the study of this model.

\subsection{The minimum $\ell_2$-norm interpolant estimator} Interpolant estimators are estimators $\hat\bbeta$ that perfectly fit the data, i.e. such that $\bX \hat\bbeta=y$.
We denote by  $\norm{\cdot}_2$ the Euclidean norm on $\bR^p$ and consider the interpolant estimator having the smallest $\ell_2$-norm
\begin{equation}\label{eq:interpolant_esti}
\hat\bbeta\in\argmin_{\bbeta:\bX \bbeta = y}\norm{\bbeta}_2.
\end{equation}The relevance of the minimum $\ell_2$-norm interpolant estimator for neural network models is that, gradient descent algorithms converges to such estimators in some regimes (see Section~3.9 from \cite{DBLP:journals/corr/abs-2105-14368} and Section~3 from \cite{BMR} and references therein).

The goal is to identify situations (that is covariance matrices $\Sigma$ and signal $\bbeta^*$) for which $\hat\bbeta$ generalizes well even though it perfectly fits the data. The generalization capacity of an estimator $\hat\bbeta$ is measured via the excess risk:
\begin{equation*}
\bE[(Y-\inr{X,\hat\bbeta})^2|(y,\bX)] -\bE(Y-\inr{X,\bbeta^*})^2 = \norm{\Sigma^{1/2}(\hat \bbeta - \bbeta^*)}_2^2
\end{equation*}where $Y=\inr{X,\bbeta^*}+\xi$, $X\sim\cN(0,\Sigma)$ and $\xi\sim\cN(0,\sigma_\xi^2)$ is independent of $X$. In this model, the excess risk is an estimation risk with respect to the norm $\norm{\Sigma^{1/2}\cdot}_2$. Our goal is to prove matching high probability upper and in expectation lower bounds on the generalization error $\norm{\Sigma^{1/2}(\hat \bbeta - \bbeta^*)}_2^2$ of the minimum $\ell_2$-norm interpolant estimator $\hat \bbeta$ from \eqref{eq:interpolant_esti}. This will allow to identify cases (i.e. $\Sigma$ and $\bbeta^*$) where this bound (i.e. the rate of convergence) tends to zero as $N$ and $p$ go to infinity. These cases are all situations where the benign overfitting phenomenon takes place. The study of the BO requires making both $N$ and $p$ tend to infinity; it means that the number $p$ of features also increases with the sample size $N$. A way such asymptotic can be performed is by adding features in such a way that only the tail of the spectrum of $\Sigma$ is modified by adding eigenvalues smaller than the previous one as $p$ tends to infinity.

\subsection{Organization of the paper} 
\label{sub:organization_of_the_paper}
The paper is organized as follows. In the next section, we provide the two main geometrical tools we will use to prove the main results in the Gaussian case. In the third section, we provide a high-level description of the  proofs of our main results. In Section~\ref{sec:main_results} and \ref{sub:lower_bound_on_the_prediction_risk_of_}, we provide the four main results in the Gaussian case, identify necessary and sufficient conditions for benign overfitting in Definition~\ref{def:BO} and introduce a self-adaptive property of $\hat\bbeta$. Section~\ref{sec:the_heavy_tailed_case} is dedicated to the heavy-tailed case and the universality of the Gaussian case. All proofs are provided in the appendix.

\subsection{Notation}\label{sec:notation}
For all $q\geq 1$, $\norm{\cdot}_q$ denotes the $\ell_q$-norm in $\bR^p$, $B_q^p$ its unit ball and $\cS_q^{p-1}$ its unit sphere. The operator norm of a matrix is denoted by   $\norm{\cdot}_{op}$. For $A\in\bR^{m\times n}$, $s_1(A)\geq \cdots \geq s_{m\wedge n}(A)$ denote the singular values of $A$, its condition number is $s_1(A)/s_{m\wedge n}(A)$ and its spectrum is the vector ${\rm spec}(A)=(s_1(A), \cdots, s_{m\wedge n}(A))$. The set of integers $\{1, \ldots, p\}$ is denoted by  $[p]$. Let $\Gamma$ be a $p\times p$ semi-definite (not necessarily positive) matrix. We denote the generalized inverse of this matrix by $\Gamma^{-1}$. Recall that $\Sigma\in\bR^{p\times p}$ is the covariance matrix of design vector $X$. The SVD of $\Sigma$ is $\Sigma = UDU^\top$ and $(u_1,\cdots,u_p)$ are the row vectors of $U^\top$ and $D = \diag\left(\sigma_1,\cdots,\sigma_p\right)$ is the $p\times p$ diagonal matrix with $\sigma_1\geq \cdots\geq \sigma_p>0$ such that $\Sigma u_j = \sigma_j u_j$ for all $j\in[p]$. For all $k\in[p]$ we denote $V_{1:k} = {\rm span}(u_1, \ldots, u_k) \mbox{ and } V_{k+1:p} = {\rm span}(u_{k+1}, \ldots, u_p)$. We denote by $P_{1:k}:\bR^p\mapsto \bR^p$ ($P_{k+1:p}$ resp.) the orthogonal projection onto $V_{1:k}$ ($V_{k+1:p}$ resp.) and for all $\bbeta\in\bR^p$, we denote $\bbeta_{1:k} := P_{1:k}\bbeta$ and $\bbeta_{k+1:p}:= P_{k+1:p}\bbeta$. We denote $\Sigma_{1:k} = U D_{1:k} U^\top \mbox{ and } \Sigma_{k+1:p} = U D_{k+1:p} U^\top$ where $ D_{1:k}  = {\rm diag}(\sigma_1, \ldots, \sigma_k, 0,\ldots, 0)$ and  $ D_{k+1:p}  = {\rm diag}(0,\ldots, 0,\sigma_{k+1}, \ldots, \sigma_p)$. We also denote $X_{1:k} = \bG^{(N\times p)} \Sigma_{1:k}^{1/2}$ and  $X_{k+1:p} = \bG^{(N\times p)} \Sigma_{k+1:p}^{1/2}$ so that $\bX = X_{1:k} + X_{k+1:p}$. The two effective ranks used in \cite{MR4263288,TB21} are given by 
\begin{equation}\label{eq:effective_ranks}
r_k(\Sigma) =  \frac{\Tr(\Sigma_{k+1:p})}{\norm{\Sigma_{k+1:p}}_{op}} \mbox{ and } R_k(\Sigma) = \frac{(\Tr(\Sigma_{k+1:p}))^2}{\Tr(\Sigma_{k+1:p}^2)}.
\end{equation}

Let $\vertiii{\cdot}$ be some norm on $\bR^p$, denote by $B$ its unit ball, by $\vertiii{\cdot}_*$ its associated dual norm and by $B^*$ its unit dual ball so that for all $x\in\bR^p$, $\vertiii{x}=\sup(\inr{x,y},y\in B^*)$. The \underline{Gaussian mean width} of $B^*$ is $\ell^*(B^*)=\bE \sup(\inr{x,G}:x\in B^*)$ where $G\sim \cN(0,I_p)$. For a set $T$ in a linear space, $\conv(T)$ denotes the convex hull of $T$.

Let $k\in[p]$ and $\kappa_{DM}$ an absolute constant. Unless otherwise specified, we let
\begin{equation*}
J_1:=\left\{j\in[k]:\, \sigma_j\geq \frac{\kappa_{DM}\ell_*^2(\Sigma^{1/2}_{k+1:p}B_2^p)}{N}\right\},\quad J_2 := [k]\backslash J_1
\end{equation*}and  $\Sigma_{1,thres}^{-1/2} := U D_{1,thres}^{-1/2} U^\top$ where $U$ is the orthogonal matrix appearing in the SVD of $\Sigma$ and $D_{1,thres}^{-1/2}$ is the $p\times p$ diagonal matrix
\begin{equation*}
{\rm diag}\left( \left(\sigma_1 \vee\frac{\kappa_{DM}\ell_*^2(\Sigma^{1/2}_{k+1:p}B_2^p)}{N} \right)^{-1/2}, \ldots, \left(\sigma_k \vee\frac{\kappa_{DM}\ell_*^2(\Sigma^{1/2}_{k+1:p}B_2^p)}{N} \right)^{-1/2}, 0, \ldots, 0\right).
\end{equation*}


\subsection{Two seminal results on benign overfitting in linear regression} We state the main upper bound results  from \cite{MR4263288} and \cite{TB21} in the Gaussian linear model introduced above and for the minimum $\ell_2^p$-norm interpolant estimator \eqref{eq:interpolant_esti},  even though they have been obtained under weaker assumptions such as sub-gaussian assumptions on the design and the noise or for the ridge estimator.

\begin{Theorem}\label{theo:BLLT}[Theorem~4 in \cite{MR4263288}]
There are absolute constants $b,c,c_1>1$ for which the following holds. Define
\begin{equation}\label{eq:k_star_BLLT}
k^*_b = \min\left(k\geq0 : r_k(\Sigma)\geq b N\right),
\end{equation}where the infimum of the empty set is defined as $+\infty$. Suppose $\delta<1$ with $\log(1/\delta)<N/c$. If $k_b^*\geq N/c_1$ then $\bE \norm{\Sigma^{1/2}(\hat \bbeta - \bbeta^*)}_2\geq \sigma_\xi/c$. Otherwise, 
\begin{equation}\label{eq:upper_bound_BLLT}
\begin{aligned}
\norm{\Sigma^{1/2}(\hat \bbeta - \bbeta^*)}_2\leq &c \norm{\bbeta^*}_2 \sqrt{\norm{\Sigma}_{op}} \max\left(\left(\frac{r_0(\Sigma)}{N}\right)^{\frac14}, \sqrt{\frac{r_0(\Sigma)}{N}}, \left(\frac{\log(1/\delta)}{N}\right)^{\frac14}\right) \\
& + c \sqrt{\log(1/\delta)}\sigma_\xi\left(\sqrt{\frac{k^*_b}{N}} + \sqrt{\frac{N}{R_{k^*_b}(\Sigma)}}\right)
\end{aligned}
\end{equation}with probability at least $1-\delta$, and $\bE \norm{\Sigma^{1/2}(\hat \bbeta - \bbeta^*)}_2\geq (\sigma_\xi/c)\left(\sqrt{k^*_b/N} + \sqrt{N/R_{k^*_b}(\Sigma)}\right)$.
\end{Theorem}

Theorem~\ref{theo:BLLT} is one of the first results proving that the BO phenomenon in the linear regression model can happen by identifying situations where the upper bound in \eqref{eq:upper_bound_BLLT} tends to zero: quoting \cite{MR4263288}:  '\textit{We say that a sequence of covariance operator $\Sigma$ is benign if
\begin{equation}\label{def:benign_BLLT}
\lim_{N,p\to +\infty} \frac{r_0(\Sigma)}{N} = \lim_{N,p\to +\infty} \frac{k^*_b}{N} = \lim_{N,p\to+\infty}\frac{N}{R_{k^*_b}(\Sigma)}=0.
\end{equation}}

The analysis of the BO from Theorem~\ref{theo:BLLT} is based uniquely on the behavior of the spectrum of $\Sigma$. In particular, it does not depend on the signal $\bbeta^*$ (even though the norm $\norm{\bbeta^*}_2$ appears as a multiplying factor in \eqref{eq:upper_bound_BLLT} and this quantity may be large since $\bbeta^*$ is a vector in $\bR^p$ with $p$ being large to allow for over-parametrization). This analysis was improved and generalized to any interpolant estimator minimizing a general norm (not only the $\ell_2$-norm) in \cite{DBLP:journals/corr/abs-2106-09276}. The latter paper also introduced a tool known as the CGMT (i.e. convex gaussian minmax Theorem). This tool revealed to be efficient to study interpolant estimators \cite{DBLP:journals/corr/abs-2111-05987,DBLP:conf/icml/DonhauserRSY22} for general norms in the Gaussian case. However, the analysis from \cite{DBLP:journals/corr/abs-2106-09276} did not catch the sharp rate of convergence in the $\ell_2$ case and, because of the CGMT, it is limited to the Gaussian case. Finally, this analysis was improved in the subsequent work from \cite{TB21}. We recall the main result from \cite{TB21} only for the ridge estimator with regularization parameter $\lambda=0$ since in that case, it coincides with the minimum $\ell_2$-norm interpolant estimator $\hat\bbeta$ from \eqref{eq:interpolant_esti}.

\begin{Theorem}[\cite{TB21}]\label{theo:TB}There are absolute constants $c,c_0$ such that the following holds. Let $\delta<1-4\exp(-N/c_0)$ and $k\leq N/c_0$. We assume that with probability at least $1-\delta$, the condition number of $X_{k+1:p}X_{k+1:p}^\top$ is smaller than $L$. For all $t\in(1,N/c_0)$, with probability at least $1 - \delta - 20\exp(-t)$,
\begin{equation}\label{eq:upper_bound_TB}
\begin{aligned}
\norm{\Sigma^{1/2}(\hat \bbeta - \bbeta^*)}_2\lesssim &L^2 \left(\norm{\Sigma_{k+1:p}^{1/2}\bbeta_{k+1:p}^*}_2+ \norm{\Sigma_{1:k}^{-1/2}\bbeta_{1:k}^*}_2\frac{\Tr(\Sigma_{k+1:p})}{N}\right)  + \sigma_{\xi} \sqrt{t} L \left(\sqrt{\frac{k}{N}} + \frac{\sqrt{N \Tr(\Sigma_{k+1:p}^2)}}{\Tr(\Sigma_{k+1:p})}\right)
\end{aligned}
\end{equation}and 
\begin{equation}\label{eq:TB_reciproque_Pisier}
\frac{\Tr(\Sigma_{k+1:p})}{\norm{\Sigma_{k+1:p}}_{op}}\geq \frac{c N}{L}.
\end{equation}
\end{Theorem}

Comparing with Theorem~\ref{theo:BLLT}, Theorem~\ref{theo:TB} shows that this is not the $\ell_2$-norm of the entire vector $\bbeta^*$ that has to be paid by $\hat\bbeta$ in the bias term but the weighted $\ell_2$ norm $\norm{\Sigma_{k+1:p}^{1/2}\bbeta_{k+1:p}^*}_2$ of the signal $\bbeta^*$ restricted to $V_{k+1:p}$, the space spanned by the $p-k$ smallest eigenvectors of $\Sigma$ as well as $\norm{\Sigma_{1:k}^{-1/2}\bbeta_{1:k}^*}_2\Tr(\Sigma_{k+1:p})/N$. Theorem~\ref{theo:TB} also introduces a key idea that the features space $\bR^p$ should be decomposed  as $V_{1:k}\oplus^\perp V_{k+1:p}$. This decomposition of $\bR^p$ is actually associated with a decomposition of the estimator $\hat\bbeta$ into a \textit{prediction component} and an \textit{overfitting component} following the \textit{self-induced regularization} idea from \cite{BMR} (see also the 'spiked-smooth' estimates from \cite{DBLP:journals/jmlr/WynerOBM17}). This decomposition is at the heart of our analysis.

\subsection{Questions}\label{sec:questions}
It can be observed from Theorem~\ref{theo:TB} that for every value of $k$, there exists a decomposition of the feature space $\bR^p$. This decomposition leads to a decomposition of $(\Sigma,\bbeta^*)$ and thus provides an upper bound for the estimation error. Since $k$ is a free parameter, meaning that it is not included in the definition of $\hat \bbeta$, there exists an optimal value of $k$ that minimizes the estimation error bound for all potential values of $k$. This optimal value of $k$ determines the optimum decomposition of $(\Sigma,\bbeta^*)$. A priori, the optimal choice of $k$ should depend on both $\Sigma$ and $\bbeta^*$, as the upper bound shown in Equation~\eqref{eq:upper_bound_TB} is a function determined by both $\Sigma$ and $\bbeta^*$.
\begin{itemize}
    \item According to \cite{TB21}, the optimal selection of the parameter $k$ can be determined by $k_b^*$ as defined in Equation~\eqref{eq:k_star_BLLT}, which is commonly referred to as the ``effective dimension.'' and is \emph{independent} with $\bbeta^*$. The assertion is grounded in the demonstration that the upper bound of the variance term, as derived from the analysis conducted by \cite{TB21}, matches its lower bound. Additionally, the upper bound of the bias term matches its lower bound of that of a \emph{random} signal. In other words, \cite{TB21} proved that there exists \emph{an} unknown signal $\bbeta^*$ for which $k=k_b^*$ is optimal. Nevertheless, we propose the subsequent questions:
    \begin{center}
        If this selection of $k$ is also optimal for \emph{every} $\bbeta^*$? 
        
        Or, is there any $\bbeta^*$ such that the optimal $k$ depends also on $\bbeta^*$?
    \end{center}An answer to this question will result in the sufficient and necessary condition for benign over-fitting for \emph{every} $(\Sigma,\bbeta^*)$.
    \item The selection of $k = k_b^*$ as stated in \cite{TB21} yields an optimal decomposition for some signal among all possible decompositions of the type $\bR^p = V_{J}\oplus^\perp V_{J^c}$, where $J = \{1,\cdots,k\}$ and $k\lesssim N$. One may wonder however
    \begin{center}
        if there are better choices of $J$ (see \cite{DBLP:journals/corr/abs-2106-09276}), especially when $k\gtrsim N$?
    \end{center}
    \item So far all the results on the BO phenomenon  have been obtained under strong exponential moments or even the Gaussian assumptions except the recent journal version of \cite{TB21} which appears concurrently with the current paper. Nevertheless,
    \begin{center}
        Can we obtain BO under weak moments assumption? for both design and noise? in this heavy-tailed setup, is it possible to achieve the same rate of convergence as in the Gaussian case?
    \end{center}Positive answers to these questions will prove the universality (in a non-asymptotic way) of the Gaussian case.
\end{itemize}


\subsection{Our contribution} The summarization of our contributions is as follows:
\begin{itemize}
    \item We improve the rates from Theorem~\ref{theo:BLLT} and Theorem~\ref{theo:TB}, and improve the probability deviation from constant in \cite{TB21}, to exponentially high.
    \item We extend the results to $k\gtrsim N$.   This allows a possibly better choice to be made regarding the splitting parameter $k$ that can lower the price of overfitting (overfitting happens on the space $V_{k+1:p}$).
    \item We also extend the features space decomposition (see Section~\ref{sub:on_the_choice_of_}) beyond that of the form $\bR^p = V_{1:k}\oplus^\perp V_{k+1:p}$ used in Theorem~\ref{theo:TB} to any features space decomposition of the form $\bR^p = V_{J}\oplus^\perp V_{J^c}$ -- where $V_J = {\rm span}(u_j:j\in J)$ -- with no restriction on the dimension $|J|$ of $V_J$ (where estimation happens).
    \item We improve the lower bound results from \cite{MR4263288,TB21}, by removing unnecessary assumptions (see Section~\ref{sub:lower_bound_on_the_prediction_risk_of_} for more detail) as well as obtaining a result for the risk itself and not a Bayesian one as in \cite{TB21}. Indeed, we obtain a lower bound on the expected prediction risk which shows that the best choice of $J$ is $\{1, \ldots, k^*_b\}$ and yielding to matching upper and lower bounds (up to constants) and so to a necessary and sufficient condition for BO that we state in Definition~\ref{def:BO} below (not that this definition depends on both $\Sigma$ and $\bbeta^*$). As a result, we fully answer the first two questions proposed in subsection~\ref{sec:questions}.

    Unlike previous works, we highlight that the benign overfitting property of $\hat\bbeta$ depends on both $\Sigma$ and $\bbeta^*$ and, in particular, on their interplay, i.e. the coordinates of $\bbeta^*$ in a basis of eigenvectors of $\Sigma$ (and not only on the spectrum of $\Sigma$). 
    
    \item Our main contribution  also lies in a technical point of view: our proofs are based upon the \textit{self-regularization property} \cite{BMR} of $\hat\bbeta$ and not on a bias/variance decomposition as in \cite{TB21,MR4263288}. We also show that the Dvoretsky-Milman theorem is a key result to handle the behavior of the design matrix, restricted to $V_{k+1:p}$, where overfitting happens. In particular, we recover the effective rank $r_k(\Sigma)$ from \cite{MR4263288,TB21} as a natural geometric parameter of the problem: the Dvoretsky dimension of the ellipsoid $\Sigma_{k+1:p}^{-1/2}B_2^p$ (i.e. that associated with the covariance matrix $\Sigma$ restricted to the subspace of $\bR^p$ where overfitting happens). The other geometrical property that we introduce, is a restricted isomorphic property of the design matrix $\bX$, restricted to $V_{1:k}$ which is needed when $k\gtrsim N$. These two tools are at the basis of our 'geometrical viewpoint' on the benign overfitting phenomenon. They provide some insights on sufficient properties the design vectors $X_1, \ldots, X_N$ may satisfy to make BO possible.

    \item Finally, following our goal to universality of the Gaussian case in a non-asymptotic way,  we extend all these results to the heavy-tailed case. It is widely believed that going beyond the Gaussian case is a significant challenge. We show that BO does not require such strong concentration properties to hold. One key point in our results under heavy-tailed assumptions is that the rate and conditions are the same as in the Gaussian case. In particular, we don't pay any other type of statistical complexities than the Gaussian mean width. This type of result is important since it proves universality of the Gaussian setup for the BO phenomenon which was not known previously. However there is a technical cost to achieve this result since most of the geometrical properties that we are using to prove BO are known under only strong exponential moment or even the Gaussian assumption. In particular, we have to prove a probabilistic Dvoretsky-Milman theorem in the anisotropic and heavy-tailed case to prove BO under heavy-tailed cases. This result of independent interest may be found in appendix. With this universality result, we answer the last question proposed in subsection~\ref{sec:questions}.
\end{itemize}

\subsection{Our main result in the Gaussian case} 
\label{sub:our_main_result_in_the_gaussian_case}
In this section, we provide our findings in the Gaussian case (i.e. when the design vector and the noise have a Gaussian distribution as introduced in Section~\ref{subsec:model_gauss}) in cases where $k^*_b\lesssim N$ and $\sigma_1 N \gtrsim \Tr(\Sigma_{k^*_b+1:p})$ -- all other cases may be found below. It is a straightforward consequence of Theorem~\ref{theo:main} below (when $k=k^*_b$ and $b=4/\kappa_{DM}$, where $\kappa_{DM}$ is the absolute constant appearing in Dvoretsky-Milman's Theorem~\ref{theo:DM} below).  We set
\begin{equation}\label{eq:optimal_rate_4_terms}
r^*:=\max \bigg\{\sigma_\xi\sqrt{\frac{k^*_b}{N}}, \sigma_\xi\frac{ \sqrt{N\Tr(\Sigma_{k_b^*+1:p}^2)}}{\Tr(\Sigma_{k_b^*+1:p})},\norm{\Sigma_{k^*_b+1:p}^{1/2}\bbeta_{k^*_b+1:p}^*}_{2},\norm{\Sigma_{1:k^*_b}^{-1/2}\bbeta_{1:k^*_b}^*}_2\left(\frac{\Tr(\Sigma_{k^*_b+1:p})}{N}\right)\bigg\}.
\end{equation}

\begin{Theorem}\label{theo:main_intro_gauss}
There are absolute positive constants $c_0$, $c_1$, $c_2$ and $c_3$ such that the following holds. We assume that the design vector $X$ is distributed according to the Gaussian $\cN(0,\Sigma)$ and that the noise is independent of $X$ with a mean zero  and variance $\sigma_\xi$ Gaussian distribution. Let $b= 4/\kappa_{DM}$. We assume that $k^*_b\lesssim N$ and $\sigma_1 N \gtrsim  \Tr(\Sigma_{k^*_b+1:p})$. With probability at least $1-c_0\exp(-c_1 k^*_b)$, $\norm{\Sigma^{1/2}(\hat\bbeta-\bbeta^*)}_2\leq c_2 r^*$. Moreover, $\bE \norm{\Sigma^{1/2}(\hat\bbeta-\bbeta^*)}_2\geq c_3 r^*$.
\end{Theorem}

Comments and comparison with the existing literature for Theorem~\ref{theo:main_intro_gauss} are postponed after the more general results given below which are Theorem~\ref{theo:main}, \ref{theo:main_2} and \ref{theo:main_3} for the upper bounds and Theorem~\ref{theo:lower_bound} for the lower bound.

\subsection{Our main result in the heavy-tailed case} 
\label{sub:our_main_result_in_the_heavy_tailed_case}
The result in this section extends the previous one to a heavy-tailed case. The key feature of the following result is that the rate of convergence is the same as in the Gaussian case, it is up to an absolute constant given by $r^*$ in \eqref{eq:optimal_rate_4_terms} and that it holds under exactly the same condition on the number of data. The only difference with the Gaussian case is on the deviation probability parameter. All of that happens under the only existence of $\log N$ moments for all marginals of $X$ and $r>4$ moments for the noise.  More comments may be found after the more general versions of this result which are in Section~\ref{sec:the_heavy_tailed_case}.

\begin{Theorem}\label{theo:main_intro_weak_moment}
We assume that the design vector $X$ is such  that $X=\Sigma^{1/2}Z$ where $Z$ is a symmetric random variable with independent coordinates. We assume that there exists some $\alpha\leq 2$, $R>0$ and $L>0$ such that for all $2\leq q\leq R \log(N)$ and all $\bv\in\bR^p$, $\norm{\inr{X,\bv}}_{L_q}\leq L q^{1/\alpha}\norm{\inr{X,\bv}}_{L_2}$. We assume that the noise vector $\bxi = (\xi_i)_{i=1}^N$ is a random vector with independent mean zero and variance $\sigma_\xi$ coordinates such that for all $i$'s,  $\norm{\xi_i}_{L_r}\leq \kappa \sigma_\xi$ for some $\kappa>0$ and $r>4$.

 Then, there are  constants $\kappa_{iso}$, $\kappa_{DM}$, $c_0$, $c_1$, $c_2$ and $c_3$ (depending only on $\kappa, L$ and $R$) such that the following holds. If $k^*_b\leq \kappa_{iso} N$ for $b = 4/\kappa_{DM} $, then the following holds: with probability at least
 \begin{equation}
  1-\left(\frac{c_2}{k^*_b}\right)^{r/4} - \left(\frac{c_3 \log p}{N}\right)^{\frac{\log(2p)}{2}} - \frac{1}{N^{c_4}} - \left(\frac{c_3}{N}\right)^{\frac{R\log p - 1}{4}}- \frac{1}{k^*_b},
  \end{equation} 
\begin{align*}
\norm{\Sigma^{1/2}(\hat\bbeta-\bbeta^*)}_2\leq c_0 r^*
\end{align*}where $r^*$ is the rate obtained in the Gaussian case and defined in \eqref{eq:optimal_rate_4_terms}. 
\end{Theorem}

\section{Two geometrical properties of Gaussian matrices} 
\label{sec:geometrical_properties_of_gaussian_matrices}
We will use two properties of Gaussian matrices related to (approximately) Euclidean sections of convex bodies. The first one is the Dvoretsky-Milman's theorem (which implies the existence of Euclidean sections of a body via a section of that body by the span of a Gaussian matrix) and the second one is a restricted isomorphy property of Gaussian matrices (which implies Euclidean sections via the kernel of a Gaussian matrix). We recall these two tools because they play a central role in our analysis of the minimum $\ell_2$-norm interpolant estimator $\hat\bbeta$.

\subsection{The Dvoretsky-Milman's theorem for ellipsoids} \label{sec:DM}
 The original aim of the Dvoretsky-Milman (DM) theorem is to show the existence of Euclidean sections of a convex body $B$ and to identify the largest possible dimension (up to  a constant) of such Euclidean sections. Even though the terms of the problem are purely deterministic, the only way known so far to obtain such sections is via the range or the kernel of some random matrices such as standard Gaussian matrices. That is the reason why we will use this theorem for the properties of Gaussian random matrices, they also use.  We provide the general form of the DM theorem even though we will use it only for ellipsoids.

\begin{Definition}\label{def:dm_dimension}
  The \underline{Dvoretsky dimension} of $B$ is 
  \begin{equation*}
  d_*(B) = \left(\frac{\ell^*(B^*)}{{\rm diam}(B^*, \ell_2^p)}\right)^2
  \end{equation*}where ${\rm diam}(B^*, \ell_2^p) = \sup\left(\norm{v}_2:v\in B^*\right)$.
  \end{Definition}

  The Dvoretsky dimension $d_*(B)$ is up to absolute constants, the largest dimension of an Euclidean section of $B$ (see \cite{MR1036275} for both upper and lower bounds). For instance, $d_*(B_1^p) = p$ means that there exists a subspace $E$ of $\bR^p$ of dimension $c_0p$ where $c_0$ is some absolute constant such that $(c_1/\sqrt{p})B_2^p\cap E \subset B_1^p \cap E \subset (c_2/\sqrt{p})B_2^p\cap E$. In other words the section $B_1^p \cap E$ of the convex body $B_1^p$ by $E$ looks like the  Euclidean ball $B_2^p\cap E$ with radius $(1/\sqrt{p})$ up to absolute constants. Another example that we use below is the case of ellipsoids $\Sigma^{-1/2}B_2^p=\{v\in\bR^p:\norm{\Sigma^{1/2}v}_2\leq1\}$ where $\Sigma\in\bR^{p\times p}$ is a PSD matrix. In that case, the Dvoretsky's dimension of $\Sigma^{-1/2}B_2^p$ satisfies 
\begin{equation}\label{eq:dvoretsky_dim_ellipsoid}
 \frac{{\rm Tr}(\Sigma)}{4\norm{\Sigma}_{op}} \leq \frac{(\Tr(\Sigma)-2\norm{\Sigma}_{op})\vee \norm{\Sigma}_{op}}{\norm{\Sigma}_{op}}\leq  d_*(\Sigma^{-1/2}B_2^p) \leq \frac{{\rm Tr}(\Sigma)}{\norm{\Sigma}_{op}}
\end{equation}because $\ell^*((\Sigma^{-1/2}B_2^p)^*)=\ell^*(\Sigma^{1/2}B_2^p)  = \bE\norm{\Sigma^{1/2}G}_2\leq \sqrt{\bE\norm{\Sigma^{1/2}G}_2^2}= \sqrt{{\rm Tr}(\Sigma)}$ and $\bE\norm{\Sigma^{1/2}G}_2\geq \sqrt{\Tr(\Sigma)/2}$ (see Remark~1 in \cite{DBLP:journals/corr/abs-2106-09276}). The quantity appearing on the right-hand side of \eqref{eq:dvoretsky_dim_ellipsoid} has been used several times in the literature on interpolant estimators and was called there, the effective rank or effective dimension \cite{MR4263288,TB21} (see $r_k(\Sigma)$ in \eqref{eq:effective_ranks}). As a consequence, assuming that $ \norm{\Sigma}_{op} N\lesssim {\rm Tr}(\Sigma)$ is equivalent to assuming the existence of an Euclidean section of dimension on the order of $N$ of the ellipsoid $\Sigma^{-1/2}B_2^p$.

   Let us now state the general form of the Dvoretsky-Milman Theorem that will be useful for our purpose. It makes a standard Gaussian matrix appear explicitly. The following form of Dvoretsky theorem is due to V.Milman (up to the minor modification that \cite{MR0293374} uses the surface measure on the Euclidean sphere $\cS_2^{p-1}$ instead of the Gaussian measure on $\bR^p$).

  \begin{Theorem}\label{theo:DM}There are absolute constants $\kappa_{DM}\leq1$ and $c_0$ such that the following holds. Let $\vertiii{\cdot}$ be some norm on $\bR^p$ and denote by $B$ its unit ball. Denote by $\bG:=\bG^{(N\times p)}$, the $N\times p$ standard Gaussian matrix with i.i.d. $\cN(0,1)$ Gaussian entries.  Assume that $N\leq \kappa_{DM} d_*(B)$. Then with probability at least $1-\exp(-c_0 d_*(B))$, for every $\blambda\in\bR^N$, 
\begin{equation}\label{eq:dvoretsky}
\frac{1}{\sqrt{2}}\norm{\blambda}_2 \ell^*(B^*) \leq \vertiii{\bG^\top \blambda}\leq \sqrt{\frac{3}{2}}\norm{\blambda}_2 \ell^*(B^*). 
\end{equation}
  \end{Theorem}
 Classical proofs of the Dvoretsky-Milman (DM) theorem  hold only with constant probability as in \cite{MR3837109} or \cite{MR1036275} because DM's theorem is mainly used to prove an existence result: that ${\rm Im}(\bG^\top)$ realizes an Euclidean section of the convex body $B$. However, we will need DM's theorem to hold with large probability and it is straightforward to get the correct probability deviation as announced in Theorem~\ref{theo:DM} above, for instance from the proofs in \cite{MR3837109} or \cite{MR1036275}. Theorem~\ref{theo:DM}  also has implications on standard Gaussian matrices. We state such an outcome for anistropic Gaussian matrices that will be useful for our analysis.

 \begin{Proposition}\label{prop:DM_ellipsoid}
 Let $\bX_2 = \bG^{(N\times p)} \Gamma^{1/2}$ where $\bG^{(N\times p)}$ is an $N\times p$ standard Gaussian matrice and $\Gamma$ is a semi-definite matrix.  Assume that $N\leq \kappa_{DM} d_*(\Gamma^{-1/2}B_2^p)$. With probability at least $1-\exp(-c_0 d_*(\Gamma^{-1/2}B_2^p))$, $$\norm{\bX_2 \bX_2^\top - \ell^*(\Gamma^{1/2}B_2^p)^2I_N}_{op}\leq (1/2)\ell^*(\Gamma^{1/2}B_2^p)^2$$ which implies that
 \begin{align*}
 \sqrt{s_1(\bX_2 \bX_2^\top)} & = s_1(\bX_2)\leq \sqrt{3/2}\ell^*(\Gamma^{1/2}B_2^p)\leq \sqrt{3\Tr(\Gamma)/2}\\
  \sqrt{s_N(\bX_2 \bX_2^\top)} & = s_N(\bX_2)\geq (1/\sqrt{2})\ell^*(\Gamma^{1/2}B_2^p)\geq \sqrt{\Tr(\Gamma)}/2
 \end{align*}and 
 \begin{equation*}
 \frac{2}{ \sqrt{\Tr(\Gamma)}}\geq s_1[\bX_2^\top(\bX_2\bX_2^\top)^{-1}]\geq s_N[\bX_2^\top(\bX_2\bX_2^\top)^{-1}]\geq \sqrt{\frac{2}{3 \Tr(\Gamma)}}.
 \end{equation*}
 \end{Proposition}
 \beginproof

It follows from Theorem~\ref{theo:DM} for $\vertiii{\cdot} = \norm{\Gamma^{1/2}\cdot}_2$, (so that $B=\Gamma^{-1/2}B_2^p$ and $B^*=\Gamma^{1/2}B_2^p$) that with probability at least $1-\exp(-c_0 d_*(\Gamma^{-1/2}B_2^p))$,
\begin{equation*}
\norm{\bX_2 \bX_2^\top - \ell^*(\Gamma^{1/2}B_2^p)^2I_N}_{op} = \sup_{\norm{\blambda}_2=1}\left|\norm{\bX_2^\top \blambda}_2^2 - \ell^*(\Gamma^{1/2}B_2^p)^2\norm{\blambda}_2^2\right|\leq \frac{\ell^*(\Gamma^{1/2}B_2^p)^2}{2}.
\end{equation*}Moreover, we already saw that $\Tr(\Gamma)/2\leq\ell^*(\Gamma^{1/2}B_2^p)^2\leq \Tr(\Gamma)$ which proves the first statements. We finish the proof by using the fact that for all $j=1,\ldots, N, s_j[\bX_2^\top(\bX_2\bX_2^\top)^{-1}] = \sqrt{s_j[(\bX_2\bX_2^\top)^{-1}]}$.
 {\mbox{}\nolinebreak\hfill\rule{2mm}{2mm}\par\medbreak}

 A consequence of Proposition~\ref{prop:DM_ellipsoid} is that in Dvoretsky's regime (i.e. $N\leq \kappa_{DM} d_*(\Gamma^{-1/2}B_2^p)$), $\bX_2\bX_2^\top= \bG \Gamma \bG^\top$ (where $\bG=\bG^{(N\times p)}$) behaves almost like the homothety $\ell^*(\Gamma^{1/2}B_2^p)^2I_N$. DM's theorem holds for sufficiently small values of $N$; which means from a statistical viewpoint, that it is a high-dimensional property. 

\begin{Remark}[Effective rank and the Dvoretsky dimension of an ellispoid]
The effective rank $r_{k}(\Sigma)$ (see \eqref{eq:effective_ranks}) and the Dvoretsky dimension of $\Sigma_{k+1:p}^{-1/2}B_2^p$ are  the same quantity up to absolute constants. Note also that when the condition number $L$ from Theorem~\ref{theo:TB} is a constant, the lower bound from \eqref{eq:TB_reciproque_Pisier} is a consequence of a general result on the upper bound on the dimension of an Euclidean section of a convex body, which implies that Eqn. \eqref{eq:TB_reciproque_Pisier} holds  in the case of Ellipsoid (see Proposition~4.6 in \cite{MR1036275}). Indeed, when $L$ is like a constant, $X_{k+1:p}X_{k+1:p}^\top$ behaves like an isomorphy and so $\Sigma_{k+1:p}^{-1/2}B_2^p$ has an Euclidean section of dimension $N$ given by the range of $X_{k+1:p}^\top$. Hence, from Proposition~4.6 in \cite{MR1036275}, the Dvoretsky dimension of $\Sigma_{k+1:p}^{-1/2}B_2^p$  has to be larger than $N$ up to some absolute constant and so  Eqn. \eqref{eq:TB_reciproque_Pisier} holds.  The Dvoretsky-Milman theorem is the basis for several results previously used in literature about BO such as Lemma~9 and Lemma~10 in \cite{MR4263288} and Lemma~4.3 in \cite{BMR}.
\end{Remark}   

 There are cases where we want to avoid the condition '$N\leq \kappa_{DM} d_*(\Gamma^{-1/2}B_2^p)$' in the DM theorem. Since this assumption is only used to obtain the left-hand side of \eqref{eq:dvoretsky} (which is sometimes called a lower isomorphic property), we can state a result on the right-hand side inequality of \eqref{eq:dvoretsky} without this assumption. The proof of the following statement is standard and follows the same scheme as that of Theorem~\ref{theo:DM}. We did not find a reference for it and so we sketch its proof for completeness.

 \begin{Proposition}\label{prop:dvoretsky_upper_bound} Let $\vertiii{\cdot}$ be some norm on $\bR^p$ and denote by $B^*$ its unit dual ball. Denote by $\bG:=\bG^{(N\times p)}$ the $N\times p$ standard Gaussian matrix with i.i.d. $\cN(0,1)$ Gaussian entries.  With probability at least $1-\exp(-N)$, for every $\blambda\in\bR^N$, 
\begin{equation*}
\vertiii{\bG^\top \blambda}\leq 2\left(\ell^*(B^*) + \sqrt{2N(1+\log(10))}\diam(B^*, \ell_2^p)\right)\norm{\blambda}_2. 
\end{equation*}
 \end{Proposition}

 \textbf{Proof of Proposition~\ref{prop:dvoretsky_upper_bound}.} It is essentially the proof of the DM theorem which is followed here.  Denote by $G_1, \ldots, G_N$,  the $N$ i.i.d. $\cN(0,I_p)$ row vectors of $\bG$. Let $\blambda\in\cS_2^{N-1}$. It follows from Borell's inequality (see Theorem~7.1 in \cite{Led01} or p.56-57 in \cite{MR2814399} - the dual ball is compact and therefore separable), that with probability at least $1-\exp(-t)$, 
\begin{equation*}
 \vertiii{\bG^\top\blambda}=\vertiii{G}\leq \bE\vertiii{G} + \sqrt{2t}\sup_{x\in B^*}\sqrt{\bE \inr{G, x}^2} =  \ell^*(B^*) + \sqrt{2t}\diam(B^*, \ell_2^p)
\end{equation*}where $G=\sum_{i=1}^N \lambda_i G_i\sim \cN(0,I_p)$. Let $0<\eps<1/2$ and $\Lambda_\eps$ be an $\eps$-net of $\cS_2^{N-1}$ w.r.t. the $\ell_2^N$-norm. Let $\blambda^*\in \cS_2^{N-1}$ be such that $\vertiii{\bG^\top\blambda^*} = \sup\big(\vertiii{\bG^\top \blambda}:\blambda\in\cS_2^{N-1}\big):=M$ and let $\blambda^*_\eps\in \Lambda_\eps$ be such that $\norm{\blambda^* - \blambda^*_\eps}_2\leq \eps$. We have
\begin{equation*}
M = \vertiii{\bG^\top\blambda^*}\leq \vertiii{\bG^\top\blambda^*_\eps} + \vertiii{\bG^\top(\blambda^*-\blambda^*_\eps)}\leq \sup_{\blambda_\eps\in \Lambda_\eps}\vertiii{\bG^\top\blambda_\eps} + M\eps
\end{equation*}and so $M\leq (1-\eps)^{-1}\sup_{\blambda_\eps\in \Lambda_\eps}\vertiii{\bG^\top\blambda_\eps}$.
It follows from the bound above, derived from Borell's inequality, and a union bound, that with  probability at least $1-|\Lambda_\eps|\exp(-t)$, $M\leq (1-\eps)^{-1}\left(\ell^*(B^*) + \sqrt{2t}\diam(B^*, \ell_2^p)\right)$. It follows from a volume argument \cite{MR1036275} that $|\Lambda_\eps|\leq (5/\eps)^N$ hence, for $\eps=1/2$ and $t=N(1+\log(10))$ with probability at least $1-\exp(-N)$, 
\begin{equation*}
\sup\left(\vertiii{\bG^\top \blambda}:\blambda\in\cS_2^{N-1}\right)\leq 2\left(\ell^*(B^*) + \sqrt{2N(1+\log(10))}\diam(B^*, \ell_2^p)\right).
\end{equation*}
{\mbox{}\nolinebreak\hfill\rule{2mm}{2mm}\par\medbreak}

\subsection{Isomorphy and Restricted isomorphy properties} We use our budget of data, i.e. $N$, to ensure some isomorphic or restricted isomorphic property of a matrix  $\bX_1=\bG^{(N\times p)}\Gamma^{1/2}$ for some semi-definite matrix $\Gamma$. Let $B$ be the unit ball of some norm $\norm{\cdot}$. Let $\kappa_{RIP}$ be some absolute constant. We define the following family of fixed points: for all $\rho>0$,
\begin{equation}\label{eq:fixed_point}
r(\rho) = \inf\left(r>0:\ell^*\left((\rho \Gamma^{1/2}B)\cap (r B_2^p)\right)\leq \kappa_{RIP} r \sqrt{N}\right)
\end{equation}where $\ell^*\left((\rho \Gamma^{1/2}B)\cap (r B_2^p)\right)$ is the Gaussian mean width of $(\rho \Gamma^{1/2}B)\cap (r B_2^p)$. We see that for all $\rho>0$, $r(\rho) = \rho R_N(\Gamma^{1/2}B)$ where
\begin{equation}\label{eq:fixed_point_without_rho}
R_N(\Gamma^{1/2}B) = \inf\left(R>0:\ell^*\left((\Gamma^{1/2}B)\cap (R B_2^p)\right)\leq \kappa_{RIP} R \sqrt{N}\right).
\end{equation}The relevance of the fixed point $R_N(\Gamma^{1/2}B)$ may be seen in the following result. This result follows from a line of research in empirical process theory on the quadratic process \cite{MR2149924,MR2373017,MR3565471,dirksen2015tail,bednorz}.

\begin{Theorem}\label{theo:rip}
There are absolute constants $0<\kappa_{RIP}<1$ and $c_0$ such that for the fixed point $R_N(\Gamma^{1/2}B)$ defined in \eqref{eq:fixed_point_without_rho} and $\bX_1=\bG^{(N\times p)}\Gamma^{1/2}$, the following holds. With probability at least $1-\exp(-c_0N)$, 
\begin{equation}\label{eq:full_isomorphy}
\frac{1}{2}\norm{\Gamma^{1/2} v}_2^2\leq \frac{1}{N}\norm{\bX_1 v}_2^2\leq \frac{3}{2}\norm{\Gamma^{1/2} v}_2^2,
\end{equation}for all $v\in\bR^p$ such that $R_N(\Gamma^{1/2}B)\norm{v}\leq \norm{\Gamma^{1/2}v}_2$.
\end{Theorem}

\textbf{Proof of Theorem~\ref{theo:rip}.} For all $v\in\bR^p$, we set $f_v:x\in\bR^p\mapsto \inr{x,v}$ and denote by $F:=\{f_v:v\in\bR^p,  \norm{v}\leq R_N(\Gamma^{1/2}B)^{-1} \mbox{ and }\norm{\Gamma^{1/2}v}_2=1\}$ when $R_N(\Gamma^{1/2}B)\neq 0$ and $F:=\{f_v:v\in\bR^p,\norm{\Gamma^{1/2}v}_2=1\}$ when $R_N(\Gamma^{1/2}B)= 0$. We denote by $\mu=\cN(0,\Gamma)$, the probability measure of the i.i.d. $p$-dimensional rows vectors $Y_1,\cdots, Y_N$ of $\bX_1$. According to Theorem~5.5 in \cite{dirksen2015tail},  there is an absolute constant $C\geq1$ such that for all $t\geq 1$, with probability at least $1-\exp(-t)$,
\begin{equation}\label{eq:sjoerd}
\begin{aligned}
 \sup_{f\in F}\left|\frac{1}{N}\sum_{i=1}^N f^2(Y_i) - \bE f^2(Y_1)\right|\leq C&\bigg( \frac{\diam(F, L_2(\mu))\ell^*(F)}{\sqrt{N}} + \frac{\ell^*(F)^2}{N}+ \diam(F, L_2(\mu))^2\left(\sqrt{\frac{t}{N}} + \frac{t}{N}\right)\bigg)
\end{aligned}
 \end{equation} where $\diam(F,L_2(\mu)):=\sup(\norm{f}_{L_2(\mu)}:f\in F)$. Note that we used the majorizing measure theorem \cite{MR3184689} to get the equivalence between Talagrand's $\gamma_2$-functional and the Gaussian mean width and that for the Gaussian measure $\mu$, the Orlicz space $L_{\psi_2}(\mu)$ is equivalent to the Hilbert space $L_2(\mu)$.  One can check that $\diam(F,L_2(\mu))=1$ and the result follows for $t=N/(64C^2)$,  $\kappa_{RIP} =  1/[8\sqrt{C}]$ and the definition of $R_N(\Gamma^{1/2}B)$ in \eqref{eq:fixed_point_without_rho}.
 {{\mbox{}\nolinebreak\hfill\rule{2mm}{2mm}\par\medbreak}}

Let us give some insight about Theorem~\ref{theo:rip} in the case that is interesting to us i.e. for $\Gamma$ of rank $k$ and $B=B_2^p$. When $\kappa_{RIP}^2 N\geq  k$, we have for all $R>0$, $\ell^*\left((\Gamma^{1/2}B_2^p)\cap (R B_2^p)\right)\leq \ell^*\left(R B_2^p \cap {\rm range}(\Gamma)\right)\leq R\sqrt{k} \leq \kappa_{RIP} R \sqrt{N}$ (where we used the fact that $R B_2^p \cap {\rm range}(\Gamma)$ is the Euclidean ball of radius $R$ in the $k$-dimensional space ${\rm range}(\Gamma)$ and so its Gaussian mean width is less than $R\sqrt{k}$)  and so one has $R_N(\Gamma^{1/2}B_2^p)=0$ which means that $\bX_1$ satisfies an isomorphic property (see \eqref{eq:full_isomorphy}) over the entire space ${\rm range}(\Gamma)$. We will use this property to study $\hat\bbeta$ in the case where $k\leq \kappa_{RIP}N$; we therefore state it in the following result.

\begin{Corollary}\label{cor:isomorphy}There are absolute constants $c_0$ and $c_1$ such that the following holds. If $\Gamma$ is a semi-definite $p\times p$ matrix of rank $k$ such that $k\leq \kappa_{RIP}N$, then for $\bX_1 = \bG^{(N\times p)}\Gamma^{1/2}$, with probability at least $1-c_0\exp(-c_1 N)$, for all $v\in {\rm range}(\Gamma)$,  $$(1/\sqrt{2})\norm{\Gamma^{1/2} v}_2\leq\norm{\bX_1 v}_2/\sqrt{N}\leq \sqrt{3/2}\norm{\Gamma^{1/2} v}_2.$$
\end{Corollary}

Corollary~\ref{cor:isomorphy} can be proved using a straightforward epsilon net argument and there is no need for Theorem~\ref{theo:rip} in that case. However the other case '$k>N$' does not follow in general from such an argument and requires a multi-scale analysis as in Theorem~5.5 in \cite{dirksen2015tail} which is the main ingredient in the proof of Theorem~\ref{theo:rip}. 

When $p>  4\kappa_{RIP}^2N$, we necessarily have $R_N(\Gamma^{1/2}B) >0$ because there exists an $R^*>0$ such that $R^*B_2^p \cap {\rm range}(\Gamma) \subset \Gamma^{1/2}B_2^p\cap {\rm range}(\Gamma)$ (take $R^*$ for instance to be the smallest non-zero singular value of $\Gamma$) and so for all $0<R\leq R^*$, we have $\ell^*\left((\Gamma^{1/2}B_2^p \cap {\rm range}(\Gamma))\cap (R B_2^p)\right)=\ell^*(RB_2^p\cap {\rm range}(\Gamma)) \geq R \sqrt{k}/2>c_0R \sqrt{N}$ and so one necessarily has $R_N(\Gamma^{1/2}B)\geq R^*>0$ (we used here that the Gaussian mean width of the euclidean ball of $\bR^k$ is greater than $\sqrt{k}/2$). This is something expected because when $k>N$, $\bX_1$ has a none trivial kernel and therefore $\bX_1$ cannot be an isomorphy over the entire space ${\rm range}(\Gamma)$. However, Theorem~\ref{theo:rip} shows that $\bX_1$ can still act as an isomorphy over a {\it{restricted}} subset of ${\rm range}(\Gamma)$ and, given the homogeneity property of an isomorphic relation, this set is a cone. In the setup of Theorem~\ref{theo:rip} this cone is that endowed by $B_2^p\cap R_N(\Gamma^{1/2}B)\Gamma^{-1/2} S_2^{p-1}$ since one can check that 
\begin{equation}\label{eq:def_cone}
\begin{aligned}
\cC&:={\rm cone}\left(B_2^p\cap R_N(\Gamma^{1/2}B)\Gamma^{-1/2} S_2^{p-1}\right)= \left\{v\in {\rm range}(\Gamma): R_N(\Gamma^{1/2}B_2^p)\norm{v}\leq \norm{\Gamma^{1/2}v}_2\right\}
\end{aligned}
\end{equation}(where we denote ${\rm cone}(T)=\{\lambda x:x\in T, \lambda\geq0\}$ for all $T\subset \bR^p$). We  therefore speak about a 'restricted isomorphy property' (RIP) in reminiscent of the classical RIP \cite{CaT06} used in Compressed sensing.

As we said, Theorem~\ref{theo:rip} is used below only for  $k$-dimensional ellipsoids -- that is when $\Gamma$ is of rank $k$ and $B=B_2^p$. For this type of convex body, there exists a sharp (up to absolute constants) computation of their Gaussian mean width intersected with an unit Euclidean ball: denote by $\sigma_1\geq\cdots\geq\sigma_k$ the $k$ non-zero singular values of $\Gamma$ associated with $k$ eigenvectors $f_1, \ldots, f_k$. It follows from Proposition~2.5.1 from \cite{MR3184689} that there exists some absolute constant $C_0\geq1$ such that
\begin{equation}\label{eq:talagrand_ellipsoid}
\ell^*\left({\Gamma^{1/2}B_2^p \cap R B_2^p}\right) \leq C_0 \sqrt{\sum_{j=1}^k \min(\sigma_j, R^2)}
\end{equation} (this estimate is sharp up to absolute constants) and so
\begin{equation*}
    R_N(\Gamma^{1/2}B_2^p)\leq \inf\left(R:\sum_{j=1}^k \min(\sigma_j, R^2)\leq c_0 R^2 N\right)
\end{equation*}
where $c_0=(\kappa_{RIP}/C_0)^2$. We may now identify three regimes for the geometric parameter  $R_N(\Gamma^{1/2}B_2^p)$:
\begin{itemize}
  \item[A)] when $k\leq  c_0 N$ (which is up to absolute constants, the case studied in Corollary~\ref{cor:isomorphy});
  \item[B)] when  there exists $k_0\in\{1,\cdots, \lfloor c_0N\rfloor\}$  such that  $\sum_{j= k_0}^k \sigma_j \leq  (c_0N-k_0+1)\sigma_{k_0}$ and in that case, we define 
\begin{equation}\label{eq:k_star_star}
k^{**}=\max\left(k_0\in\{1,\cdots, \lfloor c_0N\rfloor\}: \sum_{j= k_0}^k \sigma_j \leq  (c_0N-k_0+1)\sigma_{k_0} \right);
\end{equation}
  \item[C)] when for all $k_0\in\{1,\cdots, \lfloor c_0N\rfloor\}$ we have $\sum_{j= k_0}^k \sigma_j >  (c_0N-k_0+1)\sigma_{k_0}$.
\end{itemize} Depending on these three cases, we find the following upper bound  on the complexity fixed point
\begin{equation}\label{eq:choice_r_rho_ell2}
R_N(\Gamma^{1/2}B_2^p) \leq \inf\left(R>0: \sum_{j=1}^k \min(\sigma_j, R^2) \leq c_0^2 R^2 N\right) = \left\{
\begin{array}{cc}
0 & \mbox{ in case A)}\\
\sqrt{\sigma_{k^{**}}} & \mbox{ in case B)}\\
\sqrt{\frac{{\Tr(\Gamma) }}{c_0 N}} & \mbox{ in case C).} 
\end{array}
\right.
\end{equation}

We note that $k^{**}$ (when it exists) has a geometrical interpretation: if $\sum_{j =  k^{**}}^k \sigma_j \sim  (c_0N-k^{**}+1)\sigma_k^{**} $ and $k^{**}\leq c_0 N/2$ then $k-k^{**}$ is the smallest dimension of an ellipsoid of the form $\Gamma_{k^{**}:k}^{-1/2}B_2^p$ having an Euclidean section of dimension on the order of $N$ since we see that the Dvoretsky dimension of $\Gamma_{k^{**}:k}^{-1/2}B_2^p$ is $d_*(\Gamma_{k^{**}:k}^{-1/2}B_2^p) \sim {\rm Tr}(\Gamma_{k^{**}:k})/\norm{\Gamma_{k^{**}:k}}_{op} =  \sum_{j= k^{**}}^k \sigma_j/\sigma_{k^{**}}$. So $N\sim d_*(\Gamma_{k^{**}:k}^{-1/2}B_2^p)$ which is the regime of existence of Euclidean sections of dimension on the order of $N$ for $\Gamma_{k^{**}:k}^{-1/2}B_2^p$. 

We also note that in \textit{case B)},  ${\rm span}(f_1, \ldots, f_{k^{**}})\subset \cC$ where $\cC$ is the cone defined in \eqref{eq:def_cone} and in \textit{case C)}, ${\rm span}(f_1, \ldots, f_{r})\subset \cC$ where $r$ is the largest integer such that $c_0 N \sigma_r \geq \Tr(\Gamma)$. Hence, a way to understand the cone $\cC$, is to look at it as a cone surrounding a space endowed by the top eigenvectors of $\Gamma$ up to some dimension $k^{**}$ or $r$ depending on case B) or C). On that cone, $\bX_1=\bG^{(N\times p)}\Gamma^{1/2}$ behaves like an isomorphy.

\section{The approach: $\hat\bbeta$ is a sum of a ridge estimator and an overfitting component} 
\label{sec:our_approach_and_techniques_}
Our proof strategy relies on a decomposition of $\hat\bbeta$ that has been observed in several works. It is called the \textit{self-induced regularization phenomenon}  in \cite{BMR}: the minimum $\ell_2$-norm  interpolant estimator $\hat\bbeta$ can  be written as $\hat\bbeta = \hat\bbeta_{1:k} + \hat\bbeta_{k+1:p}$ where  $\hat\bbeta_{1:k}$ is the projection of $\hat\bbeta$ onto the space $V_{1:k}$ spanned by the top $k$ eigenvectors of $\Sigma$ and $\hat\bbeta_{k+1:p}$ is the projection onto the space spanned by the last $p-k$ eigenvectors of $\Sigma$. Both  $\hat\bbeta_{1:k}$ and  $\hat\bbeta_{k+1:p}$ play very different roles: $\hat\bbeta_{1:k}$ is used for estimation whereas $\hat\bbeta_{k+1:p}$ is used for overfitting. In the next result, a key feature of the $\ell_2$-norm is used to write $\hat\bbeta_{1:k}$ in a way which makes explicit its role as a ridge estimator of $\bbeta_{1:k}^*$ and of $\hat\bbeta_{k+1:p}$ as a minimum $\ell_2$-norm estimator of $X_{k+1:p}\bbeta_{k+1:p}^* + \bxi$. The intuition behind this result comes from Eq.(39) in \cite{BMR} and are similar up to the operator $X_{k+1:p}^\top (X_{k+1:p} X_{k+1:p}^\top)^{-1}$ (see \eqref{eq:ridge_estimator} below for a statement equivalent to Eq.(39) in \cite{BMR}).

\begin{Proposition}\label{prop:decomp}
We have $\hat\bbeta = \hat\bbeta_{1:k} + \hat\bbeta_{k+1:p}$ where
\begin{equation}\label{eq:hat_beta_1_argmin}
\hat\bbeta_{1:k}\in\argmin_{\bbeta\in\bR^p}\left(\norm{X_{k+1:p}^\top (X_{k+1:p} X_{k+1:p}^\top)^{-1}(y-X_{1:k}\bbeta)}_2^2 + \norm{\bbeta}_2^2\right)
\end{equation}and 
\begin{equation*}
\hat\bbeta_{k+1:p} = X_{k+1:p}^\top (X_{k+1:p} X_{k+1:p}^\top)^{-1}(y-X_{1:k}\hat\bbeta_{1:k}).
\end{equation*}
\end{Proposition}
\beginproof
Since $\hat\bbeta$ minimizes the $\ell_2$-norm over the set of vectors $\bbeta\in\bR^p$ such that $\bX\bbeta = y$ and since for all $\bbeta$, $\norm{\bbeta}_2^2 = \norm{\bbeta_{1:k}}_2^2 + \norm{\bbeta_{k+1:p}}_2^2$ it is clear that  $\hat\bbeta = \hat\bbeta_{1:k} + \hat\bbeta_{k+1:p}$ where
\begin{align*}
(\hat\bbeta_{1:k}, \hat\bbeta_{k+1:p})&\in\argmin_{(\bbeta_1,\bbeta_2)\in V_{1:k}\times V_{k+1:p}}\left(\norm{\bbeta_{1}}_2^2 + \norm{\bbeta_{2}}_2^2: X_{1:k}\bbeta_1 + X_{k+1:p}\bbeta_2=y\right)\\
&\in\argmin_{(\bbeta_1,\bbeta_2)\in \bR^p\times \bR^p}\left(\norm{\bbeta_{1}}_2^2 + \norm{\bbeta_{2}}_2^2: X_{1:k}\bbeta_1 + X_{k+1:p}\bbeta_2=y\right)
\end{align*}The result follows by  optimizing separately, first in $\bbeta_2$ and then in $\bbeta_1$.
{\mbox{}\nolinebreak\hfill\rule{2mm}{2mm}\par\medbreak}

Proposition~\ref{prop:decomp} is our starting point to the analysis of the prediction properties of $\hat\bbeta$. It is also a key result to understand the roles played by $\hat\bbeta_{1:k}$ and $\hat\bbeta_{k+1:p}$.

\subsection{Risk decomposition} We prove our main theorems below by using the risk decomposition that follows from the estimator decomposition $\hat\bbeta = \hat\bbeta_{1:k} + \hat\bbeta_{k+1:p}$:  
\begin{equation}\label{eq:risk_decomposition}
\norm{\Sigma^{1/2}(\hat\bbeta - \bbeta^*)}_2^2 = \norm{\Sigma_{1:k}^{1/2}(\hat\bbeta_{1:k} - \bbeta^*_{1:k})}_2^2 +  \norm{\Sigma_{k+1:p}^{1/2}(\hat\bbeta_{k+1:p} - \bbeta^*_{k+1:p})}_2^2. 
\end{equation}Estimation results will follow from high probability upper bounds on the two terms in the right-hand side of \eqref{eq:risk_decomposition}. As in \cite{BMR}, the \textit{prediction component} $\hat\bbeta_{1:k}$ is expected to be a good estimator of $\bbeta_{1:k}^*$, i.e. of the $k$ components of $\bbeta^*$ in the basis of the top $k$ eigenvectors of $\Sigma$. These $k$ components are the most important ones to estimate because they are associated with the largest weights in the prediction norm $\norm{\Sigma^{1/2}\cdot}_2$. We will see that $\hat\bbeta_{1:k}$ estimates $\bbeta_{1:k}^*$ as a ridge estimator. On the contrary, the \textit{overfitting component} $\hat\bbeta_{k+1:p}$ is not expected to be a good estimator of anything (and not of $\bbeta_{k+1:p}^*$ in particular). It is here to (over)fit the data (and in particular the noise) and to make $\hat\bbeta$ an interpolant estimator.

\subsection{$\hat\bbeta_{1:k}$: a ridge estimator of $\bbeta_{1:k}^*$}  One of the two components of our analysis based on the \textit{self-induced regularization} property of $\hat\bbeta$ as exhibited in Proposition~\ref{prop:decomp}, is to prove that $\hat\bbeta_{1:k}$ is a good estimator of $\bbeta^*_{1:k}$ and to that end, we rely on the following observation: if $k$ is chosen so that $N\leq \kappa_{DM} d_*(\Sigma_{k+1:p}^{-1/2}B_2^p)$, then, it follows from DM theorem (see Theorem~\ref{theo:DM}), that with large probability, $\norm{X_{k+1:p}^\top (X_{k+1:p} X_{k+1:p}^\top)^{-1}\cdot}_2$ is isomorphic (i.e. equivalent up to absolute constants) to $(\Tr(\Sigma_{k+1:p}))^{-1/2}\norm{\cdot}_2$ and so, according to Proposition~\ref{prop:decomp}, $\hat\bbeta_{1:k}$ is expected to behave like 
\begin{equation}\label{eq:ridge_estimator}
\argmin_{\bbeta_1\in V_{1:k}}\left(\norm{y-X_{1:k}\bbeta_{1}}_2^2 + \Tr(\Sigma_{k+1:p})\norm{\bbeta_{1}}_2^2\right),
\end{equation}which is a \textbf{ridge estimator} with regularization parameter $\Tr(\Sigma_{k+1:p})$ and with  a random design matrix $X_{1:k}$ (\eqref{eq:ridge_estimator} is like Eq.(39) in \cite{BMR} with a regularization parameter coming from the DM theorem).

When $k\lesssim N$, with high probability, the spectrum of $X_{1:k}$ is up to absolute constants given by $\{\sqrt{N\sigma_1}, \ldots, \sqrt{N\sigma_k},0,\cdots, 0\}$ (because, w.h.p. for all $\bbeta_{1}\in V_{1:k}, \norm{X_{1:k}\bbeta_{1}}_2\sim \sqrt{N}\norm{\Sigma_{1:k}^{1/2}\bbeta_{1}}_2$ and for all $\bbeta_2\in V_{k+1:p}, X_{1:k}\bbeta_2=0$). As a consequence the ridge estimator with regularization parameter $\Tr(\Sigma_{k+1:p})$ from \eqref{eq:ridge_estimator}, will be like an OLS only when $\sqrt{N\sigma_k}\gtrsim \Tr(\Sigma_{k+1:p})$. But since, we chose $k$ such that $\sqrt{N\sigma_{k+1}}\lesssim \Tr(\Sigma_{k+1:p})$ (to apply the DM's theorem), $\hat\bbeta_{1:k}$ will be an OLS estimator of $\bbeta_{1:k}^*$ only when $N\sim \Tr(\Sigma_{k+1:p}) / \sigma_{k+1} \sim d_*(\Sigma_{k+1:p}^{-1/2}B_2^p)$ unless there is a big gap between $\sigma_k$ and $\sigma_{k+1}$. Otherwise, in general, we only have $N\leq \kappa_{DM} d_*(\Sigma_{k+1:p}^{-1/2}B_2^p)$ and so (in general) $\hat\bbeta_{1:k}$ behaves like a ridge estimator and the regularization term $\Tr(\Sigma_{k+1:p})\norm{\bbeta_{1:k}}_2^2$ has an impact on the estimation properties of $\bbeta_{1:k}^*$ by $\hat\bbeta_{1:k}$. This differs from the comment from \cite{BMR} because in \cite{BMR}  $k=k^*_b$ so the ridge regularization term in \eqref{eq:ridge_estimator} has almost no effect because the regularization parameter $\Tr(\Sigma_{k+1:p})$ is smaller than the square of the smallest singular value of $X_{1:k}$ restricted to $V_{1:k}$. 

Our proof strategy is therefore, to analyze $\hat\bbeta_{1:k}$ as a ridge estimator with a random anisotropic design with the extra two difficulties (compared to the classical ridge estimator): the operator $X_{k+1:p}^\top (X_{k+1:p} X_{k+1:p}^\top)^{-1}$ appearing in \eqref{eq:hat_beta_1_argmin} and the output $y$ equals $\bX\bbeta^*+\bxi$ and not $X_{1:k}\bbeta_{1:k}^*+\bxi$. We will handle the first difficulty by using the DM theorem which implies that $X_{k+1:p}^\top (X_{k+1:p} X_{k+1:p}^\top)^{-1}$ is an 'isomorphic' operator and  the second difficulty will be handled by looking at $y$ as $y=X_{1:k}\bbeta_{1:k}^* + (X_{k+1:p}\bbeta_{k+1:p}^*+\bxi)$; which means that $X_{k+1:p}\bbeta_{k+1:p}^*+\bxi$ is considered as noise. In fact, in all our analysis $X_{k+1:p}\bbeta_{k+1:p}^*+\bxi$ is considered as a noise even from the viewpoint of $\hat\bbeta_{k+1:p}$. In particular, $\hat\bbeta$ does not aim at estimating $\bbeta_{k+1:p}^*$ well, even via $\hat\bbeta_{k+1:p}$ as we are now explaining thanks to Proposition~\ref{prop:decomp}.

\subsection{$\hat\bbeta_{k+1:p}$: a minimum $\ell_2$-norm interpolant estimator of $X_{k+1:p}\bbeta_{k+1:p}^*+\bxi$} The \textit{overfitting component} $\hat\bbeta_{k+1:p}$ (see \cite{BMR}), even though it appears as an estimator of $\bbeta_{k+1:p}^*$ in the decomposition \eqref{eq:risk_decomposition}, is not expected to be a good estimator of $\bbeta_{k+1:p}^*$. In fact, the remaining space $V_{k+1:p}$ is (automatically) used by  $\hat\bbeta$ to interpolate $X_{k+1:p}\bbeta_{k+1:p}^*+\bxi$; in particular, it interpolates the noise). 

We see from Proposition~\ref{prop:decomp} that if $\hat\bbeta_{1:k}$ was an exact estimator of $\bbeta_1^*$ then $\hat\bbeta_{k+1:p}$ would be equal to $X_{k+1:p}^\top (X_{k+1:p} X_{k+1:p}^\top)^{-1}(X_{k+1:p}\bbeta_{k+1:p}+\bxi)$, which is the minimum $\ell_2$-norm interpolant of $X_{k+1:p}\bbeta_{k+1:p}^*+\bxi$ -- and not just of $X_{k+1:p}\bbeta_{k+1:p}^*$.  Subspace $V_{k+1:p}$ is therefore, the place where overfitting takes place. This overfitting property of $\hat\bbeta_{k+1:p}$ and so of $\hat\bbeta$, has a price in terms of generalization which can be measured by the price to pay for 'bad estimation' of $\bbeta^*_{k+1:p}$ by  $\hat\bbeta_{k+1:p}$ in the second term $\norm{\Sigma_{k+1:p}^{1/2}(\hat\bbeta_{k+1:p} - \bbeta^*_{k+1:p})}_2$ of \eqref{eq:risk_decomposition}. However, this price  is expected to be small because this estimation error  is associated with the smallest weights $(\sigma_j)_{j\geq k+1}$ in the (weighted $\ell_2^p$) prediction norm $\norm{\Sigma^{1/2}\cdot}_2$. This will indeed be the case under an extra assumption on the spectrum of $\Sigma$, that $N \Tr(\Sigma_{k+1:p}^2) = o(\Tr(\Sigma_{k+1:p})^2)$ which essentially says that the spectrum of $\Sigma_{k+1:p}$ needs to be well-spread, i.e. that it cannot be well approximated by an $N$-sparse vector. 

We will therefore, use the estimation properties of $\hat\bbeta_{1:k}$ and the bound 
\begin{equation}\label{eq:decomp_beta_2_star}
\norm{\Sigma_{k+1:p}^{1/2}(\hat\bbeta_{k+1:p} - \bbeta^*_{k+1:p})}_2\leq \norm{\Sigma_{k+1:p}^{1/2}\hat\bbeta_{k+1:p}}_2 + \norm{\Sigma_{k+1:p}^{1/2}\bbeta^*_{k+1:p}}_2
\end{equation}to handle this second term. In particular, it is clear from \eqref{eq:decomp_beta_2_star}, that we will not try to estimate $\bbeta^*_{k+1:p}$ with  $\hat\bbeta_{k+1:p}$.

\subsection{The key role of parameter $k$} Decomposition of the features space $\bR^p = V_{1:k}\oplus^\perp V_{k+1:p}$ is natural when $\bbeta^*$ is estimated with respect to the prediction/weighted norm $\norm{\Sigma^{1/2}\cdot}_2$ and when $\bbeta^*$ has most of its energy in $V_{1:k}$. This features space decomposition is at the heart  of the decomposition in \eqref{eq:risk_decomposition} and is the one used in \cite{TB21}. 

In the approach described above, parameter $k$ is the parameter of a trade-off between the estimation of $\bbeta_{1:k}^*$ (by $\hat\bbeta_{1:k}$) and the lack of estimation of $\bbeta_{k+1:p}^*$ (by $\hat\bbeta_{k+1:p}$) that permits overfitting. Both properties happen simultaneously inside $\hat\bbeta$ and so $k$ needs to be chosen so that the price for the estimation of $\bbeta_1^*$ and the price of overfitting have equal magnitude. 

From this viewpoint, there is a priori, no reason to take $k$ smaller than $N$; in particular, $\norm{\Sigma_{k+1:p}^{1/2}\bbeta_{k+1:p}^*}_2$ will be part of the price of overfitting (since $\bbeta_{k+1:p}^*$ is not estimated by $\hat\bbeta_{k+1:p}$) and so having $k$ large, will be beneficial for this term --  in particular, when the signal is  strong in $V_{k+1:p}$. Hence, we will explore the properties of $\hat\bbeta$ beyond the case $k\lesssim N$. 

This is different from the previous works on the benign overfitting phenomenon of $\hat\bbeta$ \cite{MR4263288,TB21} which do not study this case.  This is the reason why there are two subsections in the next section covering the two cases  '$k\lesssim N$' and '$k\gtrsim N$'. From a stochastic viewpoint, the two regimes are different because in the first case, $X_{1:k}$ behaves like an isomorphism over the entire space $V_{1:k}$ whereas, in the second case, $X_{1:k}$ cannot be an isomorphism  over the entire space $V_{1:k}$ anymore (because it has a none trivial kernel) but it is an isomorphism restricted to a cone as proved in Theorem~\ref{theo:rip}; a geometric property we will use to prove Theorem~\ref{theo:main_2}, that is for the case $k\gtrsim N$.

\section{Main estimation results in the Gaussian case} 
\label{sec:main_results}
In this section, we provide two estimation results of $\bbeta^*$ by $\hat\bbeta$ with respect to $\norm{\Sigma^{1/2}\cdot}_2$ depending on the value of the parameter $k$ driving the space decomposition $\bR^p=V_{1:k}\oplus^\perp V_{k+1:p}$ with respect to the sample size $N$. We start with the case  $k\lesssim N$ and then we state the result in the case $k\gtrsim N$. We show that these two results also hold for any features space decomposition of the form $\bR^p=V_{J}\oplus^\perp V_{J^c}$ in Section~\ref{sub:on_the_choice_of_}. Given that the minimum $\ell_2$-norm interpolant estimator does not depend on any parameter nor on any features space decomposition this shows that it can find the 'best features space decomposition' by itself. In Section~\ref{sub:lower_bound_on_the_prediction_risk_of_}, we obtain a lower bound on the expected prediction risk which matches the upper bound for the choice $J=\{1, \ldots, k^*_b\}$ where $k^*_b$ has been introduced in \cite{MR4263288} (see \eqref{eq:k_star_BLLT}) for some constant $b=4/\kappa_{DM}$. This shows that choices of features space decomposition as well as the introduction of $k^*_b$ from \cite{MR4263288,TB21}, are the right choices to make regarding the minimum $\ell_2$-norm interpolant estimator in the Gaussian linear model.

\subsection{The small dimensional case $k\lesssim N$} 
\label{sub:the_small_dimensional_case_}


\begin{Theorem}\label{theo:main}[\textbf{the $k\lesssim N$ case.}]
There are absolute constants $c_0$, $c_1$, $c_2$ and $c_3$ such that the following holds.
We assume that there exists $k\leq \kappa_{iso} N$ such that $N\leq \kappa_{DM}d_*(\Sigma_{k+1:p}^{-1/2}B_2^p)$, then the following holds for all such $k$'s. 
With probability at least $1-c_0\exp\left(-c_1\left(|J_1| + N \left(\sum_{j\in J_2}\sigma_j\right)/\left(\Tr(\Sigma_{k+1:p})\right)\right)\right)$, 
\begin{equation*}
\norm{\Sigma^{1/2}(\hat\bbeta-\bbeta^*)}_2\leq \square +  \sigma_\xi\frac{c_2\sqrt{N\Tr(\Sigma_{k+1:p}^2)}}{\Tr(\Sigma_{k+1:p})} + \norm{\Sigma_{k+1:p}^{1/2}\bbeta^*_{k+1:p}}_2
\end{equation*}where $J_1, J_2$ and $\Sigma_{1, thres}$ have been defined in Section~\ref{sec:notation} and
\begin{itemize}
  \item[i)] if $\sigma_1 N< \kappa_{DM}\ell_*^2(\Sigma^{1/2}_{k+1:p}B_2^p)$ then  
  \begin{equation}\label{eq:square_1}
  \begin{aligned}
  \square = c_3\max \bigg\{&\sigma_\xi\sqrt{\frac{\Tr(\Sigma_{1:k})}{\Tr(\Sigma_{k+1:p})}}, \sqrt{\frac{N\sigma_1}{\Tr(\Sigma_{k+1:p})}}\norm{\Sigma_{k+1:p}^{1/2}\bbeta_{k+1:p}^*}_{2},\norm{\bbeta_{1:k}^*}_2\sqrt{\frac{\Tr(\Sigma_{k+1:p})}{N}}\bigg\}
  \end{aligned}
\end{equation}
\item[ii)]   if $\sigma_1 N\geq \kappa_{DM}\ell_*^2(\Sigma^{1/2}_{k+1:p}B_2^p)$ then
\begin{equation}\label{eq:square_main_theo_2}
\begin{aligned}
\square = c_3\max \bigg\{&\sigma_\xi \sqrt{\frac{\left|J_1\right|}{N}},  \sigma_\xi\sqrt{\frac{\sum_{j\in J_2}\sigma_j}{\Tr(\Sigma_{k+1:p})}},\norm{\Sigma_{k+1:p}^{1/2}\bbeta_{k+1:p}^*}_{2},  \norm{\Sigma_{1,thres}^{-1/2}\bbeta_{1:k}^*}_2\frac{\Tr(\Sigma_{k+1:p})}{N}\bigg\}.
\end{aligned}
\end{equation}
\end{itemize}
\end{Theorem}

Several comments on Theorem~\ref{theo:main} are in order and we list some of them below.

\subsection{Effective dimension and the Dvoretsky dimension of an ellipsoid} 
We recall that $\Tr(\Sigma_{k+1:p})/(4\sigma_{k+1})\leq d_*(\Sigma_{k+1:p}^{-1/2}B_2^p)\leq \Tr(\Sigma_{k+1:p})/\sigma_{k+1}$ so that the choice of $k$ in Theorem~\ref{theo:main} is to take it so that $N$ is smaller than the effective rank $r_k(\Sigma)$ as in \cite{MR4263288,TB21} (see  \eqref{eq:effective_ranks}). Theorem~\ref{theo:main} holds for all such $k$'s so that the self-induced regularization phenomenon holds true for every space decomposition $\bR^p = V_{1:k}\oplus^\perp V_{k+1:p}$ simultaneously such that this condition holds. One can therefore optimize either the rate or the deviation parameter over this parameter $k$. This is an interesting adaptive property of the minimum $\ell_2$-norm interpolant estimator that will be commented on later.

\subsection{Price of overfitting} The price paid by $\hat\bbeta$ for overfitting on the training set in terms of generalization capacity, is measured by  the error $$\norm{\Sigma_{k+1:p}^{1/2}(\hat\bbeta_{k+1:p} - \bbeta^*_{k+1:p})}_2.$$ Because, unlike not interpolant and reasonable estimators that would try either to estimate $\bbeta_{k+1:p}^*$ or to take values $0$ on $V_{k+1:p}$ (such as a hard thresholding estimator which avoids an unnecessary variance term on $V_{k+1:p}$), $\hat\bbeta_{k+1:p}$ is in fact interpolating the data (if we had $\hat\bbeta_{1:k}=\bbeta_{1:k}^*$ then we would still have $X_{k+1:p}\hat\bbeta_{k+1:p} = X_{k+1:p}\bbeta_{k+1:p}^* + \bxi$). One may find an upper bound on this rate in Proposition~\ref{prop:price_noise_interpolation}, Appendix~\ref{sub:upper_bounds_on_the_price_for_noise_interpolation}. In the final estimation rate obtained in Theorem~\ref{theo:main} only  the sum remains
\begin{equation}\label{eq:price_noise_interpolant}
\sigma_\xi\frac{17\sqrt{N\Tr(\Sigma_{k+1:p}^2)}}{\Tr(\Sigma_{k+1:p})} + \norm{\Sigma_{k+1:p}^{1/2}\bbeta^*_2}_2
\end{equation}which is what we call the price for overfitting from the convergence rate in Theorem~\ref{theo:main}.

\subsection{Persisting regularization in the ridge estimator $\hat\bbeta_{1:k}$} Another way to look at the self-induced regularization property of $\hat\bbeta$ from \cite{BMR} is to look at it as a persisting regularization property of $\hat\bbeta$ because $\hat\bbeta_{1:k}$ is a ridge estimator. Indeed,  $\hat\bbeta$ is the limit of ridge estimators when the regularization parameter tends to zero. Hence, since $\hat\bbeta_{1:k}$ is also a ridge estimator, $\hat\bbeta$, still keeps a part of the space $\bR^p$ over which it performs a ridge regularization.   

In particular, there are two regimes in Theorem~\ref{theo:main} ('$\sigma_1 N\lesssim \Tr(\Sigma_{k+1:p})$' or '$\sigma_1 N\gtrsim \Tr(\Sigma_{k+1:p})$') because there are two regimes for a ridge estimator:  either the regularization parameter $\Tr(\Sigma_{k+1:p})$ in the ridge estimator \eqref{eq:ridge_estimator} is larger than the square of the largest singular value of $X_{1:k}$, i.e. $\sigma_1 N$ (with high probability and up to absolute constants) or it is not. In the first case '$\sigma_1 N\lesssim \Tr(\Sigma_{k+1:p})$', the regularization parameter in \eqref{eq:ridge_estimator} is so large that it is mainly the regularization term '$\Tr(\Sigma_{k+1:p})\norm{\bbeta_{1:k}}_2^2$' which is minimized. So for $\hat\bbeta_{1:k}$ to be a good estimator of $\bbeta_{1:k}^*$, $\bbeta_{1:k}^*$ is required  to be close to zero and that is why we pay a price proportional to $\norm{\bbeta_{1:k}^*}_2$ in that case in Theorem~\ref{theo:main}. In particular, in that case, BO will require that   $\norm{\bbeta_{1:k}^*}_2^2=o\left(\Tr(\Sigma_{k+1:p})/N\right)$. However, this case should be looked at a pathological one since it is a case where the regularization parameter of a ridge estimator is too large (larger than the square of the largest singular value of the design matrix) so that the data fitting term $\norm{y-X_{1:k}\bbeta_{1:k}}_2^2$ does not play an important role in the definition \eqref{eq:ridge_estimator} of the ridge estimator compared with the regularization term. In particular, sources of generalization errors for $\hat\bbeta$ are due to a bad estimation of $\bbeta_{1:k}^*$ (when $\bbeta_{1:k}^*$ is not close to $0$) as well as overfitting. Since our aim is to identify when overfitting is benign, this case adds some extra difficulties which are not at the heart of the purpose and so we look at it as pathological even though it is possible to also obtain convergence rates in that case from Theorem~\ref{theo:main}.

In the general case where the ridge regularization parameter is not too large, i.e. the second case  '$\sigma_1 N\geq \Tr(\Sigma_{k+1:p})$', then the regularization term appears in the rate through the two sets $J_1$ and $J_2$ as well as in the thresholded matrix $\Sigma_{1,thres}^{-1/2}$. This is the interesting case, because it shows that benign overfitting happens when $\sigma_1 N\gtrsim \Tr(\Sigma_{k+1:p})$,  $N\Tr(\Sigma_{k+1:p}^2) = o\left(\Tr^2(\Sigma_{k+1:p})\right)$,
\begin{equation}\label{eq:situations_BO}
\begin{aligned}
|J_1|=o(N),\sum_{j\in J_2}\sigma_j = o\left(\Tr(\Sigma_{k+1:p})\right), \norm{\Sigma_{k+1:p}^{1/2}\bbeta_{k+1:p}^*}_{2} = o\left(1\right) \mbox{ and }   \norm{\Sigma_{1,thres}^{-1/2}\bbeta_{1:k}^*}_2\frac{\Tr(\Sigma_{k+1:p})}{N} = o(1).
\end{aligned}
\end{equation}In particular, a situation where BO happens depends on both the behavior of $\Sigma$ as well as $\bbeta^*$. Comparing with Theorem~1 from \cite{TB21} (recalled in  Theorem~\ref{theo:TB}), we observe that the convergence rate from Theorem~\ref{theo:main} is better because it fully exploits the thresholding effect of the spectrum of $X_{1:k}$ by the ridge regularization; indeed, compared with the rate from Theorem~\ref{theo:TB}, we have
\begin{equation*}
 \sqrt{\frac{\left|J_1\right|}{N}} + \sqrt{\frac{\sum_{j\in J_2}\sigma_j}{\Tr(\Sigma_{k+1:p})}} \leq \sqrt{\frac{k}{N}} \mbox{ and } \norm{\Sigma_{1,thres}^{-1/2}\bbeta_{1:k}^*}_2 \leq \norm{\Sigma_{1:k}^{-1/2}\bbeta_{1:k}^*}_2
\end{equation*}and the other terms are the same in Theorem~\ref{theo:BLLT} and Theorem~\ref{theo:main} so that Theorem~\ref{theo:main} indeed improves the upper bound result from \cite{TB21} for any $k\lesssim N$ (it also improves the deviation rate, see below). We may also check that the two results coincide when $J_2=\emptyset$; which is what happens for the choice of $k=k^*_b$ for $b=4/\kappa_{DM}$. However, even in that case our proof is different from the one in \cite{TB21} since it is based on the self-induced regularization property of $\hat\bbeta$ and not a bias/variance trade-off. It also improves the deviation parameter from constant to exponentially small in $k^*_b$: $1-c_0\exp(-c_1 k^*_b)$ is the deviation estimate expected for an OLS in the Gaussian linear model over $\bR^{k^*_b}$ and that is the one obtained in Theorem~\ref{theo:main} (see Theorem~\ref{theo:main_intro_gauss} for this statement).

\subsection{The benign overfitting phenomenon happens with large probability} The generalization bounds obtained in Theorem~\ref{theo:main} hold with exponentially large probability. This shows that when  $\Sigma$ and $\bbeta^*$ are so that benign overfitting holds then it is very likely to see it happening on data that is, to see the good generalization property of $\hat\bbeta$ on a test set (even though it interpolates the training data). In Theorem~\ref{theo:BLLT} or Theorem~\ref{theo:TB} the rates are multiplied by the deviation parameter $t$ so that they essentially hold  with constant probability -- unless one is willing to increase the rates --  and does not explain why benign overfitting happens very often in practice.

After the results of this paper were completed, we learned that the recent journal version of \cite{TB21} improved the deviation of the results (see Theorem~5 in \cite{TB21}) in their arxiv version from constant deviation to $1-\exp(-c_0 N)$ for some absolute constant $c_0$. However, they are no estimator achieving consistency  under such a large deviation; in particular, BO cannot happen with such a deviation because of the following\textbf{} lower bound result from \cite[Theorem A']{lecue2013learning}:
\begin{Lemma}\label{lemma:lower_minimax}
    Let $\mu$ be the probability distribution of $X$. There exists an absolute constant $c$ for which the following holds. Let $\cF\subset L_2(\mu)$ be a class that is star-shaped around one of its points (i.e., for some $f_0\in\cF$ and for every $f\in\cF$, $\conv\{f_0,f\}\subset\cF$),  $\cY = \{Y^{f^*}= f^*(X)+\xi:\, f^*\in\cF\}$ where $\xi$ is a centered Gaussian random variable with variance $\sigma_\xi^2$ and is independent with $X$. For any estimator $\hat f$ attains accuracy $\eps_N$ with confidence level $\delta_N$ for any target $Y^{f^*}\in\cY$, we have
    \begin{align*}
        \eps_N \geq \min\left\{ c\sigma_\xi^2\frac{\log\left(1/\delta_N\right)}{N},\, \frac{1}{4}\left(\diam\left(\cF, L_2(\mu)\right)\right)^2 \right\}.
    \end{align*}Here, we say $\hat f$ performs with accuracy $\eps_N$ and confidence $1-\delta_N$ with respect to $\cF$, if for any $Y\in\cY$, with probability (w.r.t. the product measure endowed by the joint distribution of $X$ and $Y$) larger than $1-\delta_N$, $\norm{\hat f - Y}_{L_2}^2 \leq \inf_{f\in\cF}\norm{f - Y}_{L_2}^2+\eps_N$.
\end{Lemma}

Let $\cF = \{f(\cdot) = \left<\cdot,\bbeta\right>: \bbeta\in\bR^p\}$, then $\cF$ is star-shaped around $0$ and $\cY = \{ Y^{\bbeta^*}= \left<X, \bbeta^*\right>+\xi: \bbeta^*\in\bR^p \}$, the minimum $\ell_2$-norm interpolant estimator $\hat \bbeta$ satisfies that: for any $Y^{\bbeta^*}\in\cY$, with probability larger than $1-\delta_N$, $\norm{\left<\hat\bbeta,X\right> - Y^{\bbeta^*}}_{L_2}^2 - \inf_{\bbeta\in\bR^p}\norm{\left<\bbeta,X\right>-Y^{\bbeta^*}}_{L_2}^2 =\norm{\Sigma^{1/2}(\hat\bbeta - \bbeta^*)}_2^2 \leq \eps_N$. However, once we assume that there exists an absolute constant $c_0$ such that  $\delta_N = \exp(-c_0 N)$, by Lemma~\ref{lemma:lower_minimax}, we necessarily have $\eps_N\geq cc_0\sigma_\xi^2$, which indicates that $\hat\bbeta$ is not consistent and that there is no benign overfitting.


\subsection{The large dimensional case $k\gtrsim N$} 
\label{sub:the_large_dimensional_case_}


\begin{Theorem}\label{theo:main_2}[\textbf{the $k\gtrsim N$ case.}]
There  are absolute constants $c_0,c_1,c_2,c_3$ and $C_0$ such that the following holds. We assume that there exists $k\in[p]$ such that $N\leq \kappa_{DM}d_*(\Sigma_{k+1:p}^{-1/2}B_2^p)$ and $R_N(\Sigma_{1:k}^{1/2}B_2^p)\leq \ell_*(\Sigma_{k+1:p}^{1/2}B_2^p)\sqrt{\kappa_{DM}/N}$. The following result then holds for all such $k$'s: with probability at least $1-c_0\exp\left(-c_1\left(|J_1| + N \left(\sum_{j\in J_2}\sigma_j\right)/\left(\Tr(\Sigma_{k+1:p})\right)\right)\right)$, 
\begin{equation*}
\norm{\Sigma^{1/2}(\hat\bbeta-\bbeta^*)}_2\leq \square +  \sigma_\xi\frac{c_2\sqrt{N\Tr(\Sigma_{k+1:p}^2)}}{\Tr(\Sigma_{k+1:p})} + c_3\norm{\Sigma_{k+1:p}^{1/2}\bbeta^*_{k+1:p}}_2
\end{equation*}where $J_1$ and $J_2$ are defined in Section~\ref{sec:notation} and
\begin{itemize}
  \item[i)] if $\sigma_1 N< \kappa_{DM}\ell_*^2(\Sigma_{k+1:p}^{1/2}B_2^p)$ and $\Tr(\Sigma_{1:k})\leq N \sigma_1$ then $\square$ is defined in \eqref{eq:square_1}  
\item[ii)]   if $\sigma_1 N\geq \kappa_{DM}\ell_*^2(\Sigma_{k+1:p}^{1/2}B_2^p)$ and $\sum_{j\in J_2}\sigma_j\leq \kappa_{DM}\ell_*^2(\Sigma_{k+1:p}^{1/2}B_2^p)\left(1-|J_1|/N\right)$ then $\square$ is defined in  \eqref{eq:square_main_theo_2}.
\end{itemize}
\end{Theorem}

As mentioned previously, Theorem~\ref{theo:main_2} is the main result if one wants to lower the price of overfitting by considering features space decomposition $\bR^p = V_{1:k}\oplus^\perp V_{k+1:p}$ beyond the case $k\lesssim N$, which was the only case studied in the literature so far to our knowledge. In particular, we can now identify situations where benign overfitting happens thanks to  Theorem~\ref{theo:main_2} but before that, let us comment on the assumptions in Theorem~\ref{theo:main_2}.

\subsection{Assumptions in Theorem~\ref{theo:main_2}} Both assumptions '$N\leq \kappa_{DM}d_*(\Sigma_{k+1:p}^{-1/2}B_2^p)$' and '$R_N(\Sigma_{1:k}^{1/2}B_2^p)\lesssim \ell_*(\Sigma_{k+1:p}^{1/2}B_2^p)\sqrt{\kappa_{DM}/N}$' are geometrical in nature: the first one involved the Dvoretsky dimension of the ellipsoid $\Sigma_{k+1:p}^{-1/2}B_2^p$ and the second one is used to define the cone \eqref{eq:def_cone} onto which $X_{1:k}$ is an isomorphy. Following \eqref{eq:choice_r_rho_ell2},  the latter condition is implied by the slightly stronger condition:
\begin{itemize}
    \item $\kappa_{DM}\ell_*^2(\Sigma_{k+1:p}^{1/2}B_2^p)\geq N \sigma_{k^{**}}$, when
    \begin{equation*}
        k^{**}=\max\left(k_0\in\{1,\cdots, \lfloor c_0N\rfloor\}: \sum_{j= k_0}^k \sigma_j \leq  (c_0N-k_0+1)\sigma_{k_0} \right)
    \end{equation*}exists.
    \item $c_0 \Tr(\Sigma_{k+1:p})\geq \Tr(\Sigma_{1:k})$, when for all  $k_0\in\{1,\cdots, \lfloor c_0N\rfloor\}, \sum_{j= k_0}^k \sigma_j >  (c_0N-k_0+1)\sigma_{k_0} $.
\end{itemize}

In case \textit{ii)} of Theorem~\ref{theo:main_2} which is the most interesting case for us, the extra condition '$\sum_{j\in J_2}\sigma_j\lesssim \kappa_{DM}\ell_*^2(\Sigma_{k+1:p}^{1/2}B_2^p)\left(1-|J_1|/N\right)$' (compared with the case $k\lesssim N$) is a very weak one given that the term $\sqrt{\sum_{j\in J_2}\sigma_j/ \Tr(\Sigma_{k+1:p})}$ appears in $\square$ and is therefore required to tend to $0$ to see the BO happening.

\subsection{Comparison with the previous results} 
The reader may be interested in the benign overfitting phenomenon when $k\gtrsim N$, which seems to be contradictory to the lower bound of Theorem 4 in \cite{MR4263288} (recalled in Theorem~\ref{theo:BLLT}), i.e. for $k^*_b$ recalled in \eqref{eq:k_star_BLLT} and some absolute constants $b, c>1$:
\begin{equation*}
\bE\norm{\Sigma^{1/2}(\hat\bbeta - \bbeta^*)}_2 \geqslant \frac{\sigma_\xi}{c}\left(\sqrt{\frac{k^*_b}{N}}+ \frac{\sqrt{N\Tr(\Sigma_{k_b^*+1:p}^2)}}{\Tr(\Sigma_{k_b^*+1:p})} \right).
\end{equation*}

However, this is not the case. This is because our choice of $k$ is different from theirs. We choose $k$ such that $\Tr(\Sigma_{k+1:p})\gtrsim N\norm{\Sigma_{k+1:p}}_{op}$ and $R_N(\Sigma_{1:k}^{1/2}B_2^p)\lesssim \sqrt{\Tr(\Sigma_{k+1:p})/N}$. It is therefore, (not yet) optimized and not taken equal to $k^*_b$ or $k^{**}+1$ -- it is a free parameter only required to satisfy the conditions of Theorem~\ref{theo:main_2}. In fact, $k^{**}+1$ (where $k^{**}$ is defined in Equation \eqref{eq:choice_r_rho_ell2}) plays a role similar to $k^*_b$. Indeed, by definition, $k^{**} \in \left\{1,\cdots, c_0 N\right\}$ satisfies
\begin{equation*}
\frac{\Tr(\Sigma_{k^{**}+1:p})}{\norm{\Sigma_{k^{**}+1:p}}_{op}} = \frac{\Tr(\Sigma_{k^{**}+1:k})}{\norm{\Sigma_{k^{**}+1:p}}_{op}}+\frac{\Tr(\Sigma_{k+1:p})}{\norm{\Sigma_{k^{**}+1:p}}_{op}}> c_0 N - k^{**} + 1 + \frac{\Tr(\Sigma_{k+1:p})}{\norm{\Sigma_{k^{**}+1:k}}_{op}} >  N
\end{equation*}and so $k^*_b\leq k^{**}+1$ for $b=1$. Therefore, the $k^*_b$ defined in \cite{MR4263288}, plays a similar role to the geometrical parameter $k^{**}+1$, which is also assumed to be smaller than $c_0N+1$ in one case considered in Theorem~\ref{theo:main_2}.

As remarked previously in \cite{MR4263288}, a necessary condition for benign overfitting, is the existence of a dimension $k^*_b\lesssim N$ such that $\Tr(\Sigma_{k^*_b+1:p})\gtrsim N \sigma_{k^*_b+1}$. However, the existence of such a $k^*_b$ does not necessarily force us to take it to be equal to $k$ in the decomposition $\bR^p = V_{1:k}\oplus^\perp V_{k+1:p}$. A priori, one could find a better features space decomposition $\bR^p=V_{1:k}\oplus^\perp V_{k+1:p}$  that in particular lowers the price for overfitting. Since $\norm{\Sigma_{k+1:p}^{1/2}\bbeta_{k+1:p}^*}_2$ is part of this price \eqref{eq:price_noise_interpolant} and it decreases when $k$ increases, one may look for larger $k$'s, and Theorem~\ref{theo:main_2} shows that it is possible to take $k$ larger than $N$ (up to an absolute constant).


\subsection{Benign overfitting in the regime $k\gtrsim N$} It follows from Theorem~\ref{theo:main_2}, that BO happens for choices of $k$ larger than $N$. Such a situation holds when the spectrum of $\Sigma$ is as follows: let $k_0<N < k$, $a<b<c$, $0<\alpha<1$ and
\begin{equation}\label{eq:example_BO_k_large}
 \sigma_j=a, \forall j=1,\ldots, k_0; \quad \sigma_j = \frac{b}{j^\alpha}, \forall j=k_0+1,\ldots, k \mbox{ and } \sigma_j = c, \forall j\geq k+1.
 \end{equation} It follows from Theorem~\ref{theo:main_2} that $\norm{\Sigma^{1/2}(\hat\bbeta - \bbeta^*)}_2$ tends to zero when $N,p\to +\infty$ when $c\sim bN/[k_0^\alpha p]$, $(k/k_0)^\alpha< p/N$, $a> b/k_0^\alpha$, $k_0=o(N)$, $k^{1-\alpha}k_0^\alpha=o(1)$ and $N=o(p-k)$ (which is the well-known 'over-parametrized' regime needed for BO) and
 \begin{equation}\label{eq:condi_BO_signal_large_k}
 \norm{\bbeta_{1:k_0}^*}_2 = o\left(\frac{\sqrt{a}k_0^\alpha}{b}\right), \quad  \norm{\bbeta_{k_0+1:k}^*}_2 = o\left(\sqrt{\frac{k_0^\alpha}{b}}\right) \mbox{ and } \norm{\bbeta_{k+1:p}^*}_2 = o\left(\sqrt{\frac{k_0^\alpha p}{bN}}\right) 
 \end{equation}where we used that $\norm{\Sigma_{1,thres}^{-1/2}\bbeta_{1:k}^*}_2^2  \leqslant \norm{\bbeta_{J_1}^*}_2^2/a + \left(N/\Tr(\Sigma_{k+1:p})\right)\norm{\bbeta_{J_2}^*}_2^2$ when $a> \Tr(\Sigma_{k+1:p})/N$ and $b <\Tr(\Sigma_{k+1:p})/N$ so that $J_1=\{1,\ldots, k_0\}$ and $J_2=\{k_0+1, \ldots, k\}$. Note that under these assumptions, $k^{**}$ (defined in \eqref{eq:k_star_star}) exists and is on the order of $N$ so that one can apply Theorem~\ref{theo:main_2} in that case. In general, we do not want to make any assumption on $\bbeta_{1:k_0}^*$, which is where we expect most of the information on the signal $\bbeta^*$ to lie, hence, for the $k_0$-dimensional vector  $\bbeta_{1:k_0}^*$, we expect $\norm{\bbeta_{1:k_0}^*}_2$ to be on the order of $\sqrt{k_0}$. This will be the case when  $\alpha>1/2$ or when $a>>b^2 k_0^{1/2-\alpha}$. The remarkable point here is that BO happens even for values of  $k$ larger than $N$. Other situations of BO in the case $k\gtrsim N$ may be found with various speeds of decay on the three regimes introduced in \eqref{eq:example_BO_k_large}.

 A priori, given the bound obtained in Theorem~\ref{theo:BLLT}, \ref{theo:TB}, \ref{theo:main} or \ref{theo:main_2}, the benign overfitting phenomenon depends on both $\Sigma$ and $\bbeta^*$. Ideal situations for BO are when for some $k\in[p]$ (not necessarily smaller than $N$), we have 
\begin{itemize}
  \item[(PO1)] $\bbeta^*$ is mostly supported on $V_{1:k}$ so that $\norm{\Sigma_{k+1:p}^{1/2}\bbeta_{k+1:p}^*}_2=o(1)$,
    \item[(PO2)] the spectrum of $\Sigma_{k+1:p}$, denoted by ${\rm spec}(\Sigma_{k+1:p})$, has to be such that its $\ell_2/\ell_1$-ratio is negligible compared to $1/\sqrt{N}$. This type of condition means that ${\rm spec}(\Sigma_{k+1:p})$ is not compressible: it cannot be well approximated by an $N$-sparse vector. In other word, ${\rm spec}(\Sigma_{k+1:p})$ is required to be a well-spread vector. 
  \item[(PE1)] the cardinality of $J_1$  (the set of eigenvalues of $\Sigma$ larger than $\Tr(\Sigma_{k+1:p})/N$) is negligible in front of $N$ and the remainder of the spectrum of $\Sigma_{1:k}$ is such that $\sum_{j\in J_2}\sigma_j=o(\Tr(\Sigma_{k+1:p}))$,
  \item[(PE2)] top eigenvalues of $\Sigma$ are large so that $\sigma_1 N \geq \Tr(\Sigma_{k+1:p})$ and $\norm{\Sigma_{1,thres}^{-1/2}\bbeta_{1:k}^*}_2$ is negligible in front of $N/\Tr(\Sigma_{k+1:p})$.
\end{itemize}

The price for not estimating $\bbeta^*_{k+1:p}$ is part of the price for overfitting as well as the bias term in the estimation of $\bbeta_{1:k}^*$ by the 'ridge' estimator $\hat\bbeta_{1:k}$; the term $\norm{\Sigma_{k+1:p}^{1/2}\bbeta_{k+1:p}^*}_2$ appearing in both components of the risk decomposition \eqref{eq:risk_decomposition}. One way to lower it, is to require for a condition like \textit{(PO1)}. Together with \textit{(PO2)}, they are the two conditions for BO that come from the price for overfitting. The other \textit{(PE)} conditions come from the estimation of $\bbeta_{1:k}^*$ by $\hat\bbeta_{1:k}$: \textit{(PE1)} is a condition on its variance term and \textit{(PE2)} as well as \textit{(PO1)} are the conditions for BO coming from the control of its bias term.

\subsection{On the choice of $k$ and the self-adaptive property} 
\label{sub:on_the_choice_of_}

The choice of $k$ in \cite{TB21} is limited by the constraint $k\lesssim N$ (see Theorem~\ref{theo:TB}) and so the recommendation from \cite{TB21} is to take $k$ for which $\rho_k:=\Tr(\Sigma_{k+1:p})/[N\norm{\Sigma_{k+1:p}}_{op}]$ is on the order of a constant and if such $k$ doesn't exist, one should take the smallest $k$ for which $\rho_k$ is larger than a constant. However, this recommendation does not take into account all the quantities depending on the signal even though $\norm{\Sigma_{k+1:p}^{1/2}\bbeta_{k+1:p}^*}_2$ appears explicitly in the upper bound from  \cite{TB21} (see Theorem~\ref{theo:TB}).  

A consequence of Theorem~\ref{theo:main} and Theorem~\ref{theo:main_2} is that there is no constraint to choose $k\lesssim N$ and so all features space splitting $\bR^p=V_{1:k}\oplus^\perp V_{k+1:p}$ are allowed even for $k\gtrsim N$ as long as the two geometrical conditions of Theorem~\ref{theo:main_2} hold. In particular, one can chose any $k$ satisfying the geometrical assumptions of these theorems and optimize the upper bound, including signal dependent terms such as $\norm{\Sigma_{k+1:p}^{1/2}\bbeta_{k+1:p}^*}_2$. The best choice of $k$ a priori is  making a trade-off between three terms coming from the estimation of $\bbeta_{J_1}^*$ by the OLS part of $\hat\bbeta_{1:k}$, from the estimation of $\bbeta_{J_2}^*$, by the  'over-regularized' part of $\hat\bbeta_{1:k}$ (where $J_1$ and $J_2$ have been introduced in Theorem~\ref{theo:main}) and the none-estimation of $\bbeta_{k+1:p}^*$ by the overfitting component $\hat\bbeta_{k+1:p}$. On top of that, this trade-off is particularly subtle since the ridge regularization parameter of $\hat\bbeta_{1:k}$ is on the order of $\Tr(\Sigma_{k+1:p})$ and therefore depends on $k$ (as well as the spectrum of $\Sigma$). Fortunately $\hat\bbeta$ does this trade-off by itself.


It follows from the analysis of BO from the last  subsection, that cases where BO happens depend on the coordinates of $\bbeta^*$ in the basis  of eigenvectors of $\Sigma$ and that the best cases are obtained when $\bbeta^*$ is sparse in this basis with all its energy supported on the first top $k$ eigenvectors. However, such a configuration may not be a typical situation for real-world data. Fortunately, the decomposition of the features space $\bR^p$ as $V_{1:k}\oplus^\perp V_{k+1:p}$ is arbitrary and can in fact be more adapted to $\bbeta^*$. It is indeed possible to obtain all the results (Proposition~\ref{prop:decomp} and Theorems~\ref{theo:main} and \ref{theo:main_2}) for any decomposition  $\bR^p = V_{J}\oplus^\perp V_{J^c}$ where $J\subset [p]$, $J^c=[p]\backslash J$ and $V_J={\rm span}(f_j, j\in J)$ as in \cite{DBLP:journals/corr/abs-2106-09276}. The key observation here is that Proposition~\ref{prop:decomp} still holds for this decomposition: one can still write $\hat\bbeta$ as a sum $\hat\bbeta = \hat\bbeta_J + \hat\bbeta_{J^c}$ where 
\begin{equation}\label{eq:hat_beta_J_argmin}
\begin{aligned}
\hat\bbeta_J&\in\argmin_{\bbeta \in\bR^p}\left(\norm{X_{J^c}^\top (X_{J^c} X_{J^c}^\top)^{-1}(y-X_{J}\bbeta)}_2^2 + \norm{\bbeta}_2^2\right) \mbox{ and }\\
\hat\bbeta_{J^c} &= X_{J^c}^\top (X_{J^c} X_{J^c}^\top)^{-1}(y-X_{J}\hat\bbeta_{J})
\end{aligned}
\end{equation}where $X_{J} = \bG^{(N\times p)} \Sigma_{J}^{1/2}$ and  $X_{J^c} = \bG^{(N\times p)} \Sigma_{J^c}^{1/2}$ (so that $X_{J} + X_{J^c} = \bX$)  and $\Sigma_{J} = U D_{J} U^\top$  and  $\Sigma_{J^c} = U D_{J^c} U^\top$ where $ D_{J}  = {\rm diag}(\sigma_1 I(j\in J), \ldots, \sigma_p I(j\in J))$ and  $ D_{J^c}  = {\rm diag}(\sigma_{1}I(j\in J^c), \ldots, \sigma_pI(j\in J^c))$. 

As a consequence, Theorem~\ref{theo:main} and Theorem~\ref{theo:main_2} still holds if one replaces the subsets $\{1,\cdots, k\}$ and $\{k+1,\cdots,p\}$  by $J$ and $J^c$ respectively.

\begin{Theorem}\label{theo:main_3}[\textbf{features space decomposition $\bR^p = V_J\oplus^\perp V_{J^c}$}.]
There  are absolute constants $c_0,c_1,c_2,c_3$ and $C_0$ such that the following holds. Let $J\sqcup J^c$ be a partition of $[p]$. We assume  that $N\leq \kappa_{DM}d_*(\Sigma_{J^c}^{-1/2}B_2^p)$ and $R_N(\Sigma_{J}^{1/2}B_2^p)\leq \ell^*(\Sigma_{J^c}^{1/2}B_2^p)\sqrt{\kappa_{DM}/N}$ (note that $R_N(\Sigma_{J}^{1/2}B_2^p)=0$ when $|J|\leq \kappa_{RIP}N$ so that this conditions holds trivially in that case). We define 
\begin{equation*}
J_1:=\left\{j\in J:\, \sigma_j\geq \frac{\kappa_{DM}\ell^2_*(\Sigma_{J^c}^{1/2}B_2)}{N}\right\},\quad J_2 :=J\backslash J_1
\end{equation*}and  $\Sigma_{J,thres}^{-1/2} := U D_{J,thres}^{-1/2} U^\top$ where $U$ is the orthogonal matrix appearing in the SVD of $\Sigma$ and $D_{J,thres}^{-1/2}$ the diagonal matrix
\begin{equation*}
 {\rm diag}\left( \left(\sigma_1 \vee\frac{\kappa_{DM}\ell^2_*(\Sigma_{J^c}^{1/2}B_2)}{N} \right)^{-1/2} I(1\in J), \ldots, \left(\sigma_p \vee\frac{\kappa_{DM}\ell^2_*(\Sigma_{J^c}^{1/2}B_2}{N} \right)^{-1/2}I(p\in J)\right).
\end{equation*}

 The following result then holds for all such space decomposition $\bR^p = V_J\oplus^\perp V_{J^c}$: with probability at least $1-c_0\exp\left(-c_1\left(|J_1| + N \left(\sum_{j\in J_2}\sigma_j\right)/\left(\Tr(\Sigma_{J^c})\right)\right)\right)$, 
\begin{equation*}
\norm{\Sigma^{1/2}(\hat\bbeta-\bbeta^*)}_2\leq \square(J) +  \sigma_\xi\frac{c_2\sqrt{N\Tr(\Sigma_{J^c}^2)}}{\Tr(\Sigma_{J^c})} + c_3\norm{\Sigma_{J^c}^{1/2}\bbeta^*_{J^c}}_2
\end{equation*}where
\begin{itemize}
  \item[i)]  $\square(J)$ is defined by $\square$ in \eqref{eq:square_1} where $\{1:k\}$ (resp. $\{k+1:p\}$) is replaced by $J$ (resp. $J^c$) when $|J|\lesssim N$ and $\norm{\Sigma_J}_{op} N< \kappa_{DM}\ell^2_*(\Sigma_{J^c}^{1/2}B_2)$ or when  $|J|\gtrsim N$, $\norm{\Sigma_J}_{op} N< \kappa_{DM}\ell^2_*(\Sigma_{J^c}^{1/2}B_2)$ and   $\Tr(\Sigma_{J})\leq N \norm{\Sigma_J}_{op}$;  
\item[ii)]  $\square(J)$ is defined by $\square$ in \eqref{eq:square_main_theo_2} where $\{1:k\}$ ($\{k+1:p\}$ resp. ) is replaced by $J$ ($J^c$ resp.)  when $|J|\lesssim N$  and $\norm{\Sigma_J}_{op} N\geq \kappa_{DM}\ell^2_*(\Sigma_{J^c}^{1/2}B_2)$ or $|J|\gtrsim N$, $\norm{\Sigma_J}_{op} N\geq \kappa_{DM}\ell^2_*(\Sigma_{J^c}^{1/2}B_2)$ and $\sum_{j\in J_2}\sigma_j\leq \kappa_{DM}\ell^2_*(\Sigma_{J^c}^{1/2}B_2)\left(1-|J_1|/N\right)$.
\end{itemize}
\end{Theorem}

It follows from Theorem~\ref{theo:main_3} that the range where benign overfitting  happens can be extended to cases where: a)  estimation of $\bbeta^*_J$ by $\hat\bbeta_J$ is good enough (this happens under the same conditions as for instance in \eqref{eq:condi_BO_signal_large_k} except that $\{1,\cdots, k\}$ should be replaced by $J$) and b) the price for overfitting is low, i.e. when $\norm{\Sigma_{J^c}^{1/2}\bbeta^*_{J^c}}_2$ is small and when the $\ell_2/\ell_1$ ratio of the spectrum of $\Sigma_{J^c}$ is smaller than $o(1/\sqrt{N})$.

Since $\hat\bbeta$ is a parameter free estimator, we observe that  the best features space decomposition $\bR^p = V_{J}\oplus^\perp V_{J^c}$ is performed automatically in $\hat\bbeta$. This is a remarkable property of $\hat\bbeta$ since, given the upper bound from Theorem~\ref{theo:main_3}, this best features space decomposition should a priori depends on  both $\Sigma$ and $\bbeta^*$. However, in the next subsection, we will see that in fact the best choice of $J$ is for $\{1, \ldots, k^*_b\}$ for some constant $b$ so that it only depends on $\Sigma$ and not on $\bbeta^*$.

This subsection and the previous ones lead to the following problem: what is the optimal way to decompose the features space? Does it necessarily have to be a direct sum of (two or more)  eigenspaces of $\Sigma$? It may be the case that there is a better way to entangle $\Sigma$ and $\bbeta^*$ in some $(\Sigma, \bbeta^*)$-adapted basis $\cB$ such that on the top $k$ eigenvectors of $\cB$, $\bbeta^*$ and $\Sigma$ simultaneously have most of their energy and that the restriction of  $\hat\bbeta$ to this $k$-dimensional space is a good (OLS or ridge) estimator of the restriction of $\bbeta^*$ to it and the rest of the space is used for overfitting (as long as it has the properties for benign overfitting, i.e. the restricted spectrum is well-spread and the restriction of $\bbeta^*$ has low energy). We answer these questions in the next subsection.

\section{Lower bound on the prediction risk of $\hat\bbeta$, the best features space decomposition and a definition of benign overfitting} 
\label{sub:lower_bound_on_the_prediction_risk_of_}
In this section, we obtain a lower bound on the expectation of the prediction risk of $\hat\bbeta$. This lower bound improves the lower bounds obtained previously in \cite{MR4263288} and \cite{TB21}. It removes some (a posteriori) unnecessary assumptions on the condition number of $X_{k+1:p}$ as well as the smallest singular values of $A_{-j}$ (see Lemma~3 in \cite{TB21}) and, more importantly, it is a lower bound on the prediction risk of $\hat\bbeta$ itself, not a lower bound on a Bayesian prediction risk as in \cite{MR4263288,TB21}. This lower bound shows that the best space decomposition is of the form $\bR^p=V_J\oplus^\perp V_{J^c}$ where $J=\{1,\ldots, k\}$ and for the optimal choice of $k=k^*_b$ as previously announced in \cite{MR4263288}. This answers the question asked above and proves that the two intuitions from \cite{MR4263288} and \cite{TB21} that the best split of $\bR^p$ is $V_J\oplus^\perp V_{J^c}$ for $J=\{1,\ldots, k\}$ and it is best for $k=k^*_b$ and some well chosen constant $b$. In particular, this optimal choice of features space decomposition depends only on $\Sigma$ and not at all on $\bbeta^*$, even though the convergence rate depends on $\bbeta^*$.

\begin{Theorem}\label{theo:lower_bound} There exists  absolute constants $c_0,c_1>0$ such that the following holds. If $N\geq c_0$ and $\Sigma$ is such that $k^*_b<N/4$ for some  $b\geq \max(4/\kappa_{DM}, 24)$ then
\begin{equation*}
\begin{aligned}
\bE \norm{\Sigma^{1/2}(\hat\bbeta - \bbeta^*)}_2^2\geq \frac{c_0}{b^2}\max \bigg\{ &\frac{ \sigma_\xi^2 k^*_b}{N}, \frac{\sigma_\xi^2 N\Tr(\Sigma_{k_b^*+1:p}^2)}{\Tr^2(\Sigma_{k_b^*+1:p})},\\  &\norm{\Sigma_{k^*_b+1:p}^{1/2}\bbeta_{k^*_b+1:p}^*}_{2}^2, \norm{\Sigma_{1:k^*_b}^{-1/2}\bbeta_{1:k^*_b}^*}_2^2\left(\frac{\Tr(\Sigma_{k^*_b+1:p})}{N}\right)^2\bigg\}.
\end{aligned}
\end{equation*}
\end{Theorem}

This result holds for the generalization excess risk of $\hat\bbeta $ for any  $\bbeta^*$. In particular, Theorem~\ref{theo:lower_bound} differs from the previous lower bound results from \cite{MR4263288,TB21} on some Bayesian risks   where  models on 'the signal' $\bbeta^*$ are assumed.  The upper bound on the generalization risk from Theorem~\ref{theo:main} with $\square$ from \eqref{eq:square_main_theo_2} for $k=k^*_b$, matches (up to an absolute constant) the lower bound from Theorem~\ref{theo:lower_bound}. This indeed shows that optimal features space decomposition is $\bR^p=V_J\oplus^\perp V_{J^c}$ for $J=\{1,\ldots, k^*_b\}$ and that the convergence to zero of the rate obtained in Theorem~\ref{theo:main} for $k=k^*_b$ and Theorem~\ref{theo:lower_bound}, is almost a necessary and sufficient condition for benign overfitting of $\hat\bbeta$ when consistency is defined as the convergence to zero in probability of $\norm{\Sigma^{1/2}(\hat\bbeta - \bbeta^*)}_2$. 

Thanks to Theorem~\ref{theo:main} for $k=k^*_b$ (with $\square$ from \eqref{eq:square_main_theo_2}) and Theorem~\ref{theo:lower_bound} we are in a position to define the following concept of benign overfitting in the linear regression model.

\begin{Definition}\label{def:BO}
We say that \textbf{overfitting is benign for $(\Sigma, \bbeta^*)$} when there exists $k^*_b=o(N)$ such that $\sigma_1 N\geq \Tr(\Sigma_{k^*_b+1:p})$,  $N\Tr(\Sigma_{k^*_b+1:p}^2) = o\left(\Tr^2(\Sigma_{k^*_b+1:p})\right)$,
\begin{equation}\label{eq:def_BO}
 \norm{\Sigma_{k^*_b+1:p}^{1/2}\bbeta_{k^*_b+1:p}^*}_{2} = o\left(1\right) \mbox{ and }  \norm{\Sigma_{1:k^*_b}^{-1/2}\bbeta_{1:k^*_b}^*}_2\frac{\Tr(\Sigma_{k^*_b+1:p})}{N} = o(1).
\end{equation}
\end{Definition}

 The main point in this definition is that it depends on both $\Sigma$ and $\bbeta^*$ unlike the previous one given in \cite{MR4263288} or \cite{DBLP:journals/corr/abs-2106-09276} which depend only on $\Sigma$.  Finally, we also emphasize that, once again, our proof of Theorem~\ref{theo:lower_bound} highlights the role played by the DM theorem on the overfitting part of the features space.


\section{The heavy-tailed case and the universality of the Gaussian case} 
\label{sec:the_heavy_tailed_case}
An inspection of the proof of Theorem~\ref{theo:main} reveals that the scope of the benign overfitting property of the minimum $\ell_2$-norm interpolant estimator can be extended much beyond the Gaussian case. 

As an example, we stated in the Introduction Section, Theorem~\ref{theo:main_intro_weak_moment} showing that overfitting is benign for $\hat\bbeta$  under the solely existence of $\log N$ moments for all marginals of $X$ and $r>4$ moments for the noise. A remarkable feature of this result is that the necessary and sufficient conditions for BO in the Gaussian case as defined in Definition~\ref{def:BO} are the same (up to absolute constants) as the one we obtain in this result under weak moments assumptions. This proves universality of the BO phenomenon for  $\hat\bbeta$  and that the Gaussian case is both the worst case (thanks to the lower bound result from Theorem~\ref{theo:lower_bound} matching the upper bound from Theorem~\ref{theo:main_intro_weak_moment} and the more general result below) and the typical case, even in the non-asymptotic regime. We emphasize that only the $\log N$-moment of the marginals of $X$ is needed, which is even weaker than the classical $\log p$-moment assumption in robust statistics, see \cite{lecue_sparse_2017,chinot_robust_2021}.

Our strategy to extend the result from Theorem~\ref{theo:main}  (obtained in the Gaussian case) to heavy-tailed scenarii is to identify all the properties we used on the design matrix and the noise to prove this result and to check if they still hold under weaker moments assumptions.  Hence, the first thing to do is to identify these properties. 

\subsection{Key properties of the design matrix and the noise} 
\label{sub:key_properties_of_the_design_matrix_and_the_noise_}
Regarding the design matrix $\bX$, the properties we used to prove Theorem~\ref{theo:main} are its geometrical properties as stated in Section~\ref{sec:geometrical_properties_of_gaussian_matrices} that is (with large probability and for some absolute constants $c, c'$ and $c''$): 
\begin{enumerate}
    \item For all $\vlambda\in\bR^N$, $c\sqrt{\Tr(\Sigma_{k+1:p})}\norm{\vlambda}_2 \leq \norm{\bX_{k+1:p}^\top\vlambda}_2\leq c'\sqrt{\Tr(\Sigma_{k+1:p})}\norm{\vlambda}_2$, when $N$ is smaller than the DM dimension of the ellipsoid $\Sigma_{k+1:p}^{-1/2}B_2$.
    \item For all $\vlambda\in\bR^N$, $\norm{\Sigma_{k+1:p}^{1/2}\bX_{k+1:p}^\top\vlambda}_2 \leq c^\prime\left(\sqrt{\Tr(\Sigma_{k+1:p}^2)} + \sqrt{N}\norm{\Sigma_{k+1:p}}_{op} \right)\norm{\vlambda}_2$.
    \item For all $\vv\in\bR^k$, $(1/2)\norm{\Sigma_{1:k}^{1/2}\vv}_2^2 \leq (1/N)\norm{\bX_{1:k}\vv}_2^2  \leq (3/2)\norm{\Sigma_{1:k}^{1/2}\vv}_2^2$, when $N\gtrsim k$
\end{enumerate}as well as 
\begin{itemize}
    \item[4.] for $D = X_{k+1:p}^\top A$, $\Tr\left(DD^\top\right)\leq \frac{c''N\Tr\left(\Sigma_{k+1:p}^2\right)}{\left(\Tr\left(\Sigma_{k+1:p}\right)\right)^2}$.
\end{itemize}
Regarding the noise $\bxi$, the only property we use is that if $D$ is a $p\times N$ matrix then with large probability
\begin{itemize}
    \item[5.] $\norm{D \bxi}_2\leq c_0\sigma_\xi\sqrt{\Tr(DD^\top)}$. 
\end{itemize}

Next, we need to prove that these five conditions on $\bX$ and $\bxi$ hold under weak moments assumptions on the design vector $X$ and the noise. We first identify sufficient conditions where they hold and prove that these sufficient conditions hold under weak moment assumptions.

\subsection{Sufficient conditions for properties \textit{1.} to \textit{5.}} 
\label{sub:sufficient_conditions_for_properties_it}
 We recall that $X_1,\ldots, X_N$ are $N$ i.i.d. copies of $X$ and that $P_{1:k}$ ($P_{k+1:p}$ resp.) is the projection operator onto $V_{1:k}$ ($V_{k+1:p}$ resp.). We consider the following assumptions on the design vector $X$:
\begin{enumerate}
    \item[a)] $\bP\left\{\max_{1\leq i\leq N}\left|\frac{\norm{P_{k+1:p}X_i}_2^2}{\bE \norm{P_{k+1:p}X_i}_2^2}-1\right|\leq \delta\right\}\geq 1-\gamma$ for some constants $\delta<1/2$ and $\gamma<1/2$;
    \item[b)] there are constants $R,L,\alpha\leq 2$ such that for all $\vv\in V_{k+1:p}$ and $2\leq q\leq R\log(e N)$, $\norm{\left<P_{k+1:p}X,\vv\right>}_{L_q}\leq L q^{1/\alpha}\norm{\left<P_{k+1:p}X,\vv\right>}_{L_2}$;  
    \item[c)] there exists $q>4$ and $B\geq1$ such that $\norm{\left<P_{1:k}X,\vv\right>}_{L_q} \leq B^{1/q}\norm{\left<P_{1:k}X,\vv\right>}_{L_2}$, for all $\vv\in\bR^k$;
    \item[d)] there are constants $c_0$ and $\gamma'<1$ such that with probability at least $1-\gamma'$, for all $i\in[N]$, $\norm{\Sigma_{1:k}^{-1/2}X_i}_2\leq c_0 \sqrt{k}$.
\end{enumerate}

In the next result, we show that point~\textit{1.} and point~\textit{2.} hold under \textit{a)} and \textit{b)}. This result follow from an extension of the Dvoretsky-Milman theorem to anisotropic and heavy-tailed vectors. This result of independent interest may be found in Theorem~\ref{theo:DM_weak_moments}.

\begin{Proposition}
Let $X$ be a random vector in $\bR^{p}$ and let $X_1, \ldots, X_N$ be i.i.d. copies of $X$. We assume that the projection of $X$, denoted by $P_{k+1:p}X$, on the space $V_{k+1:p}$ spanned by the $p-k$ smallest singular vectors of the covariance matrix $\Sigma$ of $X$ is symmetric and satisfies \textit{a)} and \textit{b)}. Denote by $d^*\left(\Sigma_{k+1:p}^{-1/2}B_2\right)$ the Dvoretsky-Milman dimension of the ellipsoid $\Sigma^{-1/2}_{k+1:p}B_2$. There are absolute constants $c_1,c_2,c_3$ and $c_4$ such that the following holds. If $N\leq \kappa_{DM}d^*\left(\Sigma_{k+1:p}^{-1/2}B_2\right)$ for a small enough constant $\kappa_{DM}<1$ then with probability at least $1-2\gamma-N^{-2c_1}$ for every $\vlambda\in \bR^{N}$,
\begin{eqnarray*}
    (1-2\delta)\sqrt{\Tr(\Sigma_{k+1:p})}\norm{\vlambda}_2 &\leq& \norm{\bX_{k+1:p}^\top \vlambda}_2 \leq (1+2\delta)\sqrt{\Tr(\Sigma_{k+1:p})}\norm{\vlambda}_2,
\end{eqnarray*} and
\begin{equation*}
    \norm{\Sigma_{k+1:p}^{1/2}\bX_{k+1:p}^\top\vlambda}_2\leq c_3\left( \sqrt{\Tr(\Sigma_{k+1:p}^2)} +  \sqrt{N}\norm{\Sigma_{k+1:p}}_{op} \right)\norm{\vlambda}_2.
\end{equation*}
\end{Proposition}

We move to point~\textit{3.}

\begin{Proposition}\label{prop:stochastic_argument_weak_moments} 
Let $X$ be a random vector in $\bR^{p}$ and let $X_1, \ldots, X_N$ be i.i.d. copies of $X$. We assume that the projection of $X$, denoted by $P_{1:k}X$, on the space $V_{1:k}$ spanned by the $k$ largest singular vectors of the covariance matrix $\Sigma$ of $X$ satisfies \textit{c)} and \textit{d)} and that $C k\leq  N$ where $C$ is some absolute constant depending only on $B,q$ and $c_0$ large enough, then with probability at least $1-\gamma'-1/k$, for all $\vv\in\bR^k$,
    \begin{equation*}
        \frac{1}{2}\norm{\Sigma_{1:k}^{1/2}\vv}_2^2 \leq \frac{1}{N}\norm{\bX_{1:k}\vv}_2^2  \leq \frac{3}{2}\norm{\Sigma_{1:k}^{1/2}\vv}_2^2.
    \end{equation*}
\end{Proposition}

In the next result we prove that point~\textit{4.} hold under a weak moment assumption.
\label{prop:condition_4_proposition_weak_moment}
\begin{Proposition}
  We assume that $X = \Sigma^{1/2}Z$ where $Z$ is a random vector with independent and variance one coordinates $(Z_j)_{j}$, and $X$ satisfies the $L_{R\log{eN}}-L_2$ equivalence assumption in c).
    then with probability at least $1-c_2N^{-1/2}$ with $c_2$ borrowed from Lemma~\ref{lem:moment-sum-variid-positive}, point~\textit{4.} holds true.
\end{Proposition}

Next we move to the condition on the noise we used to prove Theorem~\ref{theo:weak_moment}. We will prove it under the following assumption:
\begin{itemize}
    \item[e)]  The noise vector $\bxi = (\xi_i)_{i=1}^N$ is a random vector with independent mean zero and variance $\sigma_\xi$ coordinates such that for all $i$'s,  $\norm{\xi_i}_{L_r}\leq \kappa \sigma_\xi$ for some $\kappa>0$ and $r>4$.
\end{itemize}

\begin{Proposition}\label{prop:latala_diagonal_noise}
    We assume that the noise $\bxi$ satisfies \textit{e)}.  Then, there exists some absolute constant $C_\kappa$ (depending only on $\kappa$) such that for any matrix $D \in\bR^{p\times N}$ the following holds: if for some integer $k$,  $\sqrt{k} \norm{D}_{op}\leq \sqrt{\Tr(DD^\top)}$ then with probability at least $1-(c_0/k)^{r/4}$,
    \begin{equation*}
        \norm{D\bxi}_2 \leq  (3/2) \sigma_\xi \sqrt{\Tr(DD^\top)}.
    \end{equation*}
\end{Proposition}
The proof of Proposition~\ref{prop:latala_diagonal_noise} is given in Section~\ref{sub:property_on_the_noise}.



\subsection{Main result under weak moments assumptions}
\label{sub:main_result_under_weak_moments_assumptions}
Even though, a result similar to the one of Theorem~\ref{theo:main} for a general index $k\lesssim  N$ holds, we state the next result only for $k=k^*_b$ because it makes the result simpler to read and because we know that it is the best choice of parameter $k$. 

\begin{Theorem}\label{theo:weak_moment}
 We assume that the design vector $X$ is such that $P_{k^*_b+1:p}X$ is symmetric and satisfies \textit{a)} (with parameters $\delta$ and $\gamma$), $P_{1:k^*_b}X$ satisfies \textit{d)} (with parameters $c_0$ and $\gamma^\prime$) and that there exists some $\alpha\geq 1/2$, $R>0$ and $L>0$ such that for all $2\leq q\leq R \log(N)$ and all $\bv\in\bR^p$, $\norm{\inr{X,\bv}}_{L_q}\leq L q^{1/\alpha}\norm{\inr{X,\bv}}_{L_2}$. We assume that the noise is such that \textit{e)} holds for some parameter $\kappa>0$ and $r>4$.

 Then, there are  constants $\kappa_{iso}$, $\kappa_{DM}$, $c_1$, $c_2$ and $c_3$ (depending only on $c_0, \delta<1/4,\kappa, L$ and $R$) such that the following holds. If $k^*_b\leq \kappa_{iso} N$ for $b = 4/\kappa_{DM} $, then the following holds: with probability at least
 \begin{equation}\label{eq:deviation_heavy_tail}
  1-\left(\frac{c_2}{k^*_b}\right)^{r/4} - \left(\frac{c_3 \log N}{N}\right)^{\frac{\log(2p)}{2}}-2\gamma - \frac{1}{N^{c_4}} - \gamma^\prime - \frac{1}{k^*_b},
  \end{equation} 
\begin{align*}
\norm{\Sigma^{1/2}(\hat\bbeta-\bbeta^*)}_2\leq c_1 r^*
\end{align*}where $r^*$ is the rate obtained in the Gaussian case and defined in \eqref{eq:optimal_rate_4_terms}.
\end{Theorem}

The proof of Theorem~\ref{theo:weak_moment} may be found in Appendix~\ref{sec:extension_to_the_heavy_tailed_case}. It is a more general result than the one announced in the Introductory section. The key insight of this result is that the rate obtained in Theorem~\ref{theo:weak_moment} for $\hat\bbeta$ is the same as the one obtained in the Gaussian case and that the conditions under which it holds are also identical to the Gaussian case. In particular, there is no necessity for the introduction of an extra statistical complexity. The price paid by the heavy-tailed case is only in the probability deviation parameter \eqref{eq:deviation_heavy_tail}; we recall that in the Gaussian case the result holds with probability $1-c_0\exp(-c_1k^*_b)$.

\subsection{Proof of Theorem~\ref{theo:main_intro_weak_moment}. An example where conditions \textit{a)} and \textit{d)} hold under weak moment assumptions} 
\label{sub:an_example_where_conditions_it}The final ingredient we need to show that BO holds under weak moment assumptions is to show that \textit{a)} and \textit{d)} hold under a weak moment assumption. It is the aim of the next proposition to provide such a result. Its proof may be found in Appendix~\ref{sec:extension_to_the_heavy_tailed_case}.

\begin{Proposition}\label{prop:weak_moment_assum_a_d}
Let $X$ be a random vector in $\bR^p$ and denote by $X_1, \ldots, X_N$ i.i.d. copies of $X$. We assume that there exists $Z$ a random vector in $\bR^p$ with independent coordinates $z_1, \ldots, z_p$ so that $X=\Sigma^{1/2}Z$. We also assume that there exists $\kappa>0$ and $r>4$ so that for all $j\in[p]$,  $\norm{z_j}_{L_r}\leq \kappa$. Then, there exists an absolute constant $c_1$ depending only on $\kappa$ such that assumptions \textit{a)} and \textit{d)} from Theorem~\ref{theo:weak_moment} hold with $\delta=1/4$, $c_0=2$ and $\gamma = \gamma' = (c_1/N)^{(r-4)/4}$. 
\end{Proposition}

Together with Theorem~\ref{theo:weak_moment}, Proposition~\ref{prop:weak_moment_assum_a_d} shows that Theorem~\ref{theo:main_intro_weak_moment} holds.

\section{Conclusion} 
\label{sec:conclusion}
Our main results and their proofs provides the following ideas on the benign overfitting phenomenon:
\begin{itemize}
  \item It was known \cite{TB21,BMR} that the minimum $\ell_2$-norm interpolant estimator $\hat\bbeta$ is the sum of an estimator $\hat\bbeta_{1:k}$ and an overfitting component  $\hat\bbeta_{k+1:p}$. It is indeed the case that for the optimal features space decomposition $\bR^p=V_J\oplus^\perp V_{J^c}$ for $J=\{1,\ldots, k^*_b\}$,  $\hat\bbeta_{1:k^*_b}$ is a ridge estimator with a regularization parameter $\Tr(\Sigma_{k^*_b+1:p})$ negligible compared to the square of the smallest singular value of the design matrix so it is essentially an OLS (as announced in \cite{BMR}) and $V_{k+1:p}$ is used for interpolation. For other choices of $J$, it is in general, a ridge estimator over $V_J$ (we showed that all the results can be extended to a more general features space decomposition $\bR^p=V_{J}\oplus^\perp V_{J^c}$ even when $|J|\gtrsim N$). Our findings support the previous ideas of 'over-parametrization' (i.e. $p\gg N$) which is necessary for the existence of a space $V_{k+1:p}$ used for interpolation as well as the one that the features space contains a space $V_{1:k^{*}_b}$ with small complexity (because $k^{*}_b\lesssim N$) where  most of the estimation of $\bbeta^*$  happens. 
  \item the proofs of the main results (Theorem~\ref{theo:main}, \ref{theo:main_2} and \ref{theo:main_3}) open a new technical path for the study of interpolant estimators (with the associated technical tools such as the DM theorem) that may go beyond the minimum $\ell_2$-norm interpolant estimator. Indeed, the proof follows from the risk decomposition \eqref{eq:risk_decomposition}, which itself follows from the idea that $\hat\bbeta$ can be decomposed as the sum of an estimator and an overfitting component \cite{BMR}. In particular, our proof does not follow from a bias/variance analysis as in \cite{MR4263288,TB21} because we wanted to put forward the idea that for interpolant estimators, there is part of the features space $\bR^p$ which is not used for estimation purpose but for interpolation. This feature may probably be a specificity to all interpolant estimators and may require some special tools and techniques as the one we developed in this work.
  \item The two geometrical tools we used (i.e. Dvoretsky-Milman theorem and isomorphic and restricted isomorphic properties) are classical results from the local geometry of Banach spaces which are both using Gaussian mean widths as a complexity measure. 

  The study of the benign overfitting phenomenon in the linear regression model may not require new complexity measure tools but just a better understanding of the features space decomposition. It should also come  with a localization argument (which is used here to remove all non plausible candidates in a model).  This contrasts with the discussion from \cite{DBLP:conf/iclr/ZhangBHRV17}. 

  While isomorphic and restricted isomorphic properties have been often used in statistical learning theory (and we used it on $V_{1:k}$, the part of the features space where 'estimation happens'), we do not know of any example where the Dvoretsky-Milman theorem has been used in statistical learning theory to obtain convergence rates. It may be because  it is used to describe the behavior of the minimum $\ell_2$-norm interpolant estimator regarding its overfitting component $\hat\bbeta_{k+1:p}$ on $V_{k+1:p}$.  This part of $\hat\bbeta$ is not an estimator and therefore cannot be analyzed with classical  estimation tools (see \eqref{eq:decomp_beta_2_star}). The DM theorem may therefore be a specific tool for the study of interpolant estimators.

  The Dvoretsky-Milman theorem may be a missing tool in the 'classical theory' that may help to understand interpolant estimators and the benign overfitting phenomenon  (see \cite{DBLP:conf/iclr/ZhangBHRV17,DBLP:journals/corr/abs-2105-14368,BMR} and references therein for discussions regarding the 'failure' of classical statistical learning theory constructed during the 1990's and 2000's). 

  Both tools are however, dealing only with properties of the design matrix $\bX$ and not of the output or the signal. We believe that it is because the linear model that we are considering should be seen as a construction coming for instance, after a linear approximation of more complex models (such as the NTK approximation of neural networks in some cases). It may be also the case that for this construction of a features space (that we could look at a feature learning step coming before our linear model), one may require a complexity measure depending on both $\bX$ and $y$ (or on $\Sigma$ and $\bbeta^*$).  Somehow, the linear model considered here comes after this features learning step and does not require any other tools than that considered in this work.
   \item Even though the optimal features space decomposition, is totally independent of $\bbeta^*$, the benign overfitting phenomenon depends on both $\Sigma$ (spectrum and eigenvectors) and $\bbeta^*$ as well as their interplay; indeed benign overfitting (see the definition in \eqref{eq:def_BO}) requires that $k^*_b=o(N)$, the spectrum of $\Sigma_{k^*_b+1:p}$ is well spread, the singular values of $\Sigma_{1:k^*_b}$ should be much larger than that of $\Sigma_{k^*_b+1:p}$ and that  $\bbeta^*$ and $\Sigma$ need to be well 'aligned': ideally $\bbeta^*$ should be supported on $V_{1:k^*_b}$. 
 \item  the benign overfitting phenomenon happens with  large probability; it happens with probability at least $1-c_0\exp(-c_1k^*_b)$ in the case of $J=\{1, \ldots, k^*_b\}$ for $b=1$ in the optimal case.
 \item The minimum $\ell_2$-norm interpolant estimator $\hat\bbeta$ automatically adapts to the best decomposition of the features space: it 'finds' the best split $\bR^p = V_{1:k^*_b}\oplus^\perp V_{k^*_b+1:p}$ by itself.  However, since this optimal split is independent of the signal, it shows that $\hat\bbeta$ does not learn the best features in $\bR^p$ that could predict the output in the best possible way. In other words, $\hat\bbeta$ does not do any features learning by itself and so it needs an upstream procedure to do it for it; that is a construction of a space $\bR^p$, design matrix $\bX$ with covariance $\Sigma$ and signal $\bbeta^*$ that can well predict the output $Y$ and such that the couple $(\Sigma, \bbeta^*)$ allows for benign overfitting as defined in Definition~\ref{def:BO}.
 \item The new proof strategy we used to show BO that goes through a feature space decomposition $\bR^p=V_{1:k^*_b}\oplus^\perp V_{k^*_b+1:p}$ (following the self-induced regularization property of $\hat\bbeta$) together with the identification of the two geometrical properties of the design matrix $\bX$ adapted to this decomposition, i.e. an isomorphic property on $V_{1:k^*_b}$ and the Dvoretsky-Milman property on $V_{k^*_b+1:p}$ has revealed to be efficient to deal with the heavy-tailed case. Indeed, we were able to use available results for the isomorphic property from \cite{tikhomirov_sample_2018} and to prove a new anisotropic and heavy-tailed version of the probabilistic version of DM theorem (see Appendix~\ref{sec:extension_to_the_heavy_tailed_case} for a proof). Other proof strategies, in particular those relying on the convex Gaussian minmax theorem CGMT, may raise difficult technical problem to prove the BO phenomenon in heavy-tailed scenarii. After the results of the current paper were finished, we learned that the journal version \cite{TB21} extends their result into heavy-tailed scenarii. However, they do not give explicit results on estimation error, nor the heavy-tailed noise under $L_r$ for $r>4$ assumption.
 \item BO holds also under weak moment assumptions under the very same conditions from Definition~\ref{def:BO} as in the Gaussian case. This shows the universality of the Gaussian case as well as it extends the scenarii of benign overfitting. Proving this result required to prove a new anisotropic and heavy-tailed version of the Dvoretsky-Milman's theorem which may be of independent interest.
\end{itemize}

\appendix

\section{Proof of Theorem~\ref{theo:main}} 
\label{sec:proofs}

The appendix is organized as follows.
\begin{itemize}
 \item Section~\ref{sec:proofs} contains the proof of Theorem~\ref{theo:main}.
 \item Section~\ref{sec:proof_of_theorem_theo:main_2} contains the proof of Theorem~\ref{theo:main_2}.
 \item Section~\ref{sec:proof_of_theorem_theo:lower_bound} contains the proof of Theorem~\ref{theo:lower_bound}.
\item Section~\ref{sec:extension_to_the_heavy_tailed_case} to Section~\ref{sec:proof_theo_11} are devoted to the heavy-tailed cases. These two sections contain the proof of all the results from Section~6. 
\end{itemize}

In this section, we provide a proof of Theorem~\ref{theo:main} which relies on the prediction risk decomposition from \eqref{eq:risk_decomposition}. To make this scheme analysis described in Section~3 works we need $X_{1:k}$ to behave like an isomorphy onto $V_{1:k}$: for all $\bbeta_1\in V_{1:k}, \norm{X_{1:k}\bbeta_1}_2\sim \sqrt{N}\norm{\Sigma_{1:k}^{1/2}\bbeta_1}_2$; and we need  $\norm{X_{k+1:p}^\top (X_{k+1:p} X_{k+1:p}^\top)^{-1}\cdot}_2$ to be isomorphic to $(\Tr(\Sigma_{k+1:p}))^{-1/2}\norm{\cdot}_2$. These two properties hold on a event that we are now introducing.

\subsection{Stochastic event behind Theorem~\ref{theo:main}} We  denote by $\Omega_0$ the event onto which we have:
\begin{itemize}
  \item for all $\blambda\in\bR^N$, $(1/(2\sqrt{2}))\sqrt{\Tr(\Sigma_{k+1:p})}\norm{\blambda}_2\leq \norm{X_{k+1:p}^\top\blambda}_2\leq (3/2)\sqrt{\Tr(\Sigma_{k+1:p})}\norm{\blambda}_2$
  \item for all $\bbeta_1\in V_{1:k}$, $(1/2)\norm{\Sigma_{1:k}^{1/2}\bbeta_1}_2\leq (1/\sqrt{N})\norm{X_{1:k}\bbeta_1}_2\leq (3/2)\norm{\Sigma_{1:k}^{1/2}\bbeta_1}_2$.
\end{itemize}It follows from Theorem~\ref{theo:DM} and Corollary~\ref{cor:isomorphy} that if $N\leq \kappa_{DM}d_*(\Sigma_{k+1:p}^{-1/2}B_2^p)$ and $k\leq \kappa_{iso} N$ then $\bP[\Omega_0]\geq 1- c_0 \exp(-c_1N)$. We place ourselves on the event $\Omega_0$ up to the end of the proof.

To make the presentation simpler we denote $\bbeta_1^*=\bbeta_{1:k}^*$, $\bbeta_2^*=\bbeta_{k+1:p}^*$, $\hat\bbeta_1=\hat \bbeta_{1:k}$ and  $\hat\bbeta_2=\hat\bbeta_{k+1:p}$.  

\subsection{Estimation properties of the 'ridge estimator' $\hat\bbeta_{1:k}$} 
\label{sub:estimation_properties_of_the_}
Our starting point is \eqref{eq:hat_beta_1_argmin}:
\begin{equation}\label{eq:hat_beta_1_argmin_2}
\hat\bbeta_{1:k}\in\argmin_{\bbeta_1\in V_{1:k}}\left(\norm{A(y-X_{1:k}\bbeta_1)}_2^2 + \norm{\bbeta_1}_2^2\right)
\end{equation}where we set $A=X_{k+1:p}^\top (X_{k+1:p} X_{k+1:p}^\top)^{-1}$ and where we used that the minimum over $\bR^p$ in \eqref{eq:hat_beta_1_argmin} is actually achieved in $V_{1:k}$ (see the proof of Proposition~3). Next, we use a 'quadratic + multiplier + regularization decomposition' of the excess regularized risk associated with the RERM \eqref{eq:hat_beta_1_argmin_2} similar to the one from \cite{MR3782379} but with the difference that the regularization term in \eqref{eq:hat_beta_1_argmin_2} is the square of a norm and not directly a norm. This makes a big difference (otherwise we will need a lower bound on the quantity $\Tr(\Sigma_{k+1:p})$ which will play the role of the regularization term in \eqref{eq:hat_beta_1_argmin_2}, a condition we want to avoid here). We therefore write for all $\bbeta_1\in V_{1:k}$,
\begin{align}
\notag &\cL_{\bbeta_1} = \norm{A(y-X_{1:k}\bbeta_1)}_2^2 + \norm{\bbeta_1}_2^2 - \left(\norm{A(y-X_{1:k}\bbeta_1^*)}_2^2 + \norm{\bbeta_1^*}_2^2\right)\\
& = \norm{AX_{1:k}(\bbeta_1-\bbeta_1^*)}_2^2 + 2 \inr{A(y-X_{1:k}\bbeta_1^*), AX_{1:k}(\bbeta_1^*-\bbeta_1)} + \norm{\bbeta_1}_2^2 - \norm{\bbeta_1^*}_2^2 \label{eq:classil_Q+M+R}\\
& = \norm{(X_{k+1:p}X_{k+1:p}^\top)^{-1/2}X_{1:k}(\bbeta_1-\bbeta_1^*)}_2^2 \notag\\
&+2 \inr{X_{1:k}^\top (X_{k+1:p}X_{k+1:p}^\top)^{-1}(X_{k+1:p}\bbeta_2^*+\bxi) - \bbeta_1^*,\bbeta_1-\bbeta_1^*} + \norm{\bbeta_1-\bbeta_1^*}_2^2 \label{eq:new_Q+M+R}
\end{align}where we used that for all $\blambda \in\bR^N$, $\norm{A\blambda}_2 = \norm{(X_{k+1:p}X_{k+1:p}^\top)^{-1/2}\blambda}_2$, $A^\top A = (X_{k+1:p}X_{k+1:p}^\top)^{-1}$ and $\norm{\bbeta_1}_2^2 - \norm{\bbeta_1^*}_2^2  = \norm{\bbeta_1-\bbeta_1^*}_2^2 -2\inr{\bbeta_1^*, \bbeta_1^*-\bbeta_1}$.

The last equality is a modification on the regularization term '$\norm{\bbeta_1}_2^2 - \norm{\bbeta_1^*}_2^2$' as well as the multiplier term (i.e. the second term in \eqref{eq:classil_Q+M+R}) in the classical 'quadratic + multiplier + regularization decomposition' of the excess regularized risk  written in \eqref{eq:classil_Q+M+R}. As mentionned previously, this modification is key to our analysis (we may look at it as a 'quadratic + multiplier decomposition' of the excess regularization term).

We will use the excess risk decomposition from \eqref{eq:new_Q+M+R} to prove that with high probability 
\begin{equation}\label{eq:aim}
\norm{\Sigma_{1:k}^{1/2}(\hat\bbeta_1 - \bbeta_1^*)}_2\leq \square \mbox{ and } \norm{\hat\bbeta_1 - \bbeta_1^*}_2\leq \triangle
\end{equation}where $\square$ and $\triangle$ are two quantities that we will choose later. In other words, we want to show that $\hat\bbeta_1\in \bbeta_1^* + B$ where $B= \left\{\bbeta\in V_{1:k}: \vertiii{\bbeta}\leq 1\right\}$ and for all $\bbeta\in \bR^p$,
\begin{equation}\label{eq:interpol_norm_1}
\vertiii{\bbeta} := \max\left(\frac{ \norm{\Sigma_{1:k}^{1/2}\bbeta}_2}{\square} , \frac{ \norm{\bbeta}_2}{\triangle}\right).
\end{equation}  To do that it is enough to show that if $\bbeta_1$ is a vector in $V_{1:k}$ such that $\bbeta_1\notin \bbeta_1^*+ B$ then necessarily $\cL_{\bbeta_1}>0$, because, by definition of $\hat\bbeta_1$, we have $\cL_{\hat\bbeta_1}\leq0$. This is what we are doing now and we start with an homogeneity argument similar to the one in \cite{lecue2013learning}.

 Denote by $\partial  B$ the border of $B$ in $V_{1:k}$. Let  $\bbeta_1 \in V_{1:k}$ be such that $\bbeta_1\notin \bbeta_1^*+ B$.  There exists $\bbeta_0\in \partial  B$ and $\theta>1$ such that $\bbeta_1 - \bbeta_1^* = \theta(\bbeta_0 - \bbeta_1^*)$. Using \eqref{eq:new_Q+M+R}, it is clear that $\cL_{\bbeta_1}\geq \theta \cL_{\bbeta_0}$. As a consequence, if we prove that $\cL_{\bbeta_1}>0$ for all $\bbeta_1\in\bbeta_1^*+\partial B$ this will imply that  $\cL_{\bbeta_1}>0$ for all $\bbeta_1\notin\bbeta_1^*+ B$. Hence, we only need to show the positivity of the excess regularized risk $\cL_{\bbeta_1}$ on the border $\bbeta_1^*+\partial B$ and since an element $\bbeta_1$ in the border $\bbeta_1^*+\partial  B$ may have two different behaviors, we introduce two cases:
 \begin{itemize}
   \item[a)] $\norm{\Sigma_{1:k}^{1/2}(\bbeta_1-\bbeta_1^*)}_2=\square$ and $\norm{\bbeta_1-\bbeta_1^*}_2\leq \triangle$;
   \item[b)] $\norm{\Sigma_{1:k}^{1/2}(\bbeta_1-\bbeta_1^*)}_2\leq \square$ and $\norm{\bbeta_1-\bbeta_1^*}_2 =  \triangle$. 
 \end{itemize}All that remains now is to show that in the two cases $\cL_{\bbeta_1}>0$. 

 If we look at \eqref{eq:new_Q+M+R} among the three terms in this equation only the 'multiplier term' $$\cM_{\bbeta_1}:=2 \inr{X_{1:k}^\top (X_{k+1:p}X_{k+1:p}^\top)^{-1}(X_{k+1:p}\bbeta_2^*+\bxi) - \bbeta_1^*,\bbeta_1-\bbeta_1^*}$$ can possibly be negative whereas the two others quadratic term $$\cQ_{\bbeta_1}:=\norm{(X_{k+1:p}X_{k+1:p}^\top)^{-1/2}X_{1:k}(\bbeta_1-\bbeta_1^*)}_2^2$$ and 'regularization term' $\cR_{\bbeta_1}:= \norm{\bbeta_1-\bbeta_1^*}_2^2$ are positive. As a consequence, we will show that $\cL_{\bbeta_1}>0$ because either $\cQ_{\bbeta_1}> |\cM_{\bbeta_1}|$ (this will hold in case \textit{a)}) or $\cR_{\bbeta_1}>|\cM_{\bbeta_1}|$ (this will hold in case \textit{b)}). 

 Let us first control the multiplier term $\cM_{\bbeta_1}$ for $\bbeta_1\in \bbeta_1^*+\partial B$. We have
 \begin{align*}
 &|\cM_{\bbeta_1}|/2=|\inr{X_{1:k}^\top (X_{k+1:p}X_{k+1:p}^\top)^{-1}(X_{k+1:p}\bbeta_2^*+\bxi) - \bbeta_1^*,\bbeta_1-\bbeta_1^*}|\\
 &\leq \sup_{v\in B}|\inr{X_{1:k}^\top (X_{k+1:p}X_{k+1:p}^\top)^{-1}(X_{k+1:p}\bbeta_2^*+\bxi) - \bbeta_1^*,v}|
 \end{align*}where we recall that $B$ is the unit ball of $\vertiii{\cdot}$ intersected with $V_{1:k}$. It is straightforward to check that for all $\bbeta\in V_{1:k}$, $\vertiii{\bbeta}\leq \norm{\tilde{\Sigma}_1^{1/2}\bbeta}_2\leq \sqrt{2}\vertiii{\bbeta}$ where $\tilde{\Sigma}_1^{1/2} = U \tilde{D}_1^{1/2}U^\top$ and
\begin{equation}\label{eq:tilde_D}
 \tilde{D}_1^{1/2} = {\rm diag}\left(\max\left(\frac{\sqrt{\sigma_1}}{\square}, \frac{1}{\triangle}\right), \ldots, \max\left(\frac{\sqrt{\sigma_k}}{\square}, \frac{1}{\triangle}\right), 0, \ldots, 0\right).
 \end{equation} Therefore, $\vertiii{\cdot}$'s dual norm $\vertiii{\cdot}_*$ is also equivalent to $\norm{\tilde{\Sigma}_1^{1/2}\cdot}_2$'s dual norm  which is given by $\norm{\tilde{\Sigma}_1^{-1/2}\cdot}_2$: for all $\bbeta\in V_{1:k}$, $(1/\sqrt{2})\vertiii{\bbeta}_*\leq \norm{\tilde{\Sigma}_1^{-1/2}\bbeta}_2\leq \vertiii{\bbeta}_*$. Hence, we have for all  $\bbeta_1\in \bbeta_1^*+\partial B$,
 \begin{align}\label{eq:multiplier_decomp}
  &|\cM_{\bbeta_1}|\leq 2 \sqrt{2}\norm{\tilde \Sigma_1^{-1/2}\left(X_{1:k}^\top (X_{k+1:p}X_{k+1:p}^\top)^{-1}(X_{k+1:p}\bbeta_2^*+\bxi) - \bbeta_1^*\right)}_2\\
 &\leq 2 \sqrt{2}\bigg(\norm{\tilde \Sigma_1^{-1/2}X_{1:k}^\top (X_{k+1:p}X_{k+1:p}^\top)^{-1}X_{k+1:p}\bbeta_2^*}_2 + \norm{\tilde \Sigma_1^{-1/2}X_{1:k}^\top (X_{k+1:p}X_{k+1:p}^\top)^{-1}\bxi}_2  + \norm{\tilde \Sigma_1^{-1/2}\bbeta_1^*}_2\bigg)\notag.
 \end{align}Next, we handle the first two terms in \eqref{eq:multiplier_decomp} in the next two lemmas.

 \begin{Lemma}\label{lem:first_beta_2_star_term}
Assume that $N\leq \kappa_{DM}d_*(\Sigma_{k+1:p}^{-1/2}B_2^p)$ and $k\leq \kappa_{iso} N$. With probability at least $1-c_0\exp(-c_1N)$,
 \begin{equation*}
 \norm{\tilde \Sigma_1^{-1/2}X_{1:k}^\top (X_{k+1:p}X_{k+1:p}^\top)^{-1}X_{k+1:p}\bbeta_2^*}_2 \leq \frac{18N \sigma(\square, \triangle)}{\Tr(\Sigma_{k+1:p})}\norm{\Sigma_{k+1:p}^{1/2}\bbeta_2^*}_2
 \end{equation*}where 
\begin{equation}\label{eq:def_sigma_square_triangle}
 \sigma(\square, \triangle) :  = \left\{
\begin{array}{cc}
\square & \mbox{ if } \triangle \sqrt{\sigma_1}\geq \square\\
\triangle \sqrt{\sigma_1} & \mbox{ otherwise.}
\end{array}
\right. 
\end{equation}

 \end{Lemma}

\proof It follows from Bernstein's inequality that with probability at least $1-c_0\exp(-c_1N)$,
we have $\norm{X_{k+1:p}\bbeta_2^*}_2\leq \frac{3\sqrt{N}}{2}\norm{\Sigma_{k+1:p}^{1/2}\bbeta_2^*}_2.$

On the event $\Omega_0$, we have $\norm{\tilde \Sigma_{1:k}^{-1/2}X_{1:k}^\top}_{op} = \norm{X_{1:k} \tilde \Sigma_{1:k}^{-1/2}}_{op}\leq (3/2)\sqrt{N}\norm{\Sigma_{1:k}^{1/2}\tilde\Sigma_{1:k}^{-1/2}}_{op}$ because of the isomorphic property of $X_{1:k}$  and $\norm{(X_{k+1:p}X_{k+1:p}^\top)^{-1}}_{op}\leq (3/2)(\Tr(\Sigma_{k+1:p}))^{-1}$
because of the Dvoretsky-Milman's property satisfied by $X_{k+1:p}$ (see Proposition~\ref{prop:DM_ellipsoid}). As a consequence, with probability at least  $1-c_0\exp(-c_1N)$,

\begin{align*}
&\norm{\tilde \Sigma_{1:k}^{-1/2}X_{1:k}^\top (X_{k+1:p}X_{k+1:p}^\top)^{-1}X_{k+1:p}\bbeta_2^*}_2\leq \norm{\tilde \Sigma_{1:k}^{-1/2}X_{1:k}^\top}_{op} \norm{(X_{k+1:p}X_{k+1:p}^\top)^{-1}}_{op}\norm{X_{k+1:p}\bbeta_2^*}_2\\
&\leq \frac{18N}{\Tr(\Sigma_{k+1:p})}\norm{\Sigma_{1:k}^{1/2}\tilde\Sigma_{1:k}^{-1/2}}_{op}\norm{\Sigma_{k+1:p}^{1/2}\bbeta_2^*}_2
=\frac{18N \norm{\Sigma_{k+1:p}^{1/2}\bbeta_2^*}_2}{\Tr(\Sigma_{k+1:p})} \sigma(\square, \triangle)
\end{align*} because $\norm{\Sigma_{1:k}^{1/2}\tilde\Sigma_{1:k}^{-1/2}}_{op}$ is smaller than $\max\left(\max(\square\1\left(\triangle \sqrt{\sigma_j}\geq \square\right), \triangle \sqrt{\sigma_j}\1\left(\triangle \sqrt{\sigma_j}\leq \square\right)):j\in[k]\right),$
which is equal to $\sigma(\square, \triangle)$ as defined in \eqref{eq:def_sigma_square_triangle}.
{{\mbox{}\nolinebreak\hfill\rule{2mm}{2mm}\par\medbreak}}

To control the second term in \eqref{eq:multiplier_decomp}, we define the following two subsets
\begin{equation*}
J_1:=\left\{j\in[k]:\, \sigma_j\geqslant\left(\square/\triangle\right)^2\right\} \mbox{ and } J_2 := [k]\backslash J_1 = \left\{j\in[k]:\, \sigma_j<\left(\square/\triangle\right)^2\right\}.
\end{equation*}
Later in this section, we will observe that our selection of $\square/\triangle$ causes the indices set $J_1$ and $J_2$ to coincide with those defined in Section~\ref{sec:notation}.

 \begin{Lemma}\label{lem:noise_term}Assume that $N\leq \kappa_{DM}d_*(\Sigma_{k+1:p}^{-1/2}B_2^p)$ and $k\leq \kappa_{iso} N$.
 With probability at least $1-\exp(-t(\square, \triangle)/2)-c_0\exp(-c_1N)$,
 \begin{equation*}
\norm{\tilde \Sigma_1^{-1/2}X_{1:k}^\top (X_{k+1:p}X_{k+1:p}^\top)^{-1}\bxi}_2 \leq \frac{32\sqrt{N}\sigma_\xi}{\Tr(\Sigma_{k+1:p})}\sqrt{ \left|J_1 \right|\square^2 + \triangle^2\sum_{j\in J_2}\sigma_j}
 \end{equation*}where $t(\square, \triangle):=\left(\left|J_1 \right|\square^2 + \triangle^2\sum_{j\in J_2}\sigma_j\right)/\sigma^2(\square, \triangle)$ and $\sigma(\square, \triangle)$ has been defined in Lemma~\ref{lem:first_beta_2_star_term}.
 \end{Lemma}
 \proof
It follows from the Borell-TIS's inequality (see Theorem~7.1 in \cite{Led01} or p.56-57 in \cite{MR2814399} -- note that we apply these results to the separable set $B^*$) that for all $t>0$, conditionally on $\bX$, with probability (w.r.t. the randomness of the noise $\bxi$) at least $1-\exp(-t/2)$,
\begin{equation}\label{eq:borel_first_lemma}
\norm{\tilde \Sigma_1^{-1/2}X_{1:k}^\top (X_{k+1:p}X_{k+1:p}^\top)^{-1}\bxi}_2 \leq \sigma_\xi\left(\sqrt{\Tr(DD^\top)} + \sqrt{t}\norm{D}_{op}\right)
\end{equation}where $D = \tilde \Sigma_1^{-1/2}X_{1:k}^\top (X_{k+1:p}X_{k+1:p}^\top)^{-1}$. For the basis $(u_j)_{j=1}^p$ of eigenvectors of $\Sigma$, on the event $\Omega_0$, we have 
\begin{align*}
&\Tr(DD^\top) = \sum_{j=1}^p \norm{D^\top u_j}_2^2 = \sum_{j=1}^k \left(\frac{\sqrt{\sigma_j}}{\square}\vee \frac{1}{\triangle}\right)^{-2}\norm{\left(X_{k+1:p}X_{k+1:p}^\top\right)^{-1}X_{1:k}u_j}_2^2\\
& \leq \sum_{j=1}^k \left(\frac{\sqrt{\sigma_j}}{\square}\vee \frac{1}{\triangle}\right)^{-2}\norm{\left(X_{k+1:p}X_{k+1:p}^\top\right)^{-1}}_{op}^2\norm{X_{1:k}u_j}_2^2\leq \sum_{j=1}^k \left(\frac{\sqrt{\sigma_j}}{\square}\vee \frac{1}{\triangle}\right)^{-2}\frac{256N\sigma_j}{\left(\Tr{\Sigma_{k+1:p}}\right)^2}
\end{align*}and therefore
\begin{equation*}
\sqrt{\Tr(DD^\top)} \leq \frac{16\sqrt{N}}{\Tr(\Sigma_{k+1:p})}\sqrt{ \left|J_1 \right|\square^2 + \triangle^2\sum_{j\in J_2}\sigma_j}.
\end{equation*}Using similar arguments and that $\norm{\Sigma_{1:k}^{1/2}\tilde\Sigma_{1:k}^{-1/2}}_{op}$ is smaller than $\sigma(\square, \triangle)$, we also have (still on the event $\Omega_0$) 
\begin{equation*}
\begin{aligned}
\norm{D}_{op} &= \norm{D^\top}_{op} \leq  \norm{(X_{k+1:p}X_{k+1:p}^\top)^{-1}}_{op} \norm{X_{1:k} \tilde \Sigma_1^{-1/2} }_{op} \leq \frac{16\sqrt{N} \norm{\Sigma_{1:k}^{1/2}\tilde\Sigma_{1:k}^{-1/2}}_{op}}{\Tr(\Sigma_{k+1:p})}\leq  \frac{16\sqrt{N}\sigma(\square, \triangle)}{\Tr(\Sigma_{k+1:p})}.
\end{aligned}
\end{equation*}We use the later upper bounds on $\sqrt{\Tr(DD^\top)}$ and $\norm{D}_{op}$  and take $t=\frac{\left|J_1 \right|\square^2 + \triangle^2\sum_{j\in J_2}\sigma_j}{\sigma^2(\square,\triangle)}$
in \eqref{eq:borel_first_lemma} to conclude.
{{\mbox{}\nolinebreak\hfill\rule{2mm}{2mm}\par\medbreak}}

Let us now show that $\cL_{\bbeta_1}>0$ for $\bbeta_1\in\bbeta_1^*+\partial B$ using the upper bounds on the multiplier part $\cM_{\bbeta_1}$ which follows from the last two lemmas. Let us first analyze \textit{case a)}: let $\bbeta_1\in V_{1:k}$ be such that $\norm{\Sigma^{1/2}(\bbeta_1-\bbeta_1^*)}_2=\square$ and $\norm{\bbeta_1-\bbeta_1^*}_2\leq \triangle$. In that case, we now show that $\cQ_{\bbeta_1}>\cM_{\bbeta_1}$. On the event $\Omega_0$, we have
 \begin{equation}\label{eq:lower_bound_quad_case_a}
 \begin{aligned}
 \cQ_{\bbeta_1} &= \norm{(X_{k+1:p}X_{k+1:p}^\top)^{-1/2}X_{1:k}(\bbeta_1-\bbeta_1^*)}_2^2  \geq \frac{N}{32 \Tr(\Sigma_{k+1:p})}\norm{\Sigma_{1:k}^{1/2}(\bbeta_1-\bbeta_1^*)}_2^2 = \frac{N\square^2}{32\Tr(\Sigma_{k+1:p})}.
 \end{aligned}
 \end{equation}It follows from \eqref{eq:multiplier_decomp} that to get $\cQ_{\bbeta_1}>\cM_{\bbeta_1}$ (and so  $\cL_{\bbeta_1}>0$), it is enough to have
\begin{equation}\label{eq:ineq_quad_multipli}
\begin{aligned}
\frac{N\square^2}{64\sqrt{2} \Tr(\Sigma_{k+1:p})} &> \norm{\tilde \Sigma_1^{-1/2}X_{1:k}^\top (X_{k+1:p}X_{k+1:p}^\top)^{-1}X_{k+1:p}\bbeta_2^*}_2  + \norm{\tilde \Sigma_1^{-1/2}X_{1:k}^\top (X_{k+1:p}X_{k+1:p}^\top)^{-1}\bxi}_2 +
\norm{\tilde \Sigma_1^{-1/2}\bbeta_1^*}_2.
\end{aligned}
\end{equation}

Let us now assume that we choose $\square$ and $\triangle$ so that  $\sigma(\square, \triangle)=\square$, which means that $\triangle\sqrt{\sigma_1}\geq \square$ (we will explore the other case later).  Using Lemma \ref{lem:first_beta_2_star_term} and \ref{lem:noise_term}, inequality \eqref{eq:ineq_quad_multipli} holds if we take $\square$ such that
\begin{equation}\label{eq:square}
\begin{aligned}
\square \geq C_0\max \bigg\{&\sigma_\xi \sqrt{\frac{\left|J_1\right|}{N}},  \left(\triangle\sigma_\xi\sqrt{\frac{1}{N}\sum_{j\in J_2}\sigma_j} \right)^{1/2}, \norm{\Sigma_{k+1:p}^{1/2}\bbeta_2^*}_{2},\sqrt{\norm{\tilde \Sigma_1^{-1/2}\bbeta_1^*}_2\frac{\Tr(\Sigma_{k+1:p})}{N}}\bigg\}
\end{aligned}
\end{equation}where $C_0=4608\sqrt{2}$.

Next, let us analyse \textit{case b)}: let $\bbeta_1\in V_{1:k}$ be such that $\norm{\Sigma^{1/2}(\bbeta_1-\bbeta_1^*)}_2\leq\square$ and $\norm{\bbeta_1-\bbeta_1^*}_2= \triangle$. In that case, we show that $\cL_{\bbeta_1}>0$ by proving that $\cR_{\bbeta_1}>\cM_{\bbeta_1}$. Since $\cR_{\bbeta_1}=\triangle^2$, it follows from \eqref{eq:multiplier_decomp} that $\cR_{\bbeta_1}>\cM_{\bbeta_1}$ holds when
\begin{equation*}
\begin{aligned}
\triangle^2 > 2\sqrt{2} \bigg(&\norm{\tilde \Sigma_1^{-1/2}X_{1:k}^\top (X_{k+1:p}X_{k+1:p}^\top)^{-1}X_{k+1:p}\bbeta_2^*}_2 +\norm{\tilde \Sigma_1^{-1/2}X_{1:k}^\top (X_{k+1:p}X_{k+1:p}^\top)^{-1}\bxi}_2+\norm{\tilde \Sigma_1^{-1/2}\bbeta_1^*}_2\bigg).
\end{aligned}
\end{equation*}Using Lemma \ref{lem:first_beta_2_star_term} and \ref{lem:noise_term}, the latter inequality holds when
\begin{equation}\label{eq:triangle}
\begin{aligned}
\triangle^2 \geq C_1\max \bigg\{&\frac{\sigma_\xi \square \sqrt{|J_1|N}}{\Tr(\Sigma_{k+1:p})},  \frac{\sigma_\xi^2N}{\Tr(\Sigma_{k+1:p})^2}\sum_{j\in J_2}\sigma_j, \frac{N\square}{\Tr(\Sigma_{k+1:p})}\norm{\Sigma_{k+1:p}^{1/2}\bbeta_2^*}_{2},\norm{\tilde\Sigma_1^{-1/2}\bbeta_1^*}_2 \bigg\}
\end{aligned}
\end{equation}for $C_1=4608$.

We set $\triangle$ so that $\square/\triangle = \sqrt{\kappa_{DM}\ell_*^2(\Sigma_{k+1:p}^{1/2}B_2^p)/N}$ and take
\begin{equation}\label{eq:def_square}
\square \geq 4C_0/\kappa_{DM}\max \left\{\sigma_\xi \sqrt{\frac{\left|J_1\right|}{N}},  \sigma_\xi\sqrt{\frac{\sum_{j\in J_2}\sigma_j}{\Tr(\Sigma_{k+1:p})}}, \norm{\Sigma_{k+1:p}^{1/2}\bbeta_2^*}_{2}, \sqrt{\norm{\tilde \Sigma_1^{-1/2}\bbeta_1^*}_2\frac{\Tr(\Sigma_{k+1:p})}{N}}\right\}.
\end{equation}In particular, we check that for this choice, $\triangle$ satisfies \eqref{eq:triangle} and  $\triangle\sqrt{\sigma_1}\geq \square$ holds iff $\sigma_1 \geq \kappa_{DM}\ell_*^2(\Sigma_{k+1:p}^{1/2}B_2^p)/N$. As a consequence, in the case where $\sigma_1 \geq \kappa_{DM}\ell_*^2(\Sigma_{k+1:p}^{1/2}B_2^p)/N$, $\square$ is a valid upper bound on the rate of convergence of $\hat\bbeta_1$ toward $\bbeta_1^*$ with respect to the $\norm{\Sigma_{1:k}^{1/2}\cdot}_2$-norm. This upper bound holds with probabilty at least $1-c_0\exp(-c_1 N) - \exp(-t(\square, \triangle)/2)$ where in that case $t(\square, \triangle) = |J_1| + N \frac{\sum_{j\in J_2}\sigma_j}{\Tr(\Sigma_{k+1:p})}.$ We also note that because $\square/\triangle = \ell_*(\Sigma_{k+1:p}^{1/2}B_2^p)\sqrt{\kappa_{DM}}/\sqrt{N}$, one has $\tilde \Sigma_1^{-1/2} = U \tilde D_1^{-1/2}U^\top$ with $\tilde D_1^{-1/2}$ being the diagonal matrix
\begin{equation*}
\begin{aligned}
 &\square {\rm diag}\left( \left(\sigma_1 \vee\frac{\kappa_{DM}\ell_*^2(\Sigma_{k+1:p}^{1/2}B_2^p)}{N} \right)^{-1/2}, \cdots, \left(\sigma_k \vee\frac{\kappa_{DM}\ell_*^2(\Sigma_{k+1:p}^{1/2}B_2^p)}{N} \right)^{-1/2}, 0, \ldots, 0\right):= \square D_{1,thres}^{-1/2}
\end{aligned}
\end{equation*}and so one can choose $\square$ such that
\begin{equation}\label{eq:def_square_fin_1}
\square = \frac{16C_0^2}{\kappa_{DM}}\max \left\{\sigma_\xi \sqrt{\frac{\left|J_1\right|}{N}},  \sigma_\xi\sqrt{\frac{\sum_{j\in J_2}\sigma_j}{\Tr(\Sigma_{k+1:p})}}, \norm{\Sigma_{k+1:p}^{1/2}\bbeta_2^*}_{2}, \norm{\Sigma_{1,thres}^{-1/2}\bbeta_1^*}_2\frac{\Tr(\Sigma_{k+1:p})}{N}\right\}
\end{equation}where $\Sigma_{1,thres}^{-1/2} = U D_{1,thres}^{-1/2} U^\top$.

Let us now consider the other case that is when  $\square$ and $\triangle$ are chosen so that  $\sigma(\square, \triangle)=\triangle \sqrt{\sigma_{1}}$, which means that $\triangle\sqrt{\sigma_1}\leq \square$. In that case, $J_1$ is empty, $J_2=[k]$, $t(\square, \triangle)  = \Tr(\Sigma_{1})/\sigma_1$, $\tilde \Sigma_1^{1/2}=U \tilde D_1^{1/2}U^\top$ with $\tilde D_1^{1/2}=(1/\triangle){\rm diag}(1,\cdots, 1, 0, \cdots, 0)$ (with $k$ ones and $p-k$ zeros). In \textit{case a)}, it follows from Lemma~\ref{lem:first_beta_2_star_term}, Lemma~\ref{lem:noise_term} and inequality \eqref{eq:ineq_quad_multipli}  that $\cQ_{\bbeta_1}>\cM_{\bbeta_1}$ holds when 
\begin{equation}\label{eq:square_2}
\square^2 \geq C_0^\prime\max \left\{\triangle\sigma_\xi\sqrt{\frac{\Tr(\Sigma_{1:k})}{N}}, \triangle \sqrt{\sigma_1}\norm{\Sigma_{k+1:p}^{1/2}\bbeta_2^*}_{2}, \triangle\norm{\bbeta_1^*}_2\frac{\Tr(\Sigma_{k+1:p})}{N}\right\}
\end{equation}where $C_0^\prime$ is some absolute constant. In \textit{case b)}, we will have $\cR_{\bbeta_1}>\cM_{\bbeta_1}$ when 
\begin{equation}\label{eq:triangle_2}
\triangle^2 \geq C_1^\prime\max \left\{\frac{\sigma_\xi^2N \Tr(\Sigma_{1:k})}{\Tr(\Sigma_{k+1:p})^2}, \sigma_1\left(\frac{N}{\Tr(\Sigma_{k+1:p})}\norm{\Sigma_{k+1:p}^{1/2}\bbeta_2^*}_{2}\right)^2, \norm{\bbeta_1^*}_2^2 \right\}
\end{equation}for $C_1^\prime$ is some absolute constant. We choose $\triangle=\square \sqrt{N/\kappa_{DM}\ell_*^2(\Sigma_{k+1:p}^{1/2}B_2^p)}$ and take
\begin{equation}\label{eq:def_square_2}
\square = C_0^\prime\max \left\{\sigma_\xi\sqrt{\frac{\Tr(\Sigma_{1:k})}{\Tr(\Sigma_{k+1:p})}}, \sqrt{\frac{N\sigma_1}{\Tr(\Sigma_{k+1:p})}}\norm{\Sigma_{k+1:p}^{1/2}\bbeta_2^*}_{2}, \norm{\bbeta_1^*}_2\sqrt{\frac{\Tr(\Sigma_{k+1:p})}{N}}\right\}.
\end{equation}In particular, we check that for this choice, $\triangle$ satisfies \eqref{eq:triangle_2}.  Therefore, when $\triangle\sqrt{\sigma_1}\leq \square$ -- which is equivalent to $\sigma_1 N\leq \kappa_{DM}\ell_*^2(\Sigma_{k+1:p}^{1/2}B_2^p)$ -- the choice of $\square$ from \eqref{eq:def_square_2} provides an upper bound on $\norm{\Sigma_1^{1/2}(\hat\bbeta_1-\bbeta_1^*)}_2$ with probability at least $1-c_0\exp(-c_1N)$ -- note that $t(\square, \triangle)  = \Tr(\Sigma_{1})/\sigma_1\geq N\geq k$ so that $\exp(-t(\square, \triangle)/2)\leq \exp(-k/2)$. 

We are now gathering our findings on the estimation properties of $\hat\bbeta_1$ in the next result. We state the result for both norms $\norm{\Sigma_1^{1/2}\cdot}_2$ and $\norm{\cdot}_2$.

\begin{Proposition}\label{prop:esti_beta_1}There are absolute constants $c_0$, $c_1$ and $c_2$ such that the following holds.
We assume that there exists $k\leq \kappa_{iso} N$ such that $N\leq \kappa_{DM}d_*(\Sigma_{k+1:p}^{-1/2}B_2^p)$, then the following holds for all such $k$'s. 
With probability at least $1-p^*$,
\begin{equation*}
    \norm{\Sigma_1^{1/2}(\hat\bbeta_1-\bbeta_1^*)}_2\leq \square, \quad \norm{\hat\bbeta_1-\bbeta_1^*}_2\leq \square\sqrt{N/\Tr(\Sigma_{k+1:p})},
\end{equation*}
where, 
\begin{itemize}
  \item[i)] if $\sigma_1 N\leq \kappa_{DM}\ell_*^2(\Sigma_{k+1:p}^{1/2}B_2^p)$, $p^* = c_0\exp(-c_1N)$ and 
  \begin{equation*}
\square = C_0^2\max \left\{\sigma_\xi\sqrt{\frac{\Tr(\Sigma_{1:k})}{\Tr(\Sigma_{k+1:p})}}, \sqrt{\frac{N\sigma_1}{\Tr(\Sigma_{k+1:p})}}\norm{\Sigma_{k+1:p}^{1/2}\bbeta_2^*}_{2}, \norm{\bbeta_1^*}_2\sqrt{\frac{\Tr(\Sigma_{k+1:p})}{N}}\right\}
\end{equation*}
\item[ii)]   if $\sigma_1 N\geq \kappa_{DM}\ell_*^2(\Sigma_{k+1:p}^{1/2}B_2^p)$, $p^* = c_0\exp\left(-c_1\left(|J_1| + N \left(\sum_{j\in J_2}\sigma_j\right)/\left(\Tr(\Sigma_{k+1:p})\right)\right)\right)$ and
\begin{equation*}
\square = C_0^2\max \left\{\sigma_\xi \sqrt{\frac{\left|J_1\right|}{N}},  \sigma_\xi\sqrt{\frac{\sum_{j\in J_2}\sigma_j}{\Tr(\Sigma_{k+1:p})}}, \norm{\Sigma_{k+1:p}^{1/2}\bbeta_2^*}_{2}, \norm{\Sigma_{1,thres}^{-1/2}\bbeta_1^*}_2\frac{\Tr(\Sigma_{k+1:p})}{N}\right\}.
\end{equation*}
\end{itemize}
\end{Proposition}

The two cases \textit{i)} and \textit{ii)} appear naturally in the study of the ridge estimator \eqref{eq:ridge_estimator}: in \textit{case i)}, the regularization parameter $\Tr(\Sigma_{k+1:p})$ is larger than the square of the top singular value of $X_{1:k}$ and so the ridge estimator is mainly minimizing the regularization norm. Whereas in the other case,  the regularization term is doing a shrinkage on the spectrum of $X_{1:k}$ and one can see this effect through the threshold operator $\Sigma_{1, thres}^{1/2}$ appearing in \textit{case ii)} from Proposition~\ref{prop:esti_beta_1}.

\subsection{Upper bound on the price for overfitting} 
\label{sub:upper_bounds_on_the_price_for_noise_interpolation}
Following the risk decomposition \eqref{eq:risk_decomposition}, the last term we need to handle is $\norm{\Sigma_{k+1:p}^{1/2}(\hat\bbeta_{k+1:p} - \bbeta^*_{k+1:p})}_2$. As we said above, we do not expect $\hat\bbeta_{k+1:p}$ to be a good estimator of $\bbeta_2^*:=\bbeta^*_{k+1:p}$ because the minimum $\ell_2$-norm estimator $\hat\bbeta$ is using the 'remaining part' of $\bR^p$ endowed by the $p-k$ smallest eigenvectors of $\Sigma$ (we denoted this space by $V_{k+1:p}$) to interpolate the noise $\bxi$ and not to estimate $\bbeta_2^*$ that is why we call the error term $\norm{\Sigma_{k+1:p}^{1/2}(\hat\bbeta_{k+1:p} - \bbeta_2^*)}_2$ a price for noise interpolation instead of an estimation error. A consequence is that we just upper bound this term by
\begin{equation*}
\norm{\Sigma_{k+1:p}^{1/2}(\hat\bbeta_{k+1:p} - \bbeta^*_2)}_2\leq \norm{\Sigma_{k+1:p}^{1/2}\hat\bbeta_{k+1:p}}_2 + \norm{\Sigma_{k+1:p}^{1/2}\bbeta^*_2}_2.
\end{equation*}

Then we just need to find a high probability upper bound on $\norm{\Sigma_{k+1:p}^{1/2}\hat\bbeta_{k+1:p}}_2$. It follows from Proposition~\ref{prop:decomp} that $\hat\bbeta_{k+1:p} = X_{k+1:p}^\top (X_{k+1:p} X_{k+1:p}^\top)^{-1}(y-X_{1:k}\hat\bbeta_{1:k})$, hence, we have for $A:=X_{k+1:p}^\top (X_{k+1:p} X_{k+1:p}^\top)^{-1}$,
\begin{align}\label{eq:decomp_second_term}
\notag &\norm{\Sigma_{k+1:p}^{1/2}\hat\bbeta_{k+1:p}}_2= \norm{\Sigma_{k+1:p}^{1/2}A(y-X_{1:k}\hat\bbeta_{1:k})}_2\\ 
&\leq \norm{\Sigma_{k+1:p}^{1/2}AX_{1:k}(\bbeta_1^*-\hat\bbeta_{1:k})}_2 + \norm{\Sigma_{k+1:p}^{1/2}AX_{k+1:p}\bbeta_2^*}_2 + \norm{\Sigma_{k+1:p}^{1/2}A\bxi}_2
\end{align}and now we obtain high probability upper bounds on the three terms in \eqref{eq:decomp_second_term}.

We denote by $\Omega_1$ the event onto which for all $\blambda\in\bR^N$, 
\begin{equation}\label{eq:upper_dvoretzky}
\norm{\Sigma_{k+1:p}^{1/2} X_{k+1:p}^\top \blambda}_2\leq  6\left(\sqrt{\Tr(\Sigma_{k+1:p}^2)}+\sqrt{N}\norm{\Sigma_{k+1:p}}_{op}\right)\norm{\blambda}_2.
\end{equation}It follows from Proposition~\ref{prop:dvoretsky_upper_bound} that $\bP[\Omega_1]\geq 1-\exp(-N)$ (and this result holds without any extra assumption on $N$).

On $\Omega_0\cap\Omega_1$, we have
\begin{align}\label{eq:price_noise_inter_X1}
\notag &\norm{\Sigma_{k+1:p}^{1/2}AX_{1:k}(\bbeta_1^*-\hat\bbeta_{1:k})}_2 = \norm{\Sigma_{k+1:p}^{1/2} X_{k+1:p}^\top (X_{k+1:p} X_{k+1:p}^\top)^{-1} X_{1:k}(\bbeta_1^*-\hat\bbeta_{1:k})}_2\\
\notag&\leq \norm{\Sigma_{k+1:p}^{1/2} X_{k+1:p}^\top}_{op} \norm{(X_{k+1:p} X_{k+1:p}^\top)^{-1}}_{op} \norm{X_{1:k}(\bbeta_1^*-\hat\bbeta_{1:k})}_2\\
&\leq 30\sqrt{2} \frac{\left(\sqrt{N\Tr(\Sigma_{k+1:p}^2)}+N \norm{\Sigma_{k+1:p}}_{op}\right)}{\Tr(\Sigma_{k+1:p})}\norm{\Sigma_{1:k}^{1/2}(\bbeta_1^*-\hat\bbeta_{1:k})}_2
\end{align}

It follows from Bernstein's inequality that $\norm{X_{k+1:p} \bbeta_2^*}_2\leq (3/2)\sqrt{N}\norm{\Sigma_{k+1:p}^{1/2}\bbeta_2^*}_2$
holds with probability at least $1-c_0\exp(-c_1N)$. Hence, with probability at least $1-c_0\exp(-c_1N)-\bP[(\Omega_0\cap\Omega_1)^c]$, 
\begin{align*}
\norm{\Sigma_{k+1:p}^{1/2}AX_{k+1:p}\bbeta_2^*}_2&\leq \norm{\Sigma_{k+1:p}^{1/2} X_{k+1:p}^\top}_{op} \norm{(X_{k+1:p} X_{k+1:p}^\top)^{-1}}_{op} \norm{X_{k+1:p}\bbeta_2^*}_2\\
&\leq 30\sqrt{2} \frac{\sqrt{N\Tr(\Sigma_{k+1:p}^2)} + N \norm{\Sigma_{k+1:p}}_{op}}{\Tr(\Sigma_{k+1:p})}\norm{\Sigma_{k+1:p}^{1/2}\bbeta_2^*}_2.
\end{align*}

Finally, it follows from Borell's inequality that, conditionally on $\bX$, for all $t>0$ with probability at least $1-\exp(-t/2)$, $\norm{D\bxi}_2\leq \sigma_{\xi}\left(\sqrt{\Tr(DD^\top)} + \sqrt{t}\norm{D}_{op}\right)$ where $D = \Sigma_{k+1:p}^{1/2}A$. 
Let $D = \Sigma_{k+1:p}^{1/2}A$ and recall that $A = X_{k+1:p}^\top\left(X_{k+1:p}X_{k+1:p}^\top\right)^{-1}$. Then
\begin{align*}
    \Tr\left(DD^\top \right) &= \Tr\left(D^\top D\right) = \Tr\left(\left(X_{k+1:p}X_{k+1:p}^\top\right)^{-1}X_{k+1:p}\Sigma_{k+1:p} X_{k+1:p}^\top\left(X_{k+1:p}X_{k+1:p}^\top\right)^{-1}\right)\\
    &\leq \frac{\Tr\left(X_{k+1:p}\Sigma_{k+1:p}X_{k+1:p}\right)}{\norm{X_{k+1:p}X_{k+1:p}^\top}_{op}^2} = \frac{\sum_{i=1}^N\norm{\Sigma_{k+1:p}^{1/2}P_{k+1:p}X_i}_2^2}{\norm{X_{k+1:p}X_{k+1:p}^\top}_{op}^2}.
\end{align*}

As $\sum_{i=1}^N\norm{\Sigma_{k+1:p}^{1/2}P_{k+1:p}X_i}_2^2$ is sum of $N$ i.i.d. sub-exponential random variables with sub-exponential norm (up to constants) $\Tr\left(\Sigma_{k+1:p}^2\right)$, by Bernstein's inequality, with probability at least $1-\exp\left(-cN\right)$, $\sum_{i=1}^N\norm{\Sigma_{k+1:p}^{1/2}P_{k+1:p}X_i}_2^2\leq 20N\Tr\left(\Sigma_{k+1:p}^2\right)$, thus
\begin{align*}
\Tr(DD^\top) \leq \frac{\Tr\left(X_{k+1:p}\Sigma_{k+1:p}X_{k+1:p}\right)}{\norm{X_{k+1:p}X_{k+1:p}^\top}_{op}^2} \leq \frac{20N \Tr(\Sigma_{k+1:p}^2)}{\Tr^2(\Sigma_{k+1:p})}.
\end{align*}Moreover,
\begin{align*}
\norm{D}_{op} &= \norm{\Sigma_{k+1:p}^{1/2}X_{k+1:p}^\top (X_{k+1:p} X_{k+1:p}^\top)^{-1}}_{op}\leq \norm{\Sigma_{k+1:p}^{1/2}X_{k+1:p}^\top}_{op} \norm{(X_{k+1:p} X_{k+1:p}^\top)^{-1}}_{op}\\
&\leq \frac{12}{\Tr(\Sigma_{k+1:p})}\left(\sqrt{\Tr(\Sigma_{k+1:p}^2)}+\sqrt{N}\norm{\Sigma_{k+1:p}}_{op}\right).
\end{align*}hence, for all $t>0$ with probability at least $1-e^{-t/2}-c_0e^{-c_1t}$,
\begin{equation*}
\begin{aligned}
\norm{\Sigma_{k+1:p}^{1/2}A\bxi}_2\leq \sigma_\xi\bigg(&\frac{\sqrt{20N \Tr(\Sigma_{k+1:p}^2)}}{\Tr(\Sigma_{k+1:p})}  +\frac{12 \sqrt{t} }{\Tr(\Sigma_{k+1:p})}\left(\sqrt{\Tr(\Sigma_{k+1:p}^2)}+\sqrt{N}\norm{\Sigma_{k+1:p}}_{op}\right)\bigg).
\end{aligned}
\end{equation*}

Gathering the last three upper bounds in \eqref{eq:decomp_second_term}, we obtain the following result on the cost for noise interpolation.

\begin{Proposition}\label{prop:price_noise_interpolation}
There are absolute constants $c_0$ and $c_1$ such that the following holds.
We assume that there exists $k\leq \kappa_{iso} N$ such that $N\leq \kappa_{DM}d_*(\Sigma_{k+1:p}^{-1/2}B_2^p)$, then the following holds for all such $k$'s. For all $t>0$, with probability at least $1-c_0\exp(-c_1N)-\exp(-t/2)$,
\begin{align*}
&\norm{\Sigma_{k+1:p}^{1/2}(\hat\bbeta_{k+1:p} - \bbeta^*_2)}_2\leq  40\sqrt{2}  \frac{\sqrt{N\Tr(\Sigma_{k+1:p}^2)} +  N \norm{\Sigma_{k+1:p}}_{op}}{\Tr(\Sigma_{k+1:p})}\bigg(\norm{\Sigma_{k+1:p}^{1/2}\bbeta_2^*}_2 +\norm{\Sigma_{1:k}^{1/2}(\bbeta_1^*-\hat\bbeta_{1:k})}_2\bigg)\\
&+\norm{\Sigma_{k+1:p}^{1/2}\bbeta^*_2}_2+ \sigma_\xi\bigg(\frac{\sqrt{20N \Tr(\Sigma_{k+1:p}^2)}}{\Tr(\Sigma_{k+1:p})}  + \frac{12 \sqrt{t} }{\Tr(\Sigma_{k+1:p})}\bigg(\sqrt{\Tr(\Sigma_{k+1:p}^2)}+\sqrt{N}\norm{\Sigma_{k+1:p}}_{op}\bigg)\bigg).
\end{align*}
\end{Proposition}

\subsection{End of the proof of Theorem~\ref{theo:main}} 
\label{sub:end_of_the_proof_of_theorem_theo:main}
Parameter $k$ is chosen so that $N\norm{\Sigma_{k+1:p}}_{op}\leq \Tr(\Sigma_{k+1:p})$ (because $\kappa_{DM}\leq1$), in particular, $\sqrt{N \Tr(\Sigma_{k+1:p}^2)}\leq \Tr(\Sigma_{k+1:p})$. As a consequence, under the assumption of Proposition~\ref{prop:price_noise_interpolation}, we have, with probability at least $1-c_0\exp(-c_1N)-\exp(-t/2)$, 
\begin{align}\label{eq:final_beta_2_term}
\notag &\norm{\Sigma_{k+1:p}^{1/2}(\hat\bbeta_{k+1:p} - \bbeta^*_2)}_2\leq  80\sqrt{2} \norm{\Sigma_{1:k}^{1/2}(\bbeta_1^*-\hat\bbeta_{1:k})}_2  +(80\sqrt{2}+1)\norm{\Sigma_{k+1:p}^{1/2}\bbeta^*_2}_2\notag\\
&+  \sigma_\xi\frac{\sqrt{20N\Tr(\Sigma_{k+1:p}^2)}}{\Tr(\Sigma_{k+1:p})}+\sigma_\xi\left(\frac{12 \sqrt{t} }{\Tr(\Sigma_{k+1:p})}\left(\sqrt{\Tr(\Sigma_{k+1:p}^2)}+\sqrt{N}\norm{\Sigma_{k+1:p}}_{op}\right)\right).
\end{align}Theorem~\ref{theo:main} follows from the last result by choosing $t=N \Tr(\Sigma_{1:k})/\Tr(\Sigma_{k+1:p})$ when $\sigma_1 N \leq \Tr(\Sigma_{k+1:p})$ and $t=|J_1| + N \left(\sum_{j\in J_2}\sigma_j\right)/\left(\Tr(\Sigma_{k+1:p})\right)$ when $\sigma_1 N\geq \Tr(\Sigma_{k+1:p})$. Because, in the first case, when $t=N \Tr(\Sigma_{1:k})/\Tr(\Sigma_{k+1:p})$ we have $t\leq N$ and
\begin{equation*}
\frac{\sqrt{tN}\norm{\Sigma_{k+1:p}}_{op}}{\Tr(\Sigma_{k+1:p})}\leq \frac{\sqrt{N\Tr(\Sigma_{k+1:p}^2)}}{\Tr(\Sigma_{k+1:p})}
\end{equation*}and when $t=|J_1| + N \left(\sum_{j\in J_2}\sigma_j\right)/\left(\Tr(\Sigma_{k+1:p})\right)$ we have $t\leq N$ and  
\begin{equation*}
\frac{\sqrt{tN}\norm{\Sigma_{k+1:p}}_{op}}{\Tr(\Sigma_{k+1:p})}\leq \sqrt{\frac{\left|J_1\right|}{N}} +   \sqrt{\frac{\sum_{j\in J_2}\sigma_j}{\Tr(\Sigma_{k+1:p})}}. 
\end{equation*}Therefore, in both cases, the term in \eqref{eq:final_beta_2_term} is negligible in front of $\frac{12\sigma_\xi\sqrt{N\Tr(\Sigma_{k+1:p}^2)}}{\Tr(\Sigma_{k+1:p})}$
and  $\square$  defined in Proposition~\ref{prop:esti_beta_1}.


\section{Proof of Theorem~\ref{theo:main_2}} 
\label{sec:proof_of_theorem_theo:main_2}
The proof of Theorem~\ref{theo:main_2} relies on the decomposition given in Proposition~\ref{prop:decomp} of the estimator $\hat\bbeta=\hat\bbeta_1+\hat\bbeta_2$ like in the proof of Theorem~\ref{theo:main}. It therefore follows the same path as the one of the proof in the previous section (we will therefore use the same notation and detail only the main differences); in particular, it uses the excess risk decomposition \eqref{eq:risk_decomposition}. However, because $k> N$, $X_{1:k}$ cannot act anymore as an isomorphy over the entire space $V_{1:k}$ since $V_{1:k}$ is of dimension $k$ and $X_{1:k}$ has only $N<k$ rows and is therefore of rank at most $N$; in particular, $\ker(X_{1:k})\cap V_{1:k}$ is none trivial. 

However, following Theorem~\ref{theo:rip} there is a cone in $V_{1:k}$ onto which $X_{1:k}$ behaves like an isomorphy; it is given by $\cC:=\left\{v\in V_{1:k}: R_N(\Sigma_{1:k}^{1/2}B_2^p)\norm{v}_2\leq \norm{\Sigma_{1:k}^{1/2}v}_2\right\}$. We will use this \textit{restricted isomorphy property}  to lower bound the quadratic process $\cQ_{\bbeta_1}$ in \textit{case a)} (it is the case defined in the previous section where the quadratic process is used to dominate the multiplier process); we will therefore need $\bbeta_1-\bbeta_1^*\in\cC$ in \textit{case a)}. The other difference with the proof from the previous section (that is for the case $k<N$) deals with the upper bound on the multiplier process: we will use an upper bound on $\norm{\tilde \Sigma_{1:k}^{-1/2}X_{1:k}^\top}_{op}$ that does not follow from an isomorphic property of $X_{1:k}$ but from the upper side of DM's theorem like the one given in Proposition~\ref{prop:dvoretsky_upper_bound}. The stochastic properties we use on $X_{k+1:p}$ in the case $k>N$ are the same as the one used in the previous proof since having $k<N$ or $k>N$ does not play any role on the behavior of $X_{k+1:p}$.   

We gather all the stochastic properties we need on $X_{1:k}$ and $X_{k+1:p}$ to obtain our estimation results for $\hat\bbeta_1$ and $\hat\bbeta_2$ in the case $k>N$ in the following event.

\subsection{Stochastic event behind Theorem~\ref{theo:main_2}} We denote by $\Omega_0^\prime$ the event onto which the following three geometric properties hold:
\begin{itemize}
  \item for all $\blambda\in\bR^N$, $$(1/(2\sqrt{2}))\sqrt{\Tr(\Sigma_{k+1:p})}\norm{\blambda}_2\leq \norm{X_{k+1:p}^\top\blambda}_2\leq (3/2)\sqrt{\Tr(\Sigma_{k+1:p})}\norm{\blambda}_2$$
  \item for all $\bbeta_1$ in the cone  $\cC:=\left\{v\in V_{1:k}: R_N(\Sigma_{1:k}^{1/2}B_2^p) \norm{v}_2\leq\norm{\Sigma_{1:k}^{1/2}v}_2\right\}$, 
\begin{equation*}
  \frac{1}{2}\norm{\Sigma_{1:k}^{1/2}\bbeta_1}_2\leq \frac{1}{\sqrt{N}}\norm{X_{1:k}\bbeta_1}_2\leq \frac{3}{2}\norm{\Sigma_{1:k}^{1/2}\bbeta_1}_2,
\end{equation*}
\item for all $\blambda\in\bR^N$, 
\begin{equation}\label{eq:third_point_event_omega_prime}
\norm{\tilde\Sigma_1^{-1/2}X_{1:k}^\top\blambda}_2\leq  3\left(\sqrt{\square^2|J_1| + \triangle^2 \sum_{j\in J_2}\sigma_j}+\sqrt{N}\sigma(\square, \triangle)\right)\norm{\blambda}_2
\end{equation}where $\tilde \Sigma_1^{1/2}$ has been introduced in \eqref{eq:tilde_D} and $\sigma(\square, \triangle)$ in Lemma~\ref{lem:first_beta_2_star_term}. However, in the following we will always assume that
\begin{equation}\label{eq:diameter}
\sqrt{\square^2|J_1| + \triangle^2 \sum_{j\in J_2}\sigma_j}\leqslant\sqrt{N}\sigma(\square, \triangle) 
\end{equation}where $\sigma(\square, \triangle)$ is defined in Lemma~\ref{lem:first_beta_2_star_term}. We will explain the reason in Section~\ref{sub:estimation_properties_of_the_2}.
\end{itemize}

If $N\leq \kappa_{DM}d_*(\Sigma_{k+1:p}^{-1/2}B_2^p)$, we know from Theorem~\ref{theo:DM} and Theorem~\ref{theo:rip} that the first two points defining the event $\Omega_0^\prime$ hold simultaneously with probability at least $1-c_0\exp(-c_1N)$. It only remains to handle the last point. To that end we use Proposition~\ref{prop:dvoretsky_upper_bound}: with probability at least $1-c_0\exp(- c_1N)$, for all $\blambda\in\bR^N$, 
\begin{align*}
&\norm{\tilde\Sigma_{1:k}^{-1/2}X_{1:k}^\top\blambda}_2 = \norm{\tilde\Sigma_{1:k}^{-1/2}\Sigma_{1:k}^{1/2}\bG^\top\blambda}_2\leq 2 \left(\sqrt{\Tr\left((\tilde\Sigma_{1:k}^{-1/2}\Sigma_{1:k}^{1/2})(\tilde\Sigma_{1:k}^{-1/2}\Sigma_{1:k}^{1/2})^\top\right)}+\sqrt{c_2 N}\norm{\tilde\Sigma_{1:k}^{-1/2}\Sigma_{1:k}^{1/2}}_{op}\right)\norm{\blambda}_2
\end{align*}where $\bG$ is a $N\times p$ standard Gaussian matrix. The spectrum of $\tilde\Sigma_{1:k}^{-1/2}\Sigma_{1:k}^{1/2}$ is the same as the one of $\tilde D_{1}^{-1/2}D_{1:k}^{1/2}$ given by $p-k$ zeros and $\left(\square\sqrt{\sigma_j}/\max(\sqrt{\sigma_j}, (\square/\triangle)):j=1,\ldots, k\right)$. Therefore, the third point of event $\Omega_0^\prime$ holds with probability at least $1-\exp(-N)$. We conclude that $\bP[\Omega_0^\prime]\geq 1-2c_0\exp(-c_1N)$. 

Note that for the choice of $\square$ and $\triangle$ such that $(\square/\triangle)^2 = \kappa_{DM}\ell_*^2(\Sigma_{k+1:p}^{1/2}B_2^p)/N$ that we will make later, \eqref{eq:diameter} is equivalent to 
\begin{equation}\label{eq:equiv_eq_41}
\left\{
\begin{array}{cc}
\sum_{j\in J_2}\sigma_j \leq \kappa_{DM}\ell_*^2(\Sigma_{k+1:p}^{1/2}B_2^p) \left(1-\frac{|J_1|}{N}\right) & \mbox{ when } N \sigma_1 \geq \kappa_{DM}\ell_*^2(\Sigma_{k+1:p}^{1/2}B_2^p)\\
\Tr(\Sigma_{1:k})\leq N \sigma_1 & \mbox{ otherwise.}
\end{array}\right.
\end{equation}

We place ourselves on the event $\Omega_0^\prime$ up to the end of the proof. As in the proof of  Theorem~\ref{theo:main}, we split our analysis into three subsections: one for the study of  $\hat\bbeta_1$, one for $\hat\bbeta_2$ and the last one where the two previous sections are merged.

\subsection{Estimation properties of the 'ridge estimator' $\hat\bbeta_{1:k}$; case $k\gtrsim N$} 
\label{sub:estimation_properties_of_the_2}

 We are now providing some details on the arguments we need to prove Theorem~\ref{theo:main_2} which are different from the one of Theorem~\ref{theo:main}. Let us first handle the lower bound we need on the quadratic process in \textit{case a)} that is when $\bbeta_1$ is such that $\norm{\Sigma^{1/2}(\bbeta_1-\bbeta_1^*)}_2=\square$ and $\norm{\bbeta_1-\bbeta_1^*}_2\leq \triangle$. We can only use the restricted isomorphy property satisfied by $X_{1:k}$ over the cone $\cC$ so we have to insure that $\bbeta_1-\bbeta_1^*$ lies in that cone. That is the case when $\square/\triangle \geq R_N(\Sigma_{1:k}^{1/2}B_2^p)$ since, in that case,  we have
\begin{equation*}
\begin{aligned}
 R_N(\Sigma_{1:k}^{1/2}B_2^p) \norm{\bbeta_1-\bbeta_1^*}_2\leq  R_N(\Sigma_{1:k}^{1/2}B_2^p)\triangle =  R_N(\Sigma_{1:k}^{1/2}B_2^p)\frac{\triangle}{\square}\norm{\Sigma^{1/2}(\bbeta_1-\bbeta_1^*)}_2\leq \norm{\Sigma^{1/2}(\bbeta_1-\bbeta_1^*)}_2.
\end{aligned}
\end{equation*} That is the reason why, in the case $k>N$ we do have the extra condition $\square/\triangle \geq R_N(\Sigma_{1:k}^{1/2}B_2^p)$ (note that when $k\leq \kappa_{RIP}N$ then $R_N(\Sigma_{1:k}^{1/2}B_2^p)=0$ and so there is no need for this condition). As a consequence, when $\square/\triangle \geq R_N(\Sigma_{1:k}^{1/2}B_2^p)$ we have $\bbeta_1-\bbeta_1^*\in\cC$ in \textit{case a)}  and so (on the event $\Omega_0^\prime$), 
\begin{align*}
\cQ_{\bbeta_1} = \norm{(X_{k+1:p}X_{k+1:p}^\top)^{-1/2}X_{1:k}(\bbeta_1-\bbeta_1^*)}_2^2\geq \frac{1}{8 \Tr(\Sigma_{k+1:p})}\norm{X_{1:k}(\bbeta_1-\bbeta_1^*)}_2^2\geq  \frac{N\square^2}{32\Tr(\Sigma_{k+1:p})}.
\end{align*}We therefore recover the same lower bound as in the proof of Theorem~\ref{theo:main} (see \eqref{eq:lower_bound_quad_case_a}).

Let us now handle the multiplier process. It follows from \eqref{eq:multiplier_decomp} that we need to upper bound the two quantities
\begin{equation*}
    \norm{\tilde \Sigma_1^{-1/2}X_{1:k}^\top (X_{k+1:p}X_{k+1:p}^\top)^{-1}X_{k+1:p}\bbeta_2^*}_2\mbox{ and }
    \norm{\tilde \Sigma_1^{-1/2}X_{1:k}^\top (X_{k+1:p}X_{k+1:p}^\top)^{-1}\bxi}_2
\end{equation*}
as we did in Lemma~\ref{lem:first_beta_2_star_term} and Lemma~\ref{lem:noise_term} but without the isomorphic property of $X_{1:k}$ on $V_{1:k}$. 

We know from the third point of the event $\Omega_0^\prime$ that 
\begin{equation}\label{eq:upper_bound_largest_sing_val_prod}
\norm{\tilde\Sigma_{1:k}^{-1/2}X_{1:k}^\top}_{op}\leq  6\sqrt{N}\sigma(\square, \triangle).
\end{equation}
As a consequence, using similar argument as in the proof of Lemma~\ref{lem:first_beta_2_star_term}, we have with probability at least $1-c_0\exp(-c_1N)$,
\begin{align}\label{eq:final_1_2}
\notag&\norm{\tilde \Sigma_1^{-1/2}X_{1:k}^\top (X_{k+1:p}X_{k+1:p}^\top)^{-1}X_{k+1:p}\bbeta_2^*}_2\leq \norm{\tilde\Sigma_{1:k}^{-1/2}X_{1:k}^\top}_{op}\norm{(X_{k+1:p}X_{k+1:p}^\top)^{-1}}_{op} \cdot \\&\norm{X_{k+1:p}\bbeta_2^*}_2\leq \sigma(\square, \triangle) \frac{54N}{\Tr(\Sigma_{k+1:p})}\norm{\Sigma_{k+1:p}^{1/2}\bbeta_2^*}_2
\end{align}which is up to absolute constants the same result as in Lemma~\ref{lem:first_beta_2_star_term}.

Next, we prove a high probability  upper bound on $\norm{\tilde \Sigma_1^{-1/2}X_{1:k}^\top (X_{k+1:p}X_{k+1:p}^\top)^{-1}\bxi}_2$. It follows from Borell's inequality that for all $t>0$, with probability at least $1-\exp(-t/2)$,
\begin{equation*}
\norm{\tilde \Sigma_1^{-1/2}X_{1:k}^\top (X_{k+1:p}X_{k+1:p}^\top)^{-1}\bxi}_2 \leq \sigma_\xi\left(\sqrt{\Tr(DD^\top)} + \sqrt{t}\norm{D}_{op}\right)
\end{equation*}where $D = \tilde \Sigma_1^{-1/2}X_{1:k}^\top (X_{k+1:p}X_{k+1:p}^\top)^{-1}$. We know from Bernstein's inequality that with probability at least $1-\exp(-c_1N)$, $\sum_{i=1}^N\norm{\tilde\Sigma_{1:k}^{-1/2}P_{1:k}X_i}_2^2\leq 2N\left(\left|J_1\right|\square + \triangle^2\sum_{j\in J_2}\sigma_j\right)$. Hence, with probability at least $1-c_0\exp(-c_1N)$,
\begin{equation*}
\sqrt{\Tr(DD^\top)} = \sqrt{\Tr\left(D^\top D\right)} \leq \frac{16\sqrt{N}}{\Tr(\Sigma_{k+1:p})}\sqrt{ \left|J_1 \right|\square^2 + \triangle^2\sum_{j\in J_2}\sigma_j}.
\end{equation*}We also have
\begin{equation*}
\norm{D}_{op}  \leq \norm{\tilde \Sigma_1^{-1/2} X_{1:k}^\top }_{op}  \norm{(X_{k+1:p}X_{k+1:p}^\top)^{-1}}_{op} \leq \frac{24\sqrt{N}\sigma(\square, \triangle)}{\Tr(\Sigma_{k+1:p})}.
\end{equation*}

Finally, using the same other arguments as in the proof of Lemma~\ref{lem:noise_term} we obtain the same bound (up to absolute constants) as in Lemma~\ref{lem:noise_term}:  with probability at least $1-c_0\exp(-c_1N)-\exp(-t(\square, \triangle)/2)$, 
\begin{equation}\label{eq:final_2_2}
\norm{\tilde \Sigma_1^{-1/2}X_{1:k}^\top (X_{k+1:p}X_{k+1:p}^\top)^{-1}\bxi}_2\leq \frac{32\sqrt{N}\sigma_\xi}{\Tr(\Sigma_{k+1:p})}\sqrt{ \left|J_1 \right|\square^2 + \triangle^2\sum_{j\in J_2}\sigma_j}.
\end{equation}

As a consequence, all the machinery used in Section~\ref{sub:estimation_properties_of_the_} also applies for the same choice of $\square$ and $\triangle$(up to absolute constants) under the three extra conditions that $\square/\triangle\geqslant R_N(\Sigma_{1:k}^{1/2}B_2^p)$, that one from \eqref{eq:equiv_eq_41}. We therefore end up with almost the same result as Proposition~\ref{prop:esti_beta_1} regarding the estimation property of $\hat\bbeta_{1:k}$:

\begin{Proposition}\label{prop:esti_beta_1_2}There are absolute constants $c_0$, $c_1$ and $c_3$ such that the following holds.
We assume that there exists $k\in[p]$ such that $N\leq \kappa_{DM}d_*(\Sigma_{k+1:p}^{-1/2}B_2^p)$, \eqref{eq:equiv_eq_41} holds and $\kappa_{DM}\ell_*^2(\Sigma_{k+1:p}^{1/2}B_2^p)\geq R_N(\Sigma_{1:k}^{1/2}B_2^p)^2 N$, then the following holds for all such $k$'s. 
With probability at least $1-p^*$,
\begin{equation*}
    \norm{\Sigma_1^{1/2}(\hat\bbeta_1-\bbeta_1^*)}_2\leq \square \mbox{ and } \norm{\hat\bbeta_1-\bbeta_1^*}_2\leq \square\sqrt{\frac{N}{\kappa_{DM}\ell_*^2(\Sigma_{k+1:p}^{1/2}B_2^p)}},
\end{equation*}
where, 
\begin{itemize}
  \item[i)] if $\sigma_1 N\leq \kappa_{DM}\ell_*^2(\Sigma_{k+1:p}^{1/2}B_2^p)$, $p^* = c_0\exp(-c_1N)$ and 
  \begin{equation*}
\square = c_3\max \left\{\sigma_\xi\sqrt{\frac{\Tr(\Sigma_{1:k})}{\Tr(\Sigma_{k+1:p})}}, \sqrt{\frac{N\sigma_1}{\Tr(\Sigma_{k+1:p})}}\norm{\Sigma_{k+1:p}^{1/2}\bbeta_2^*}_{2}, \norm{\bbeta_1^*}_2\sqrt{\frac{\Tr(\Sigma_{k+1:p})}{N}}\right\}
\end{equation*}
\item[ii)]   if $\sigma_1 N\geq \kappa_{DM}\ell_*^2(\Sigma_{k+1:p}^{1/2}B_2^p)$, $p^* = c_0\exp\left(-c_1\left(|J_1| + N \left(\sum_{j\in J_2}\sigma_j\right)/\left(\Tr(\Sigma_{k+1:p})\right)\right)\right)$ and
\begin{equation*}
\square = c_3\max \left\{\sigma_\xi \sqrt{\frac{\left|J_1\right|}{N}},  \sigma_\xi\sqrt{\frac{\sum_{j\in J_2}\sigma_j}{\Tr(\Sigma_{k+1:p})}}, \norm{\Sigma_{k+1:p}^{1/2}\bbeta_2^*}_{2}, \norm{\Sigma_{1,thres}^{-1/2}\bbeta_1^*}_2\frac{\Tr(\Sigma_{k+1:p})}{N}\right\}.
\end{equation*}
\end{itemize}

\end{Proposition}

\begin{Remark}
If we use $\norm{\tilde\Sigma_1^{-1/2}X_{1:k}^\top}_{op}\leq  6\left(\sqrt{\square^2|J_1| + \triangle^2 \sum_{j\in J_2}\sigma_j}\right),$ which holds when $N\leq \kappa_{DM}t(\square, \triangle)$ instead of \eqref{eq:third_point_event_omega_prime}, we will end up with
\begin{equation*}
\begin{aligned}
\square = C_0^2 \max\bigg\{\left(\sigma_\xi\vee \norm{\Sigma_{k+1:p}^{1/2}\bbeta_2^*}_2\right)\sqrt{\frac{\left|J_1\right|}{N}}, \left(\sigma_\xi\vee \norm{\Sigma_{k+1:p}^{1/2}\bbeta_2^*}_2\right)\sqrt{\frac{\sum_{j\in J_2}\sigma_j}{\Tr(\Sigma_{k+1:p})}},\norm{\Sigma_{1, thres}^{-1/2}\bbeta_1^*}_2 \frac{\Tr(\Sigma_{k+1:p})}{N}\bigg\}.
\end{aligned}
\end{equation*}
However, the latter rate does not showcase the benign overfitting phenomenon. Indeed,  if one has $N\leq t(\square, \triangle)$ then
\begin{itemize}
\item when $\sigma(\square, \triangle) = \square$: in this case, $t(\square, \triangle)=\left|J_1\right|+N\sum_{j\in J_2}\sigma_j/\Tr(\Sigma_{k+1:p})\geq N$, which leads to $\frac{\sum_{j\in J_2}\sigma_j}{\Tr(\Sigma_{k+1:p})}\geqslant 1-\frac{\left|J_1\right|}{N}.$
As $\left|J_1\right|/N$ tends to $0$, we have $\sum_{j\in J_2}\sigma_j/\Tr(\Sigma_{k+1:p})$ tends to $1$, on which we will never have BO.
\item when $\sigma(\square, \triangle) = \triangle \sqrt{\sigma_1}$, which is equivalent to $\sigma_1 N \lesssim \Tr(\Sigma_{k+1:p})$. In this case, $\frac{\Tr(\Sigma_{1:k})}{\sigma_1} = t(\square, \triangle)  \geq N.$
Therefore, $\Tr(\Sigma_{1:k})/\Tr(\Sigma_{k+1:p})\geqslant \sigma_1 N /\Tr(\Sigma_{k+1:p})$. However, as $\Tr(\Sigma_{k+1:p})/N$ tends to $0$, this will not converge.
\end{itemize}
\end{Remark}

\subsection{Upper bound on the price for noise interpolation; the case $k\gtrsim N$} 
\label{sub:upper_bounds_on_the_price_for_noise_interpolation_larger}
The aim of this section is to obtain a high probability bound on  $\norm{\Sigma_{k+1:p}^{1/2}(\hat\bbeta_{k+1:p} - \bbeta^*_{k+1:p})}_2$. We use the same decomposition as in \eqref{eq:decomp_beta_2_star} together with the closed-form of $\hat\bbeta_{k+1:p}$ given in Proposition~\ref{prop:decomp}:
\begin{align}\label{eq:decomp_beta_2_star_2}
\notag &\norm{\Sigma_{k+1:p}^{1/2}(\hat\bbeta_{k+1:p} - \bbeta^*_{k+1:p})}_2\leq \norm{\Sigma_{k+1:p}^{1/2}A(y-X_{1:k}\hat\bbeta_{1:k})}_2 + \norm{\Sigma_{k+1:p}^{1/2}\bbeta^*_{k+1:p}}_2\\
&\leq  \norm{\Sigma_{k+1:p}^{1/2}AX_{1:k}(\bbeta_1^*-\hat\bbeta_{1:k})}_2 + \norm{\Sigma_{k+1:p}^{1/2}AX_{k+1:p}\bbeta_2^*}_2 +  \norm{\Sigma_{k+1:p}^{1/2}A\bxi}_2+ \norm{\Sigma_{k+1:p}^{1/2}\bbeta^*_{k+1:p}}_2
\end{align}where $A=X_{k+1:p}^\top (X_{k+1:p} X_{k+1:p}^\top)^{-1}$. We will handle the last three terms of \eqref{eq:decomp_beta_2_star_2} as in the Section~\ref{sub:upper_bounds_on_the_price_for_noise_interpolation} because they do not depend on $X_{1:k}$. However, for the first term in \eqref{eq:decomp_beta_2_star_2}, we used in the proof of Theorem~\ref{theo:main} the isomorphic property of $X_{1:k}$ in \eqref{eq:price_noise_inter_X1}: $\norm{X_{1:k}(\bbeta_1^*-\hat\bbeta_{1:k})}_2\leq 3\sqrt{N} \norm{\Sigma_{1:k}^{1/2}(\bbeta_1^*-\hat\bbeta_{1:k})}_2/2$. We cannot use this isomorphic property in the case $k\gtrsim N$ since it does not hold on the entire space $V_{1:k}$. However, this last inequality still holds if one can show that $\bbeta_1^*-\hat\bbeta_{1:k}$ lies in the cone $\cC$; that is to show that $R_N(\Sigma_{1:k}^{1/2}B_2^p) \norm{\bbeta_1^*-\hat\bbeta_{1:k}}_2\leq\norm{\Sigma_{1:k}^{1/2}(\bbeta_1^*-\hat\bbeta_{1:k})}_2$. This type of condition usually holds when we regularize by a norm (this is for instance the case of the LASSO) but it is not clear if this holds when the regularization function is the \textit{square of a norm} as it is the case for the 'ridge estimator' $\hat\bbeta_1$. Therefore, we cannot use this argument here. The way we will handle this issue is by showing an upper bound on  $\norm{X_{1:k}(\bbeta_1^*-\hat\bbeta_{1:k})}_2$ directly without going through a norm equivalence with $\norm{\Sigma_1^{1/2}\cdot}_2$.

We know from Proposition~\ref{prop:decomp} that $\cL_{\hat\bbeta_1}\leq 0$ where $\cL_{\hat\bbeta_1} = \cQ_{\hat\bbeta_1} + \cM_{\hat\bbeta_1} + \cR_{\hat\bbeta_1}$ and for all $\bbeta_1\in V_{1:k}$,
\begin{equation*}
\begin{aligned}
\cQ_{\bbeta_1} = \norm{(X_{k+1:p}X_{k+1:p}^\top)^{-1/2}X_{1:k}(\bbeta_1-\bbeta_1^*)}_2^2, \cM_{\bbeta_1}= 2 \inr{X_{1:k}^\top (X_{k+1:p}X_{k+1:p}^\top)^{-1}(X_{k+1:p}\bbeta_2^*+\bxi) - \bbeta_1^*,\bbeta_1-\bbeta_1^*}
\end{aligned}
\end{equation*} and  $\cR_{\bbeta_1}= \norm{\bbeta_1-\bbeta_1^*}_2^2$. We therefore have $\cQ_{\hat\bbeta_1} + \cR_{\hat\bbeta_1}\leq |\cM_{\hat\bbeta_1}|$. We are now proving a lower bound on $\cQ_{\hat\bbeta_1}$ and an upper bound on $|\cM_{\hat\bbeta_1}|$. For the lower bound on $\cQ_{\hat\bbeta_1}$, we use that on $\Omega_0^\prime$,
  \begin{equation}\label{eq:lower_bound_quad_case_2}
 \cQ_{\hat\bbeta_1} = \norm{(X_{k+1:p}X_{k+1:p}^\top)^{-1/2}X_{1:k}(\hat\bbeta_1-\bbeta_1^*)}_2^2  \geq \frac{1}{32 \Tr(\Sigma_{k+1:p})}\norm{X_{1:k}(\hat\bbeta_1-\bbeta_1^*)}_2^2. 
 \end{equation} For the upper bound on $|\cM_{\hat\bbeta_1}|$, we consider the norm (restricted to $V_{1:k}$) defined for all $\bbeta\in V_{1:k}$
\begin{equation}\label{eq:interpol_norm_2}
\vertiiii{\bbeta} := \max\left( \norm{\Sigma_{1:k}^{1/2}\bbeta}_2, \sqrt{\frac{\kappa_{DM}\ell_*^2(\Sigma_{k+1:p}^{1/2}B_2^p)}{N}} \norm{\bbeta}_2\right)
\end{equation} and we set $\vertiiii{\bbeta} = 0$ for all $\bbeta\in V_{k+1:p}$. On the event $\Omega_0^{\prime}$, we have $\vertiiii{\hat\bbeta_1-\bbeta_1^*}\leq \square$ where $\square$ is defined in Proposition~\ref{prop:esti_beta_1_2}. Therefore, if we define $\theta : = \sup\left(|\cM_{\bbeta_1}|:\vertiiii{\bbeta_1-\bbeta_1^*}\leq 1\right)$  we have
\begin{equation*}
  \norm{\hat \bbeta_1-\bbeta_1^*}_2^2 + \frac{1}{32 \Tr(\Sigma_{k+1:p})}\norm{X_{1:k}(\bbeta_1-\bbeta_1^*)}_2^2\leq \theta \square.
 \end{equation*} and so we have  $\norm{X_{1:k}(\hat\bbeta_1-\bbeta_1^*)}_2\leqslant 4\sqrt{2} \sqrt{\theta \square\Tr(\Sigma_{k+1:p})}$ (we also have  $\norm{\hat\bbeta_1-\bbeta_1^*}_2\leq  \sqrt{\theta\square}$ however we will not use it here). The last result we need to prove is a high probability upper bound on $\theta$. It appears that we already did it in the previous Section~\ref{sub:estimation_properties_of_the_2} since for the choice of $\square$ and $\triangle$ such that  $(\square/\triangle)^2 = \kappa_{DM}\ell_*^2(\Sigma_{k+1:p}^{1/2}B_2^p)/N$ we have $\square\vertiii{\bbeta} =  \vertiiii{\bbeta}$ for all $\bbeta\in\bR^p$. Hence, if we denote by $\vertiiii{\cdot}_*$ the dual norm of $\vertiiii{\cdot}$, we have $\vertiii{\cdot}_* = \square\vertiiii{\cdot}_*$ and so $\vertiiii{\cdot}_*$ is equivalent to $\square^{-1}\norm{\tilde \Sigma_1^{-1/2}\cdot}_2$ where $\tilde \Sigma_1^{-1/2}$ is defined in \eqref{eq:tilde_D}. Therefore, using the lower bound for $\cQ_{\bbeta_1}/\square$(since the quadratic term dominates the multiplier term) as the upper bound for the following decomposition
\begin{equation*}
\theta/2\leq \vertiiii{\bbeta_1^*}_* + \vertiiii{X_{1:k}^\top (X_{k+1:p}X_{k+1:p}^\top)^{-1}X_{k+1:p}\bbeta_2^*}_* + \vertiiii{X_{1:k}^\top (X_{k+1:p}X_{k+1:p}^\top)^{-1}\bxi}_*
\end{equation*}we obtain the following result:

\begin{Proposition}\label{prop:esti_beta_1_X_1}
Under the same assumptions as in Proposition~\ref{prop:esti_beta_1_2}, with probability at least $1-p^*$,
\begin{equation*}
    \norm{X_{1:k}(\hat\bbeta_{1:k} - \bbeta_1^*)}_2\leq 4\sqrt{N} \square,
\end{equation*}
where $p^*$ and $\square$ are defined in Proposition~\ref{prop:esti_beta_1_2}.
\end{Proposition}

The following result follows from the decomposition \eqref{eq:decomp_beta_2_star_2} and Section~\ref{sub:upper_bounds_on_the_price_for_noise_interpolation} and  Proposition~\ref{prop:esti_beta_1_X_1}, for the control of the four terms in this decomposition. Namely, if $N\geqslant 5\log{p}$, then

\begin{equation*}
\norm{\Sigma_{k+1:p}^{1/2}AX_{k+1:p}\bbeta_2^*}_2\leq 30\sqrt{2} \frac{\sqrt{N\Tr(\Sigma_{k+1:p}^2)} + N \norm{\Sigma_{k+1:p}}_{op}}{\Tr(\Sigma_{k+1:p})}\norm{\Sigma_{k+1:p}^{1/2}\bbeta_2^*}_2,
\end{equation*}

and

\begin{equation*}
\begin{aligned}
\norm{\Sigma_{k+1:p}^{1/2}A\bxi}_2\leq \sigma_\xi\bigg(\frac{\sqrt{20N \Tr(\Sigma_{k+1:p}^2)}}{\Tr(\Sigma_{k+1:p})}  + \frac{12 \sqrt{t} }{\Tr(\Sigma_{k+1:p})}\bigg(\sqrt{\Tr(\Sigma_{k+1:p}^2)}+\sqrt{N}\norm{\Sigma_{k+1:p}}_{op}\bigg)\bigg)
\end{aligned}
\end{equation*}

holds with probability at least $1-\exp(-t/2)-c_0\exp(-c_1 N)-\bP[(\Omega_0'\cap \Omega_1)^c]$, where $\Omega_1$ is defined in \eqref{eq:upper_dvoretzky}.

\begin{Proposition}\label{prop:price_noise_interpolation_2}
There are absolute constants $c_0$ and $c_1$ such that the following holds.
We assume that there exists $k>N$ such that $N\leq \kappa_{DM}d_*(\Sigma_{k+1:p}^{-1/2}B_2^p))$, then the following holds for all such $k$'s. For all $t>0$, with probability at least $1-2c_0\exp(-c_1N)-\exp(-t(\square, \triangle)/2)-\exp(-t/2)$,
\begin{align*}
&\norm{\Sigma_{k+1:p}^{1/2}(\hat\bbeta_{k+1:p} - \bbeta^*_2)}_2\leq  40\sqrt{2}  \frac{\sqrt{N\Tr(\Sigma_{k+1:p}^2)} +  N \norm{\Sigma_{k+1:p}}_{op}}{\Tr(\Sigma_{k+1:p})}\bigg(\norm{\Sigma_{k+1:p}^{1/2}\bbeta_2^*}_2 +\norm{\Sigma_{1:k}^{1/2}(\bbeta_1^*-\hat\bbeta_{1:k})}_2\bigg)\\&+\norm{\Sigma_{k+1:p}^{1/2}\bbeta^*_2}_2+\sigma_\xi\bigg(\frac{\sqrt{20N \Tr(\Sigma_{k+1:p}^2)}}{\Tr(\Sigma_{k+1:p})}  + \frac{12 \sqrt{t} }{\Tr(\Sigma_{k+1:p})}\bigg(\sqrt{\Tr(\Sigma_{k+1:p}^2)}+\sqrt{N}\norm{\Sigma_{k+1:p}}_{op}\bigg)\bigg).
\end{align*}
\end{Proposition}

\subsection{End of the proof of Theorem~\ref{theo:main_2}} 
\label{sec:end_of_the_proof_of_theorem_theo:main_2}
We use the decomposition of the excess from \eqref{eq:risk_decomposition}, the result on the estimation property of $\hat\bbeta_1$ from Proposition~\ref{prop:esti_beta_1_2} and the one on $\hat\bbeta_2$ from Proposition~\ref{prop:price_noise_interpolation_2} to conclude. Since the upper bounds on $\norm{\Sigma_{k+1:p}^{1/2}(\hat\bbeta_{k+1:p}-\bbeta_2^*)}_2$ and $\norm{\Sigma_{1:k}^{1/2}(\hat\bbeta_{1:k}-\bbeta_1^*)}_2$ are the same up to absolute constants as that in Section~\ref{sub:estimation_properties_of_the_} and Section~\ref{sub:upper_bounds_on_the_price_for_noise_interpolation}, all the choices of $t$ from Section~\ref{sub:end_of_the_proof_of_theorem_theo:main} also applies for the case $k\gtrsim N$ under the two extra assumptions: $\Tr(\Sigma_{k+1:p})\geqslant R_N(\Sigma_{1:k}^{1/2}B_2^p)^2 N$ and \eqref{eq:equiv_eq_41}.

\section{Proof of Theorem~\ref{theo:lower_bound}} 
\label{sec:proof_of_theorem_theo:lower_bound}
The proof of Theorem~\ref{theo:lower_bound} is based on three argument: a) a slightly adapted randomization argument from Proposition~1 in \cite{DBLP:conf/aistats/RichardsMR21} that allows to replace the Bayesian risk from \cite{MR4263288,TB21} by the true risk; b) the argument from \cite{MR4263288,TB21} on $V_{1:k^*_b}$ and c) the DM theorem on $V_{k^*_b+1:p}$.  

We denote by $\Omega^*$ the event onto which for all $\blambda\in\bR^N$, 
\begin{equation*}
(1/(2\sqrt{2}))\sqrt{\Tr(\Sigma_{k^*_b+1:p})}\norm{\blambda}_2\leq \norm{X_{k^*_b+1:p}^\top\blambda}_2\leq (3/2)\sqrt{\Tr(\Sigma_{k^*_b+1:p})}\norm{\blambda}_2.
\end{equation*}It follows from Theorem~\ref{theo:DM} and \eqref{eq:dvoretsky_dim_ellipsoid} that when $b\geq 4/\kappa_{DM}$, $\bP[\Omega^*]\geq 1-\exp(-c_0 N)$ if we are in the Gaussian scenarii, and from Theorem~\ref{theo:DM_weak_moments} that when $b\geq 4/\kappa_{DM}$, $\bP[\Omega^*]\geq 1-\gamma - N^{-2c_1}$.

 We start with a bias/variance decomposition of the risk as in \cite{MR4263288,TB21}. We have 
\begin{equation}\label{eq:bias_variance_LB}
\bE \norm{\Sigma^{1/2}(\hat\bbeta - \bbeta^*)}_2^2 =  \bE\norm{\Sigma^{1/2}(\bX^\top (\bX \bX^\top)^{-1}\bX - I_p)\bbeta^*}_2^2 + \bE \norm{\Sigma^{1/2}\bX^\top (\bX \bX^\top)^{-1}\bxi }_2^2
\end{equation}where we used that $\hat\bbeta = \bX^\top (\bX \bX^\top)^{-1}y =\bX^\top (\bX \bX^\top)^{-1}(\bX \bbeta^* + \bxi)$.

For the variance term (second term in the right-hand side of  \eqref{eq:bias_variance_LB}) we use Lemma~2 and Theorem~5 from \cite{TB21} (see also Lemma~16 in \cite{MR4263288} for a similar result).

\begin{Proposition}\label{prop:lower_bound_variance}[Lemma~2 and  Theorem~5 in \cite{TB21}]There exists some absolute constant $c_0>0$ such that the following holds.
Let $b\geq 1/N$ and assume that $k^*_b<N/4$. We have
\begin{equation*}
\bE \norm{\Sigma^{1/2}\bX^\top (\bX \bX^\top)^{-1}\bxi }_2^2\geq  \frac{c_0 \sigma_\xi^2}{\max((2+b)^2, (1+2b^{-1})^2)}\left(\frac{k^*_b}{N} + \frac{N \Tr(\Sigma^2_{k^*_b+1:p})}{\Tr^2(\Sigma_{k^*_b+1:p})}\right).
\end{equation*}
\end{Proposition}

For the bias term (first term in the right-hand side of  \eqref{eq:bias_variance_LB}), we cannot use the results from \cite{MR4263288,TB21} because they hold either for one specific $\bbeta^*$ (see Theorem~4 in \cite{MR4263288}) or for some Bayesian risk (see Lemma~3 in \cite{TB21}). Moreover Lemma~3 in \cite{TB21} requires some extra assumptions on the smallest singular values of some matrix (see the condition \textit{'...for any $j>k$, w.p. at least $1-\delta$, $\mu_n(A_{-j})\geq L^{-1}(\sum_{j>k}\sigma_j)$.'} in there). The following lower bound holds for the prediction risk of $\hat\bbeta$ for the estimation of $\bbeta^*$ itself. It holds for any given $\bbeta^*$ and not a random one. It also holds under only the assumption that $k^*_b\lesssim N$ (which is a necessary condition for BO according to Theorem~1 in \cite{MR4263288}).

\begin{Proposition}\label{prop:lower_bound_bias} There is an absolute constant $c_0>0$ such that the following holds. If $N\geq c_0$ then for any $b\geq \max(4/\kappa_{DM}, 24)$, we have
\begin{equation*}
\begin{aligned}
 \bE\norm{\Sigma^{1/2}(\bX^\top (\bX \bX^\top)^{-1}\bX - I_p)\bbeta^*}_2^2 \geq \frac{1}{18 b^2} \bigg(& \norm{\Sigma_{k^*_b+1:p}^{1/2}\bbeta_{k^*_b+1:p}^*}_{2}^2,\norm{\Sigma_{1:k^*_b}^{-1/2}\bbeta_{1:k^*_b}^*}_2^2\left(\frac{\Tr(\Sigma_{k^*_b+1:p})}{N}\right)^2\bigg).
\end{aligned}
\end{equation*}
\end{Proposition}
\textit{Proof of Proposition~\ref{prop:lower_bound_bias}.} We show that the bias term can be written as a 'Bayesian bias term' using a similar argument as in Proposition~1 from \cite{DBLP:conf/aistats/RichardsMR21}. Let $U$ be an $p\times p$ orthogonal matrix which commutes with $\Sigma$. We have
\begin{align*}
\norm{\Sigma^{1/2}(\bX^\top (\bX \bX^\top)^{-1}\bX - I_p)U\bbeta^*}_2 = \norm{\Sigma^{1/2}((\bX U)^\top ((\bX U) (\bX U)^\top)^{-1}(\bX U) - I_p)\bbeta^*}_2^2. 
 \end{align*} Moreover, $\bX U$ has the same distribution as $\bX$ and so $\norm{\Sigma^{1/2}(\bX^\top (\bX \bX^\top)^{-1}\bX - I_p)U\bbeta^*}_2$ has the same distribution as $\norm{\Sigma^{1/2}(\bX^\top (\bX \bX^\top)^{-1}\bX - I_p)\bbeta^*}_2$. In particular, they have the same expectation with respect to $\bX$. Now, let us consider the random matrix $U=\sum_{j=1}^p \varepsilon_j u_j u_j^\top$ where $(\epsilon_j)_{j=1}^p$ is a family of $p$ i.i.d. Rademacher variables (and $(u_j)_j$ is a basis of eigenvectors of $\Sigma$). Since $Uu_j = \varepsilon_j u_j$ for all $j$, $U$ commutes with $\Sigma$ and it is an orthogonal matrix. Therefore, we get
 \begin{equation}\label{eq:jaouad_random_trick}
   \bE\norm{\Sigma^{1/2}(\bX^\top (\bX \bX^\top)^{-1}\bX - I_p)\bbeta^*}_2^2 = \bE_U \bE_\bX \norm{\Sigma^{1/2}(\bX^\top (\bX \bX^\top)^{-1}\bX - I_p)U\bbeta^*}_2^2
  \end{equation} where $\bE_U$ (resp. $\bE_\bX$) denotes the expectation w.r.t. $U$ (resp. $\bX$). 

  We could now use Lemma~3 from \cite{TB21} to handle the Bayesian bias term from \eqref{eq:jaouad_random_trick} (that is the right hand side term). However, we want to avoid some (unnecessary) conditions required for that result to hold. To do so we use the DM theorem on $V_{k^*_b+1:p}$. 

  We denote for all $j\in[p], \bbeta^*_j=\inr{\bbeta^*, u_j}, z_j = \sigma_j^{-1/2}\bX u_j$, that is $z_j= \bG u_j$ where we recall that $\bX =\bG \Sigma^{1/2}$. We have 
  \begin{align*}
  &\bE_U\norm{\Sigma^{1/2}(\bX^\top (\bX \bX^\top)^{-1}\bX - I_p)U\bbeta^*}_2^2 = \sum_{j=1}^p (\bbeta_j^*)^2 \norm{\Sigma^{1/2}(\bX^\top (\bX \bX^\top)^{-1}\bX - I_p)u_j}_2^2\\ 
  &\geq \sum_{j=1}^p (\bbeta_j^*)^2 \inr{u_j,\Sigma^{1/2}(\bX^\top (\bX \bX^\top)^{-1}\bX - I_p)u_j}^2 =\sum_{j=1}^p (\bbeta_j^*)^2  \sigma_j \left(1-\sigma_j \norm{(\bX \bX^\top)^{-1/2}z_j}_2^2\right)^2.
  \end{align*}We lower bound the terms from the overfitting part (i.e. $j\geq k^*_b+1$) using the DM theorem and the one from the estimating part (i.e. $1\leq j\leq k^*_b$) using the strategy from \cite{TB21}. Let us start with the overfitting part. 

  Let $j\geq k^*_b+1$. On the event $\Omega^*$, it follows from Proposition~\ref{prop:DM_ellipsoid} that 
  \begin{equation*}
  \norm{(\bX \bX^\top)^{-1/2}z_j}_2\leq s_1[(\bX \bX^\top)^{-1/2}]\norm{z_j}_2 = \frac{\norm{z_j}_2}{\sqrt{s_N[\bX \bX^\top]}}\leq \frac{2\norm{z_j}_2}{\sqrt{\Tr(\Sigma_{k^*_b+1:p})}}
  \end{equation*}where we used that
  \begin{equation*}
      \bX\bX^\top = \bX_{1:k^*_b}\bX_{1:k^*_b}^\top + \bX_{k^*_b+1:p}\bX_{k^*_b+1:p}^\top \succeq  \bX_{k^*_b+1:p}\bX_{k^*_b+1:p}^\top \succeq (\Tr(\Sigma_{k^*_b+1:p})/4)I_N
  \end{equation*}
(thanks to Proposition~\ref{prop:DM_ellipsoid}). It follows from Borell's inequality that with probability at least $1-\exp(-N)$, $\norm{z_j}_2 = \norm{\bG u_j}_2\leq 3\sqrt{N}$. Hence, we obtain that with probability at least $1-\exp(-c_0N)$, $ \norm{(\bX \bX^\top)^{-1/2}z_j}_2\leq 6 \sqrt{N/\Tr(\Sigma_{k^*_b+1:p})}$ and so  
  \begin{equation*}
\left(1-\sigma_j \norm{(\bX \bX^\top)^{-1/2}z_j}_2^2\right)^2 \geq 1-2\sigma_j  \norm{(\bX \bX^\top)^{-1/2}z_j}_2^2\geq 1-\frac{ 12\sigma_{k^*_b+1} N}{\Tr(\Sigma_{k^*_b+1:p})}\geq 1-\frac{12}{b}\geq \frac{1}{2}
  \end{equation*}when $b\geq 24$. Therefore, when $N$ is larger than some absolute constant so that $1-\exp(-c_0N)\geq 1/2$, we get that 
  \begin{equation*}
  \begin{aligned}
  \bE_\bX \left[\sum_{j=k^*_b+1}^p (\bbeta_j^*)^2   \sigma_j \left(1-\sigma_j \norm{(\bX \bX^\top)^{-1/2}z_j}_2^2\right)^2\right]&\geq \frac{1}{4}\sum_{j=k^*_b+1}^p (\bbeta_j^*)^2   \sigma_j = \frac{\norm{\Sigma_{k^*_b+1:p}^{1/2}\bbeta^*_{k^*_b+1:p}}_2^2}{8}
  \end{aligned}
  \end{equation*}which is the expected lower bound on the overfitting part of the bias term. Let us now turn to the estimation part.

  Let $1\leq j\leq k^*_b$. To obtain the desired lower bound we use the same strategy as in Lemma~3 from \cite{TB21} together with Borell's inequality and the DM theorem. It follows from the Sherman-Morrison formulae (see the proof of Lemma~15 in \cite{TB21}) that 
  \begin{equation*}
   \sigma_j \left(1-\sigma_j \norm{(\bX \bX^\top)^{-1/2}z_j}_2^2\right)^2 = \frac{\sigma_j}{\left(1+\sigma_jz_j^\top A_{-j}^{-1}z_j\right)^2}
  \end{equation*}where $A_{-j} = \bX \bX^\top - \sigma_j z_j z_j^\top$. We have  $z_j^\top A_{-j}^{-1}z_j\leq \norm{z_j}_2 \norm{A_{-j}^{-1}z_j}_2\leq \norm{z_j}_2^2 s_1[A_{-j}^{-1}]=\norm{z_j}_2^2 / s_N[A_{-j}]$. Since $j\leq k^*_b$, we see that $A_{-j}\succeq X_{k^*_b+1:p} X_{k^*_b+1:p}^\top$ and so $s_N[A_{-j}]\geq s_N[X_{k^*_b+1:p} X_{k^*_b+1:p}^\top]$. On the event $\Omega^*$,
  \begin{equation*}
      s_N[X_{k^*_b+1:p} X_{k^*_b+1:p}^\top]\geq \Tr(\Sigma_{k^*_b+1:p})/4
  \end{equation*}
  (see Proposition~\ref{prop:DM_ellipsoid}).
The last result together with Borell's inequality yields with probability at least $1-\exp(-c_0 N)$, $z_j^\top A_{-j}^{-1}z_j\leq 36 N/ \Tr(\Sigma_{k^*_b+1:p})$. Furthermore, by definition of $k^*_b$ and since $j\leq k^*_b$, we have $bN\sigma_j> \Tr(\Sigma_{j:p})\geqslant \Tr(\Sigma_{k^*_b:p})$. Gathering the two pieces together we get that with probability at least $1-\exp(-c_0N)$, 
\begin{equation*}
\frac{\sigma_j}{\left(1+\sigma_jz_j^\top A_{-j}^{-1}z_j\right)^2}\geq \frac{\sigma_j}{\left(\frac{(36+b)\sigma_j N }{\Tr(\Sigma_{k^*_b+1:p})}\right)^2} = \frac{\Tr^2(\Sigma_{k^*_b+1:p})}{(36+b)^2 \sigma_j N^2}.
\end{equation*}Therefore, when $N$ is larger than some absolute constant so that $1-\exp(-c_0N)\geq 1/2$, we get that 
  \begin{equation*}
  \bE_\bX \left[\sum_{j=1}^{k^*_b} (\bbeta_j^*)^2   \sigma_j \left(1-\sigma_j \norm{(\bX \bX^\top)^{-1/2}z_j}_2^2\right)^2\right]\geq \frac{1}{2(36+b)^2}\sum_{j=1}^{k^*_b} \frac{(\bbeta_j^*)^2}{\sigma_j} \left(\frac{\Tr(\Sigma_{k^*_b+1:p})}{N}\right)^2 
  \end{equation*}which is the expected lower bound on the estimation of $\bbeta^*_{1:k^*_b}$ part of the bias term.




\section{Proof of Proposition~4: extension to the Dvoretsky-Milman's theorem to the anisotropic and heavy-tailed case} 
\label{sec:extension_to_the_heavy_tailed_case}

In this section, we extend the following stochastic properties originally satisfied by Gaussian vectors  to a heavy-tailed setup:
\begin{enumerate}
    \item For all $\vlambda\in\bR^N$, $c\sqrt{\Tr(\Sigma_{k+1:p})}\norm{\vlambda}_2 \leq \norm{\bX_{k+1:p}^\top\vlambda}_2\leq c'\sqrt{\Tr(\Sigma_{k+1:p})}\norm{\vlambda}_2$, when $N$ is smaller than the DM dimension of the ellipsoid $\Sigma_{k+1:p}^{-1/2}B_2$.
    \item For all $\vlambda\in\bR^N$, $\norm{\Sigma_{k+1:p}^{1/2}\bX_{k+1:p}^\top\vlambda}_2 \leq c^\prime\left(\sqrt{\Tr(\Sigma_{k+1:p}^2)} + \sqrt{N}\norm{\Sigma_{k+1:p}}_{op} \right)\norm{\vlambda}_2$.
\end{enumerate}

In our approach, these  properties are the geometrical properties the empirical covariance matrix needs to satisfy on the space $V_{k+1:p}$ spanned by the $p-k$ smallest eigenvectors of $\Sigma$ -- it is the part of the features space where interpolation of the noise happens. We provide here our main result on these two properties under weak moment assumptions. This result is a restatement of Proposition~4 in the main associated document.

\begin{Theorem}\label{theo:stochastic_argument_weak_moments} Let $X$ be a random vector in $\bR^{p}$ and let $X_1, \ldots, X_N$ be i.i.d. copies of $X$. We assume that the projection of $X$, denoted by $P_{k+1:p}X$, on the space $V_{k+1:p}$ spanned by the $p-k$ smallest singular vectors of the covariance matrix $\Sigma$ of $X$ is symmetric and satisfies:
    \begin{enumerate}
    \item[a)] $\bP\left\{\max_{1\leq i\leq N}\left|\frac{\norm{P_{k+1:p}X_i}_2^2}{\bE \norm{P_{k+1:p}X_i}_2^2}-1\right|\leq \delta\right\}\geq 1-\gamma$ for some constants $\delta<1/4$ and $\gamma<1/2$;
    \item[b)] there are constants $R,L,\alpha$ such that for all $\vv\in V_{k+1:p}$ and $2\leq q\leq R\log(e N)$, $\norm{\left<X,\vv\right>}_{L_q}\leq L q^{1/\alpha}\norm{\left<X,\vv\right>}_{L_2}$.  
\end{enumerate} 
Denote by $d^*\left(\Sigma_{k+1:p}^{-1/2}B_2\right)$ the Dvoretsky-Milman dimension of the ellipsoid $\Sigma^{-1/2}_{k+1:p}B_2$. There are absolute constants $c_1,c_2,c_3$ and $c_4$ such that the following holds. If $N\leq \kappa_{DM}d^*\left(\Sigma_{k+1:p}^{-1/2}B_2\right)$ for a small enough constant $\kappa_{DM}<1$ then with probability at least $1-2\gamma-N^{-2c_1}$ for every $\vlambda\in \bR^{N}$,
\begin{eqnarray*}
    (1-2\delta)\sqrt{\Tr(\Sigma_{k+1:p})}\norm{\vlambda}_2 &\leq& \norm{\bX_{k+1:p}^\top \vlambda}_2 \leq (1+2\delta)\sqrt{\Tr(\Sigma_{k+1:p})}\norm{\vlambda}_2,
\end{eqnarray*} and
\begin{equation*}
    \norm{\Sigma_{k+1:p}^{1/2}\bX_{k+1:p}^\top\vlambda}_2\leq c_3\left( \sqrt{\Tr(\Sigma_{k+1:p}^2)} +  \sqrt{N}\norm{\Sigma_{k+1:p}}_{op} \right)\norm{\vlambda}_2.
\end{equation*}
\end{Theorem}

We provide some examples for which Theorem~\ref{theo:stochastic_argument_weak_moments} holds in section~\ref{sec:example_weak_moments} under weak moments assumptions -- that is to show that \textit{a)} holds under weak moments assumptions. The proof of Theorem~\ref{theo:stochastic_argument_weak_moments} follows from Theorem~\ref{theo:DM_weak_moments} below applied to $P_{k+1:p}X$ and Theorem~\ref{theo:UB_DM_weak_moment} below for $P_{1:k}X$.

\subsection{An embedding Dvoretzky-Milman theorem for ellipsoid via anisotropic heavy-tailed random vectors}
In this section, a (probabilistic) Dvoretzky-Milman  embedding theorem is demonstrated for  subspaces produced by random matrices whose columns are independent anisotropic random vectors other than Gaussian ones. There has been many results proving random embeddings during the last few decades with classical results such as the Johnson-Lindenstrauss lemma and the Dvoretsky-Milman theorem (see \cite{vershynin_high-dimensional_2018}). So far random matrices have been key to successful embedding with optimal dimensions. Both random matrices with independent rows or independent columns have been used to achieve this goal. The setup of Theorem~\ref{theo:DM_weak_moments} below is the one of random matrices with independent columns and is the one associated with the Dvoretsky-Milman theorem.

We denote by $A:\bR^N\to \bR^p$ a random matrix with i.i.d. columns vectors distributed like $X$ and denoted by $X_1,\ldots, X_N$, i.e. $A=[X_1|\cdots| X_N]$ (for the sake of simplicity we still denote by $X_i$ the $i$-th column of $A$ even though we will actually use the projections of the $X_i$'s later for our application). We assume that $X$ is centered and denote by $\Gamma=\bE X X^\top \in\bR^{p\times p}$ the covariance matrix of $X$.  We also set $\ell^*=\sqrt{\bE\norm{X}_2^2} = \sqrt{\Tr(\Gamma)}$.  Our main result in this section is the following embedding theorem via matrix $A$ satisfied when $N$ is smaller than the Dvoretsky-Milman dimension $d_*(\Gamma^{-1/2}B_2^p)$ of the ellipsoid $\Gamma^{-1/2}B_2^p$.

\begin{Theorem}\label{theo:DM_weak_moments}Let $X$  be a random vector in $\bR^p$ and let $X_1, \ldots, X_N$ be $N$ i.i.d. copies of $X$. We assume that $X$ is centered,  symmetric and satisfies 
\begin{enumerate}
    \item[H1)] $\bP\left\{\max_{1\leq i\leq N}\left|\frac{\norm{X_i}_2^2}{\left(\ell^*\right)^2}-1\right|\leq \delta\right\}\geq 1-\gamma$ for some constants $\delta<1/2$ and $\gamma<1$.
    \item[H2)] For every $\vv\in\bR^{p}$, $\norm{\left<X,\vv\right>}_{L_q}\leq L q^{1/\alpha}\norm{\left<X,\vv\right>}_{L_2}$ for $2\leq q\leq R\log(e N)$ for some constants $R,L,\alpha$.
\end{enumerate} Then, there are absolute constants $\kappa_{DM}$ and $c_0$ depending only on $L$ and $\alpha$ such that, when $N\leq \kappa_{DM} d_*(\Gamma^{-1/2}B_2)$, with probability at least $1-\gamma - N^{-c_0}$,  for all $\vlambda\in\bR^N$, 
\begin{equation*}
(1-2\delta)\ell^*\norm{\vlambda}_2 \leq \norm{A \vlambda}_2 \leq (1+2\delta)\ell^*\norm{\vlambda}_2.
\end{equation*}
\end{Theorem}

An interesting feature of Theorem~\ref{theo:DM_weak_moments} is that even though the random vectors $X_1,\ldots, X_N$ used for the embedding may be far from Gaussian (for instance, we will provide examples where only $\log(N)$ moments suffice to satisfy the assumptions of Theorem~\ref{theo:DM_weak_moments}), the distortion and the maximal value of $N$ are the same (up to absolute constants) as the ones for Gaussian vectors. This was recently observed in \cite{bartl_random_2022} where a random embedding theorem similar to the one from Theorem~\ref{theo:DM_weak_moments} is obtained. 

We will follow the strategy from \cite{bartl_random_2022} to prove Theorem~\ref{theo:DM_weak_moments}. There are however some differences between \cite{bartl_random_2022}  and ours. The most important difference is that our random matrix $A$ is formed by anisotropic random column vectors, theirs is generated by isotropic ones --  this issue cannot simply be solved by a renormalization trick like changing the  set $T$ in \cite{bartl_random_2022} with $\Gamma^{1/2}T$. Our result shows that anisotropicity is not an issue regarding embeddings; one just need to adjust the distortion and the maximal value of $N$ as in the (anisotropic) Gaussian case. Moreover, because our indexing set is an Euclidean sphere (unlike in \cite{bartl_random_2022} where a general set $T$ is considered), our analysis relies on a one scale net approximation argument  instead of the multi-scaled generic chaining argument from \cite{bartl_random_2022}. This simplification allows us to remove extra logarithmic factors that were previously in \cite{bartl_random_2022}.

The rest of this section provides a proof of Theorem~\ref{theo:DM_weak_moments} along the same lines as in \cite{bartl_random_2022} adapted to the anisotropic case and the Euclidean sphere. Our aim is to show that with high probability 
\begin{equation}\label{eq:target_dvoretzky}
    \underset{\lambda\in S_2^{N-1}}{\mathrm{sup}}\left|\frac{\norm{A\lambda}_2^2}{\left(\ell^*\right)^2 } - 1 \right|\leq 2\delta.
\end{equation}The proof starts  with a ``column randomization'' approach that is analogous to the symmetrization method used in empirical process theory.

\subsubsection{Proof of Theorem~\ref{theo:DM_weak_moments}: a symmetrization and a one scale approximation arguments} Let $\epsilon = (\epsilon_i)_{1\leq i\leq N}$ be a Rademacher vector independent with $X_1, \ldots, X_N$ and set
\begin{equation*}
    A_\epsilon = A\left(\begin{matrix}
        \epsilon_1 & 0 & \cdots & 0\\
        0& \epsilon_2&  \cdots & 0 \\
        0 & 0 & \ddots & 0\\
        0 & 0  & \cdots& \epsilon_N
    \end{matrix} \right)=[\eps_1X_1|\cdots|\eps_N X_N]: \bR^N\to\bR^{p}.
\end{equation*}
Since we have assumed that the $X_i$'s are symmetric, $A_\epsilon$ has the same distribution as $A$ and so to get a high probability upper bound on \eqref{eq:target_dvoretzky}, we prove such a bound on \eqref{eq:target_dvoretzky} with $A_\eps$ replacing $A$. 

To obtain an high-probability upper bound for Equation \eqref{eq:target_dvoretzky} (with $A_\eps$ replacing $A$), we first use a one scale net argument. Set $V_{\eps_0} = \left\{\pi \lambda:\, \lambda\in S_2^{N-1} \right\}$ to be an $\epsilon_0$-net of $S_2^{N-1}$ so that for all $\lambda\in S_2^{N-1}$, $\norm{\pi \lambda-\lambda}_2\leq \epsilon_0$ and $\left|V_{\epsilon_0}\right|\leq(5/\epsilon_0)^N$ (see Lemma 4.4.1 of \cite{vershynin_high-dimensional_2018}). Parameter $0<\epsilon_0<1/2$ will be chosen at the very end of the proof to compensate all the absolute constants accumulating during  the proof.

It follows from Cauchy-Schwarz inequality that for every $\lambda\in  S_2^{N-1}$,
\begin{eqnarray*}
    \left|\frac{\norm{A_\epsilon \lambda}_2^2}{\left(\ell^*\right)^2} - 1 \right| &=& \left| \frac{\norm{A_\eps\left(\lambda-\pi\lambda\right)}_2^2}{\left(\ell^*\right)^2} + \frac{\norm{A_\eps\pi\lambda}_2^2}{\left(\ell^*\right)^2} + \frac{2\left<A_\eps \left(\lambda-\pi\lambda \right), A_\eps \pi  \lambda\right>}{\left(\ell^*\right)^2}- 1 \right|\\
    &\leq& \frac{\norm{A_\eps\left(\lambda-\pi\lambda\right)}_2^2}{\left(\ell^*\right)^2} + \left|\frac{\norm{A_\eps\pi\lambda}_2^2}{\left(\ell^*\right)^2} - 1\right| + \left|\frac{2\left<A_\eps \left(\lambda-\pi \lambda \right), A_\eps \pi  \lambda\right>}{\left(\ell^*\right)^2}\right|
\end{eqnarray*}and so
\begin{equation*}
    \underset{\lambda\in  S_2^{N-1}}{\mathrm{sup}}\left|\frac{\norm{A_\eps \lambda}_2^2}{\left(\ell^*\right)^2} - 1 \right| \leq \Phi^2 + \Psi^2 + 2\Phi\sqrt{\Psi^2 + 1}
\end{equation*}where
\begin{equation*}
    \Phi^2 := \underset{\lambda\in S_2^{N-1}}{\mathrm{sup}} \frac{\norm{A_\eps\left(\lambda-\pi\lambda\right)}_2^2}{\left(\ell^*\right)^2} \mbox{ and } \Psi^2 := \sup_{\lambda\in S_2^{N-1}}\left|\frac{\norm{A_\eps\pi\lambda}_2^2}{\left(\ell^*\right)^2} - 1\right|.
\end{equation*}

Thus, we only need to bound $\Phi$ and $\Psi$ from above. To that end we start with a decoupling argument to deal with the cross terms.

\subsubsection{Proof of Theorem~\ref{theo:DM_weak_moments}: a decoupling argument}
Let $(\eta_i)_{1\leq i\leq N}$ be $N$ i.i.d. selectors (i.e. $\bP(\eta_i=0) = \bP(\eta_i = 1) = 1/2$) and define $ I_\eta := \left\{1\leq i\leq N:\, \eta_i = 1 \right\}$. For all $\vu,\vv\in\bR^N$ and $I\subset[N]$, we also define
\begin{equation*}
    V_{I,\vu,\vv} := \frac{1}{\left(\ell^*\right)^2}\left<\left(\sum_{i\in I}\epsilon_i u_i X_i \right), \, \left(\sum_{j\in I^c}\epsilon_j v_j X_j \right)\right>.
\end{equation*}
For every $\vu\in B_2^{N}$, a decoupling technique (see for instance Chapter 6 of \cite{vershynin_high-dimensional_2018}) leads to the following result
\begin{eqnarray}\label{eq:decoupling}
    \nonumber\frac{\norm{A_\eps \vu}_2^2}{\left(\ell^*\right)^2} &=& \sum_{i=1}^N\frac{\norm{X_i}_2^2u_i^2}{\left(\ell^*\right)^2}  + \sum_{i\neq j}\epsilon_i\epsilon_j u_i u_j \frac{\left<X_i, X_j\right>}{\left(\ell^*\right)^2}=\sum_{i=1}^N\frac{\norm{X_i}_2^2u_i^2}{\left(\ell^*\right)^2} + \sum_{i,j=1}^N \bE_\eta 4\eta_i {(1-\eta_j)}\epsilon_i\epsilon_j u_i u_j \frac{\left<X_i, X_j\right>}{\left(\ell^*\right)^2}\\
    &=&\sum_{i=1}^N\frac{\norm{X_i}_2^2u_i^2}{\left(\ell^*\right)^2} + \frac{4}{\left(\ell^*\right)^2}\bE_{\eta}{\left<\left(\sum_{i\in I_\eta}\epsilon_i u_i X_i\right),\, \left(\sum_{j\in I_\eta^c}\epsilon_j u_j X_j \right) \right>}=\sum_{i=1}^N\frac{\norm{X_i}_2^2u_i^2}{\left(\ell^*\right)^2} +4\bE_{\eta}V_{I_\eta,\vu,\vu},
\end{eqnarray}
where $\bE_\eta$ is the expectation with respect to $(\eta_i)_{i\in[N]}$ conditionally on all other random variables. With the decoupling technique, we are going to estimate $\Phi^2$ and $\Psi^2$ from above. We start with $\Psi^2$ and get 
\begin{align}\label{eq:PSI_1}
     \nonumber \Psi^2 =&  \underset{\lambda\in S_2^{N-1}}{\mathrm{sup}}\left|\frac{\norm{A_\eps\pi   \lambda}_2^2}{\left(\ell^*\right)^2} - 1\right| \leq \underset{\lambda\in S_2^{N-1}}{\mathrm{sup}}\left|\sum_{i=1}^N\frac{\norm{X_i}_2^2}{\left(\ell^*\right)^2} \left(\pi   \lambda\right)_i^2 -1 \right|+ 4\underset{\lambda\in S_2^{N-1}}{\mathrm{sup}}\left|\bE_\eta V_{I_\eta,\pi   \lambda,\pi   \lambda}\right|\\
     &\leq \max_{1\leq i\leq N}\left|\frac{\norm{X_i}_2^2}{\left(\ell^*\right)^2} -1 \right|+ 4\underset{\lambda\in S_2^{N-1}}{\mathrm{sup}}\left|\bE_\eta V_{I_\eta,\pi   \lambda,\pi   \lambda}\right|.
\end{align}

Next, we have  $\Phi^2=\sup\left(\norm{A_\eps(\lambda-\pi\lambda)}_2^2/\left(\ell^*\right)^2:\, \lambda\in S_2^{N-1}\right)\leq (\epsilon_0/\ell^*)^2\norm{A_\eps}_{op}^2$, thus it remains to prove an high probability upper bound on $\norm{A_\eps}_{op}$. We do it by using a net-argument: let  $V_{1/2}=\left\{\pi\mu:\, \mu\in B_2^{N}\right\}\subset B_2^{N}$ be such that $\sup\left(\norm{\mu-\pi\mu}_2:\,\mu\in B_2^{N}\right)\leq 1/2$, and $\left|V_{1/2}\right|\leq 10^N$, we have $\norm{A_\eps}_{op} \leq 2\underset{\pi\mu\in V_{1/2}}{\mathrm{sup}}\norm{A_\eps(\pi\mu)}_2$ and, by the decoupling argument from \eqref{eq:decoupling}, we have
\begin{equation}\label{eq:PHI_1}
    \Phi^2 \leq \frac{4\epsilon_0^2}{\left(\ell^*\right)^2} \underset{\pi\mu\in V_{1/2}}{\mathrm{sup}}\norm{A_\eps(\pi\mu)}_2^2 
    \leq4\epsilon_0^2\left(\underset{1\leq i\leq N}{\mathrm{max}}\frac{\norm{X_i}_2^2}{\left(\ell^*\right)^2} + 4\underset{\pi\mu\in V_{1/2}}{\mathrm{sup}}\bE_\eta V_{I_\eta, \pi\mu, \pi\mu}\right).
\end{equation}
Therefore, we only need to find a high probability upper bound on $\bE_\eta V_{I_\eta,\vu,\vu}$ uniformly for all $\vu$ in $V_{1/2}$ and $V_{\eps_0}$ (the other two terms in \eqref{eq:PSI_1} and \eqref{eq:PHI_1} involving the $\norm{X_i}_2$'s will be handled using assumption \textit{H1)}).

\subsubsection{Proof of Theorem~\ref{theo:DM_weak_moments}: uniform bounds on $\bE_\eta V_{I_\eta,\vu,\vu}$ over one scale approximation sets.} We will obtain such a bound using Corollary 2.11 from \cite{bartl_random_2022} but we first introduce some notation: for all integer $s$ such that $1\leq 2^s\leq N$ we let
\begin{eqnarray*}
    \cO_{I,2^s} := \frac{1}{\left(\ell^*\right)^2}\mathrm{max}\left(\left<A \vx, A \vy \right>:\, \vx\in S_{I,2^s}, \, \vy\in S_{I^c, 2^s} \right) \mbox{ and }
    \cO_{2^s} := \bE_\eta \cO_{I_\eta, 2^s}
\end{eqnarray*}where for all $I\subset[N]$ and $\ell\in[N]$, we set
\begin{eqnarray*}
    S_I := \left\{\lambda\in S_2^{N-1}:\, \lambda_i=0,\forall i\in I^c \right\} \mbox{ and }
    S_{I,\ell} := \left\{\lambda\in S_I:\, \left|\left\{ i\in I:\, \lambda_i\neq 0  \right\}\right|\leq \ell \right\}.
\end{eqnarray*}

Note that the $\cO_{2^s}$'s do not depend on the Rademacher random variables $\epsilon_i$'s. The following lemma can be found in \cite{bartl_random_2022}.
\begin{Lemma}[Corollary 2.11,  \cite{bartl_random_2022}]\label{lemma:moments_bernoulli}
For  every  $a>0$ there  is a constant $C_1(a)$ such that the  following holds.  For all $\vu,\vv\in\bR^N$ and $2^s \geq \log{N}$, with  $\bP_\epsilon$-probability at least  $1-\exp\left(-a 2^s\right)$, we have
\begin{equation*}
    \bE_\eta\left|V_{I_\eta,\vu,\vv} \right| + \bE_\eta\left|V_{I_\eta^c,\vu,\vv} \right| \leq C_1(a)\norm{\vu}_2\norm{\vv}_2\cO_{2^s}
\end{equation*}where $\bP_\eps$ is the probability distribution with respect to $\eps$ given all the other random variables.
\end{Lemma}

Together with an union bound, we will obtain from  Lemma~\ref{lemma:moments_bernoulli}  high probability upper bounds on $\bE_\eta V_{I_\eta, \mu, \mu}$ uniformly over all $\mu\in V_\eta$ and $\mu\in V_{1/2}$. The latter will be used to upper bound  $\Phi^2$ and $\Psi^2$ thanks to \eqref{eq:PHI_1} and \eqref{eq:PSI_1}. We set $s = \lfloor \log_2{N} \rfloor$ and $N'=2^s$ (so that $(1/2)N\leq N' \leq N$). We apply Lemma~\ref{lemma:moments_bernoulli} with $a = a(\epsilon_0) = 6\log{(5/\epsilon_0)}$ then, with an union bound over all $\pi\lambda\in V_{\epsilon_0}$ and $\pi\mu\in V_{1/2}$ -- we note that $\left|V_{\epsilon_0} \right|\leq (5/\epsilon_0)^N$ and $|V_{1/2}|\leq 10^N\leq (5/\epsilon_0)^N$ because $\eps_0\leq 1/2$ -- with probability at least $1-\exp\left(-a(\epsilon_0)N/12 \right)$,
\begin{eqnarray}\label{eq:upper_Psi_square}
    \Psi^2=\underset{\lambda\in S_2^{N-1}}{\mathrm{sup}}\left|\frac{\norm{A_\eps\pi   \lambda}_2^2}{\left(\ell^*\right)^2} - 1\right| \leq\underset{1\leq i\leq N}{\mathrm{max}}\left|\frac{\norm{X_i}_2^2}{\left(\ell^*\right)^2} - 1 \right| + C_1(\eps_0)\cO_{N'}
\end{eqnarray}and  
\begin{eqnarray}\label{eq:upper_Phi_square}
    \Phi^2 = \underset{\lambda\in S_2^{N-1}}{\mathrm{sup}}\frac{\norm{A_\eps(\lambda-\pi\lambda)}_2^2}{\left(\ell^*\right)^2} 
    \leq4\epsilon_0^2\left(\underset{1\leq i\leq N}{\mathrm{max}}\frac{\norm{X_i}_2^2}{\left(\ell^*\right)^2} + C_1(\eps_0)\cO_{N'}\right). 
\end{eqnarray}
Hence, the next step to get high probability upper bounds on $\Phi^2$ and $\Psi^2$ is to obtain one on $\cO_{N'}$.

\subsubsection{Proof of Theorem~\ref{theo:DM_weak_moments}: upper bound on $\cO_{N'}$ via the iterative scheme from \cite{bartl_random_2022}} So far we have only used the symmetry assumption on $X$. To upper bound $\cO_{N'}$, in this  paragraph, we will use  assumptions H1) and H2).

To upper bound $\cO_{N'}$ we will use some intermediary quantities: for all integer $r$ such that $1\leq 2^r \leq N$ and $I\subset[N]$, we define
\begin{equation*}
    \cM_{I, 2^r} := \frac{1}{\ell^*}\max\left(\norm{\sum_{i\in I}\lambda_i\Gamma^{1/2} X_i}_2:\, \lambda\in S_{I, 2^r}\right) \mbox{ and } \cM_{2^r} = \cM_{[N], 2^r}.
\end{equation*}

We are now using the reduction scheme from \cite{bartl_random_2022} adapted to our  anisotropic setup  to control $\cM_{N'}$ via a series of $\cM_{I, 2^r}$.  
Our aim is to obtain an upper bound on $\cO_{N'}$ using $\cM_{N'}$. We will do it by proving that the following two inequalities hold with high probability:
\begin{equation}\label{eq:upper_bound_M}
    \left(\cM_{N'}\right)^2 \leq \norm{\Gamma}_{op}\left(\frac{\max_{1\leq i\leq N}\norm{X_i}_2^2}{\left(\ell^*\right)^2} + 4\cO_{N'} \right) 
\end{equation}and
\begin{equation}\label{eq:upper_bound_O}
    \cO_{N'}\leq c\frac{\sqrt{N}}{\ell^*}\cM_{N'} 
\end{equation}
where $c>0$ is some absolute constant. Indeed, it will follow from \eqref{eq:upper_bound_M} and \eqref{eq:upper_bound_O} that 
\begin{equation*}
\cO_{N'} \leq c\frac{\sqrt{N}}{\ell^*}\cM_{N'} 
\leq c \frac{\sqrt{N}}{\ell^*} \left[ \sqrt{\norm{\Gamma}_{op}}\frac{\max_{1\leq i\leq N}\norm{X_i}_2}{\ell^*} + \norm{\Gamma}_{op} \frac{\sqrt{N}}{\ell^*} \right]
\end{equation*}so that an upper bound on $\max_{1\leq i\leq N}\norm{X_i}_2$ and for $N$ small enough compare with the DM dimension of $\Gamma^{-1/2}B_2^p$ will provide a constant upper bound on $\cO_{N'}$.

 We start with a proof of \eqref{eq:upper_bound_M}. For all $\lambda\in S_2^{N-1}$, we have $\norm{\sum_{i\in [{N'}]}\lambda_i \Gamma^{1/2} X_i}_2\leq \sqrt{\norm{\Gamma}_{op}}\cdot\norm{\sum_{i\in [{N'}]}\lambda_i X_i}_2,$
then using the same machinery as above (i.e symmetrization, decoupling argument as in \eqref{eq:decoupling}, a one scale net approximation of $S_{[N],N'}$, a union bound and Lemma~\ref{lemma:moments_bernoulli}) we obtain that with probability at least $1-\exp(-cN)$ 
\begin{eqnarray*}
    \left(\cM_{N'}\right)^2 &\leq& \frac{\norm{\Gamma}_{op}}{\left(\ell^*\right)^2} \cdot{\mathrm{max}}\left(\sum_{i=1}^N\lambda_i^2\norm{X_i}_2^2 + \sum_{i\neq j}\lambda_i \lambda_j\left<X_i, X_j\right> :\, \lambda\in S_{[N],N'}\right)\\
    &\leq& \norm{\Gamma}_{op}\left(\frac{\max_{1\leq i\leq N}\norm{X_i}_2^2}{\left(\ell^*\right)^2} + 4\cO_{N'} \right).
\end{eqnarray*}

The next step is to show that \eqref{eq:upper_bound_O} holds with high probability. To that end we use the following lemma from \cite{bartl_random_2022}. We recall that for all $\vz\in\bR^N$, $z_1^*\geq \cdots \geq z_N^*\geq0$ denotes the non-increasing rearrangement of the absolute values of the coordinates of $\vz$.

\begin{Lemma}[\cite{bartl_random_2022}]\label{lemma:reduction}
There are absolute constants $C_3$ and  $C_4$ such that the following holds. Let $I\subset[N]$. For every integer $r$ such that $1\leq 2^r \leq N$ there exists a set $S_{I,2^r}' \subset S_{I, 2^r}$  such that $\log\left|S_{I, 2^r}'\right| \leq C_3 2^r \log\left(e N/2^r\right)$ for which the following statement holds:
\begin{equation}\label{eq:reduction}
    \cO_{I,2^s} \leq C_4 \left(\cO_{I, 1}+ \sum_{r=1}^s \left( \cE_{I, 2^r} + \cF_{I^c, 2^r} \right) \right)
\end{equation}where
\begin{eqnarray*}
    \cE_{I, 2^r} &:=& \frac{1}{\left(\ell^*\right)^2}\cdot{\mathrm{max}}\left(\sqrt{\sum_{i=2^{r-1}+1}^{2^r} \left(a_I(\vy)_i^*\right)^2 }:\, \vy\in S_{I^c, 2^r}' \right),\\
    \cF_{I^c, 2^r} &:=& \frac{1}{\left(\ell^*\right)^2}\cdot{\mathrm{max}}\left(\sqrt{\sum_{j=2^{r-1}+1}^{2^r} \left(b_I(\vx)_j^*\right)^2 }:\,\vx\in S_{I, 2^{r-1}}' \right),
\end{eqnarray*} and for all $\vx$ and $\vy$ in $\bR^N$, $a_I(\vy)$ and $b_I(\vx)$ are two vectors in $\bR^N$ such that $a_I(\vy)_i= \left<X_i,\, \sum_{j\in I^c}\vy_j X_j\right>$ if $i\in I$ and $0$ otherwise and $b_I(\vx)_j = \left<X_j,\, \sum_{i\in I}\vx_i X_i\right>$ if $j\in I^c$ and $0$ otherwise.

\end{Lemma}
Lemma~\ref{lemma:reduction} provides a reduction scheme from $\cO_{I,2^s}$ to $\cO_{I, 1}$. In the following, we will condition on $(X_j)_{j\in I^c}$ to deal with $\cE_{I, 2^r}$, and condition on $(X_i)_{i\in I}$ to deal with $\cF_{I^c, 2^r}$ alternatively. For the sake of completeness, we repeat the proof of Lemma~\ref{lemma:reduction} from \cite{bartl_random_2022} in the following.

\beginproof Let $c\geq 2$ be an absolute constant and set $\delta_{2^r}:= \left((2^r\wedge |I|)/e \left|I\right|\right)^{c}$ for all $r$  such that $1\leq 2^r\leq N$. Since there are less than $(2^r\wedge |I|) \log\left(e|I|/(2^r\wedge |I|) \right)$ subsets of $\left\{1,\cdots,\left|I\right|\right\}$ of cardinality at most $2^r$ and that there are $\delta_{2r}$-nets (w.r.t. the Euclidean norm) of the unit Euclidean sphere of $\bR^{I}$ with cardinality less than $(5/\delta_{2r})^{|I|}$ we conclude that there is  a set $ S_{I, 2^r}'\subset  S_{I,2^r}$ of cardinality at most
\begin{align*}
    \exp\left((2^r\wedge |I|) \log\left(\frac{e \left|I\right|}{(2^r\wedge |I|)  \delta_{2^r}}\right)\right) & = \exp\left((1+c)(2^r\wedge |I|) \log\left(\frac{e \left|I\right|}{(2^r\wedge |I|) }\right)\right)\leq \exp\left(C_3 2^r \log\left(\frac{e N}{2^r}\right)\right),
\end{align*}
such that for all $\vy\in\bR^N$, $\underset{\vx\in S_{I,2^r}}{\mathrm{max}}\left<\vx,\vy\right> \leq \frac{1}{1-\delta_{2^r}} \underset{\vx\in S_{I,2^r}'}{\mathrm{max}}\left<\vx,\vy\right>.$
The analogous statement holds also for $ S_{I^c, 2^r}$. 

By the definition of $\cO_{I, 2^s}$ and interchanging the order of two maxima, we have
\begin{eqnarray*}
    \cO_{I, 2^s} &=& \frac{1}{\left(\ell^*\right)^2}\cdot{\mathrm{max}}\left(\left<A \vx, A \vy \right>:\, \vx\in  S_{I,2^s}\, \vy\in  S_{I^c,2^s}\right)\leq\frac{1}{\left(\ell^*\right)^2}\underset{\vx\in  S_{I,2^s}}{\mathrm{max}}\frac{1}{1-\delta_{2^s}}\underset{\vy\in  S_{I^c,2^s}'}{\mathrm{max}}\left<A^\top A \vx, \vy \right>\\
    &=&\frac{1}{1-\delta_{2^s}}\frac{1}{\left(\ell^*\right)^2}\cdot{\mathrm{max}}\left( \max\left(\left< A \vx, A \vy\right>:\,\vx\in  S_{I,2^s}\right):\,\vy\in S_{I^c,2^s}' \right)\\
    &=&\frac{1}{1-\delta_{2^s}}\frac{1}{\left(\ell^*\right)^2}\cdot{\mathrm{max}}\left( \max\left( 
\sum_{i\in I}x_i\left< X_i, \sum_{j\in I^c} y_j X_j \right> :\,\vx\in  S_{I,2^s}\right):\,\vy\in S_{I^c,2^s}' \right)\\
    &=&\frac{1}{1-\delta_{2^s}}\frac{1}{\left(\ell^*\right)^2}\cdot{\mathrm{max}}\left(\sqrt{ \sum_{i=1}^{2^s}\left((a_I(\vy))_i^*\right)^2} : \, \vy\in  S_{I^c,2^s}'\right),
\end{eqnarray*}
where we have used the fact that for any $\vu =(u_i)_{i\in [N]}$ and $\ell\in[N]$
\begin{equation}\label{eq:maximum_rearrangement}
    {\mathrm{max}}\left(\sum_{i\in I}u_i v_i :\, \vv\in S_{I,\ell}\right) = \sqrt{\sum_{i=1}^\ell\left((\vu_{I})_i^*\right)^2},
\end{equation}where $\vu_{I}$ is the restriction to $I$  of $\vu$ and for all $i\in[N]$, $(\vu_{I})_i^*$ is the $i$-th largest absolute values of the coordinates of $\vu_{I}$ (in particular, $(\vu_{I})_i^* =0$ for all $i\geq |I|+1$).

Splitting the sum over $i=1,\cdots, 2^s$ into a sum over $\left\{1,\cdots, 2^{s-1}\right\}$ and over $\left\{2^{s-1}+1,\cdots, 2^s\right\}$, then
\begin{eqnarray*}
    \frac{1}{\left(\ell^*\right)^2}{\mathrm{max}}\left(\sqrt{ \sum_{i=1}^{2^s}\left((a_I(\vy))_i^*\right)^2} : \, \vy\in  S_{I^c,2^s}'\right)&\leq& \frac{1}{\left(\ell^*\right)^2}\max\left(\sqrt{\sum_{i=1}^{2^{s-1}}\left((a_I(\vy))_i^*\right)^2}:\, \vy\in S_{I^c, 2^r}'\right)+ \cE_{I, 2^s}.
\end{eqnarray*}
Using \eqref{eq:maximum_rearrangement} again and notice that $S_{I, 2^s}'\subset S_{I, 2^s}$, we obtain
\begin{equation*}
    \cO_{I, 2^s} \leq \frac{1}{1-\delta_{2^s}}\cdot{\mathrm{max}}\left(\frac{1}{\left(\ell^*\right)^2}\left<\sum_{i\in I}\vx_i X_i, \sum_{j\in I^c}\vy_j X_j\right> + \cE_{I, 2^s}:\, \vx\in  S_{I, 2^{s-1}},\, \vy\in  S_{I^c, 2^s}'\right).
\end{equation*}
By doing so, we reduce $\vx\in S_{I, 2^s}$ to $\vx\in S_{I,2^{s-1}}$. Repeat the above analysis with the roles of $\vx$ and $\vy$ reversed, we obtain
\begin{equation*}
    \frac{1}{\left(\ell^*\right)^2}{\mathrm{max}}\left(\left<\sum_{i\in I}\vx_i X_i, \sum_{j\in I^c}\vy_j X_j\right>:\, \vx\in  S_{I, 2^{s-1}},\,\vy\in S_{I^c, 2^s} \right) \leq \frac{1}{1-\delta_{2^{s-1}}}\left(\cO_{I, 2^{s-1}} + \cF_{I^c, 2^s}\right).
\end{equation*}
Combining the above two estimates, we get
\begin{equation*}
    \cO_{I, 2^s} \leq \frac{1}{1-\delta_{2^s}}\frac{1}{1-\delta_{2^{s-1}}}\left(\cO_{I, 2^{s-1}}+\cE_{I, 2^s} + \cF_{I^c, 2^s} \right).
\end{equation*}
Finally note that
\begin{equation*}
    \prod_{1\leq 2^r\leq \left|I\right|}\left(1-\delta_{2^r}\right)\left(1-\delta_{2^{r-1}}\right)\geq c''
\end{equation*}
for a constant $c''>0$ depending on $c$. Hence we can repeat the above procedure iteratively starting from $r=s$ until $r=1$.
{{\mbox{}\nolinebreak\hfill\rule{2mm}{2mm}\par\medbreak}}

Next it follows from Lemma~\ref{lemma:reduction} above that we need to handle $\ell_2$ norms of rearrangement of some random vectors indexed by some sets uniformly over these sets (see the definition of $\cE_{I, 2^r}$ and $\cF_{I, 2^r}$ in Lemma~\ref{lemma:reduction}). For that we use the following result from \cite{bartl_random_2022}.
\begin{Lemma}[Lemma 3.3, \cite{bartl_random_2022}]\label{lemma:uniform}
Let $C, C_5$ be some absolute constants, $C_3$ be the constant defined in Lemma~\ref{lemma:reduction} and $\alpha, L$ are the constants from our stochastic assumptions \textit{H1)} and \textit{H2)}. Then there  are  constants $R=R(C_3, C_5,  C)$ and $C_6  = C_6\left(C_3, C_5, \alpha,  C, L\right)$  such  that the following  holds. Let $1\leq 2^r \leq N$ and let $\cY$ be a collection of random variables such that $ \log \left|\cY\right| \leq C_3 2^r \log\left(e N/2^r\right)$ and for all $Y\in\cY$ and $2\leq q\leq R\log(e N)$,  $\norm{Y}_{L_q} \leq L q^{1/\alpha} \norm{Y}_{L_2}$. Then for every $u\geq e$, with probability at least $1-\exp\left(-C_5C\left(\log{u}\right)\cdot 2^r \log\left(e N/2^r\right)\right)$,
we have
\begin{equation*}
    {\mathrm{max}}\left(Y_{2^r}^*:\, Y\in \cY\right) \leq C_6 u  \log^{1/\alpha}\left(\frac{e N}{2^r}\right) \cdot {\mathrm{max}}\left(\norm{Y}_{L_2}:\, Y\in\cY\right),
\end{equation*}
where for all $Y\in \cY$ ,$(Y_i)_{i=1}^N$ are  $N$  independent copies of $Y$ and  $(Y_i^*)$ is the $i$-th largest element in $\{|Y_1|, \ldots, |Y_N|\}$.
\end{Lemma}

Recall the $\cO_{N'} = \bE_\eta\cO_{I_\eta, 2^s}$ for $s=\lfloor \log_2{N} \rfloor$. Let $I = I_\eta$ in Lemma~\ref{lemma:reduction}, and integrate \eqref{eq:reduction} with respect to $\eta$, we obtain $\cO_{N'} \leq C_4 \left(\cO_{1}+ \sum_{r=1}^{\lfloor \log_2{N} \rfloor} \left( \bE_{\eta}\cE_{I_\eta, 2^r} + \bE_{\eta}\cF_{I_\eta^c, 2^r} \right) \right).$
For a fixed $I\subset[N]$, condition on $(X_j)_{j\in I^c}$, we have for all integers $1\leq r \leq \lfloor \log_2{N} \rfloor$,
\begin{align*}
    \cE_{I, 2^r} & = \frac{1}{\left(\ell^*\right)^2}\cdot{\mathrm{max}}\left(\sqrt{\sum_{i=2^{r-1}+1}^{2^r} \left(\left<X_i, \sum_{j\in I^c}\vy_j X_j\right>^*\right)^2 }:\, \vy\in S_{I^c, 2^r}' \right)\\ 
    &\leq \frac{\sqrt{2^r}}{\left(\ell^*\right)^2}{\mathrm{max}}\left(\left<X_i, \sum_{j\in I^c}\vy_j X_j\right>_{2^{r-1}}^*:\, \vy\in S_{I^c, 2^r}'\right).
\end{align*} Consider the family of random variables $\cY := \left\{\frac{1}{\ell^*}\left<X, \sum_{j\in I^c}\vy_j X_j\right>:\, \vy\in  S_{I^c,2^r}'\right\}$.
The cardinality of $\cY$ satisfies the first assumption in Lemma~\ref{lemma:uniform}, and for every $Y\in\cY$ and every $2\leq q\leq R\log(e N)$,
\begin{equation*}
    \norm{Y}_{L_q} = \frac{1}{\ell^*}\norm{\left<X,\sum_{j\in I^c}\vy_j X_j\right>}_{L_q} \leq \frac{L q^{1/\alpha}}{\ell^*}\norm{\sum_{j\in I^c}\vy_j\Gamma^{1/2} X_j}_{2} = L q^{1/\alpha}\norm{Y}_{L_2} \leq L q^{1/\alpha}\cM_{I^c, 2^r},
\end{equation*}
because of the moment assumption on the linear projections of $X$ and because $\vy\in S_{I^c, 2^r}'\subset S_{I^c, 2^r}$ (here, $L_q$ norms are taken w.r.t. $X$ conditionally on $(X_j)_{j\in I^c}$).  Applying Lemma~\ref{lemma:uniform}, for every $u\geq e$, condition on $(X_j)_{j\in I^c}$, take the product probability measure of $(X_i)_{i\in I}$(denoted as $\bP_{\bX_I}$), we get that with $\bP_{\bX_I}$-probability at least $1-\exp\left(-C_5 C \left( \log{u}\right)2^r \log\left(e N/2^r\right) \right)$,
\begin{equation*}
    \cE_{I, 2^r}\leq \frac{\sqrt{2^r}}{\ell^*}\max\left(Y^*_{2^{r-1}+1}:Y\in\cY\right)\leq C_6 u \frac{\sqrt{2^r}}{\ell^*}\log^{1/\alpha}\left(\frac{e N}{2^r}\right)\cM_{I^c, 2^r} 
\end{equation*}where we used that $\max(\norm{Y}_{L_2}:Y\in\cY)\leq \cM_{I^c, 2^r}$. Moreover, $\cM_{I^c, 2^r}\leq \cM_{2^r}$ and, by Fubini theorem we can replace $\bP_{\bX_I}$ by $\bP_{\bX}$, the product probability measure of $(X_i)_{i\in[N]}$, so that there exists absolute constants $C_7, C_8>0$ for which (using an union bound over all $1\leq r\leq  \lfloor \log_2{N}\rfloor$), with probability at least
\begin{equation*}
    1-\sum_{r=1}^{\lfloor \log_2{N}\rfloor} \exp\left( -C_5 C \left( \log{u}\right)2^r \log\left(\frac{e N}{2^r}\right)\right) \geq {1-u^{-C_7C\log(eN)}},
\end{equation*}
we have
\begin{equation*}
    \sum_{r=1}^{\lfloor \log_2{N}\rfloor} \cE_{I, 2^r} \leq C_6 u \sum_{r=1}^{\lfloor \log_2{N}\rfloor} \frac{\sqrt{2^r}}{\ell^*}\log^{1/\alpha}\left(\frac{e N}{2^r} \right) \cM_{2^r} \leq {C_8 u \frac{\sqrt{N'}}{\ell^*}\cM_{N'}}.
\end{equation*}
As our final goal is to obtain an high-probability upper bound for $\bE_\eta \cE_{I_\eta,2^r}$, we need to exchange the order of two integrals. Choose $C_7$ large enough, and a tail integration leads to
\begin{equation}\label{eq:tail_integral}
    \left(\bE_\bX \left|\frac{\sum_{r=1}^{\lfloor \log_2{N}\rfloor} \cE_{I, 2^r}}{C_8 \frac{\sqrt{N'}}{\ell^*} \cM_{N'}} \right|^q\right)^{1/q} \leq C_9,
\end{equation}
for $q = 2C\log(eN)$ and an absolute constant $C_9>0$. Since \eqref{eq:tail_integral} holds for any $I\subset\left\{1, \cdots, N\right\}$, we use it for $I_\eta$ and integrate over $\eta$, then, by Jensen's inequality applied to $\bE_\eta$ and the convex function $t\mapsto \left|t\right|^{q}$, we get the moment bound
\begin{equation*}
    \left(\bE_\bX \left|\frac{\sum_{r=1}^{\lfloor \log_2{N}\rfloor}\bE_\eta \cE_{I_\eta, 2^r}}{C_8 \frac{\sqrt{N'}}{\ell^*} \cM_{N'}} \right|^q\right)^{1/q} \leq \bE_\eta\norm{\frac{\sum_{r=1}^{\lfloor \log_2{N}\rfloor} \cE_{I_\eta, 2^r}}{C_8 \frac{\sqrt{N'}}{\ell^*} \cM_{N'}}}_{L_q(\bP_\bX)} \leq C_9
\end{equation*}where the $L_q(\bP_\bX)$ norm is the one taken with respect to $\bP_\bX$. Therefore, by Markov's inequality,
\begin{equation*}
    \bP_\bX\left(\sum_{r=1}^{\lfloor \log_2{N}\rfloor} \bE_\eta\cE_{I_\eta, 2^r}\geq C_8 C_9 e \frac{\sqrt{N'}}{\ell^*} \cM_{N'}\right) \leq \exp(-q) = N^{-2C}.
\end{equation*}
The analysis for $\cF_{I^c, 2^r}$ and {$\cO_{I, 1}$} follow by using the same argument and so we can  state the following lemma.
\begin{Lemma}\label{lem:UB_on_ON}
There is an absolute constant $C_{10}$ such that, with $\bP_\bX$-probability  at least $1-N^{-2C}$, $\cO_{N'} \leq C_{10}\frac{\sqrt{N'}}{\ell^*} \cM_{N'}$.
\end{Lemma}

\subsubsection{End of the proof of Theorem~\ref{theo:DM_weak_moments}}
It follows from Lemma~\ref{lem:UB_on_ON} that \eqref{eq:upper_bound_M} and \eqref{eq:upper_bound_O} hold with probability at least $1-N^{-2C}$. Moreover, $\max_{1\leq i\leq N}\norm{X_i}_2\leq \ell^* \sqrt{1+\delta}$ with probability at least $1-\gamma$,  so, with probability at least $1-\gamma-N^{-2C}$,
\begin{equation}\label{eq:upper_cO}
    \cO_{N'} \leq \left(C_{10}\sqrt{1+\delta}\frac{\sqrt{N'\norm{\Gamma}_{op}}}{\ell^*}\vee C_{10}^2\frac{N'\norm{\Gamma}_{op}}{\left(\ell^*\right)^2}\right).
\end{equation}Moreover, we assumed that $ N \leq \kappa_{DM} d^*\left(\Gamma^{-1/2}B_2\right)$ for some constant $\kappa_{DM}<1$ small enough and $d^*\left(\Gamma^{-1/2}B_2\right)\sim (\ell^*)^2/\norm{\Gamma}_{op}$, hence, we can bound $\cO_{ N'}$ from above by some absolute constant $C_{11} = C_{11}(\kappa_{DM})<C_{1}(\epsilon_0)^{-1}$(by choosing a small enough $\kappa_{DM}$ and recall that $C_1(\epsilon_0)$ depends on $\epsilon_0$ from Lemma~\ref{lemma:moments_bernoulli}, see Equation \eqref{eq:upper_Psi_square} and Equation \eqref{eq:upper_Phi_square}). Finally, we plug this high probability upper bound by a constant on $\cO_{N'}$  back in \eqref{eq:upper_Psi_square} and \eqref{eq:upper_Phi_square} to get 
\begin{equation*}
    \Phi^2 \leq 4\epsilon_0^2\left((1+\delta)+C_{1}(\epsilon_0)C_{11}\right)\mbox{ and } \Psi^2 \leq\delta +C_{1}(\epsilon_0)C_{11}.
\end{equation*}
As a result, with probability at least $1-\gamma-N^{-2C}$,
\begin{equation*}
    \underset{\lambda\in  S_2^{N-1}}{\mathrm{sup}}\left|\frac{\norm{A_\eps \lambda}_2^2}{\left(\ell^*\right)^2} - 1 \right| \leq 4\epsilon_0^2\left(1+\delta+C_{1}(\epsilon_0)C_{11}\right) + \delta + C_{1}(\epsilon_0)C_{11} + 4\epsilon_0\left(1+\delta+C_{1}(\epsilon_0)C_{11}\right).
\end{equation*}
Note that  $\delta<1/2$ and we can choose $\epsilon_0$ and $C_{11}$ small enough (this will only affects the constant in the deviation parameter and $\kappa_{DM}$) so that the right hand side term in the inequality above is strictly smaller than $2\delta<1$. This proves that with probability at least $1-\gamma-N^{-2C}$ for all $\vlambda\in\bR^N$, 
\begin{equation*}
(1-2\delta)\ell^*\norm{\vlambda}_2 \leq \norm{A_\eps \vlambda}_2 \leq (1+2\delta)\ell^*\norm{\vlambda}_2,
\end{equation*}
and then we conclude since $A_\eps$ and $A$ have the same distribution by symmetry of the $X_i$'s.
{{\mbox{}\nolinebreak\hfill\rule{2mm}{2mm}\par\medbreak}}

\subsection{An upper bound result for point~\textit{2.}} A close inspection of the proof of Theorem~\ref{theo:DM_weak_moments} reveals that we can get a result like the one we want in point \textit{2.}, i.e. an uniform over all $\blambda\in\cS_2^{N-1}$  upper bound on $\norm{\Gamma^{1/2}A\blambda}_2$ under no assumption on $N$. Let us state this result now.

\begin{Theorem}\label{theo:UB_DM_weak_moment}
We assume that $X_1,\ldots, X_N$ are $N$ i.i.d. copies of $X$, a random vector in $\bR^p$ satisfying the assumptions from Theorem~\ref{theo:DM_weak_moments} (i.e. centered, symmetric, \textit{H1)} and \textit{H2)}). There  are absolute constants $c_0$ and $c_1$ depending only on $L$ and $\alpha$ such that  with probability at least $1-\gamma-N^{-2c_0}$, for all $\vlambda\in \bR^N$,
\begin{equation*}
    \norm{\Gamma^{1/2}A\vlambda}_2 \leq c_1 \left( \sqrt{\Tr(\Gamma^2)} +  \sqrt{N}\norm{\Gamma}_{op}\right)\norm{\vlambda}_2. 
\end{equation*}
\end{Theorem}

\beginproof
We apply the  result above for $\bY = \Gamma^{1/2}A$, without the Dvoretzky-Milman assumption, that is, when we do not necessarily assume that $N<\kappa_{DM}d^*\left(\Gamma^{-1}B_2\right)$, then \eqref{eq:upper_cO} can be bounded from above by
\begin{equation*}
    \cO_{N'} \leq C_{10}\sqrt{1+\delta}\frac{\sqrt{N\norm{\Gamma^2}_{op}}}{\sqrt{\Tr\left(\Gamma^2\right)}} \vee  C_{10}^2\frac{N\norm{\Gamma^2}_{op}}{\Tr\left(\Gamma^2\right)}.
\end{equation*}
Using this bound in \eqref{eq:upper_Psi_square} and \eqref{eq:upper_Phi_square} together with Assumption~\textit{H1)}, we get
\begin{eqnarray*}
    \Phi^2 &\leq& 4\epsilon_0^2 \left(1+\delta+C_{1}(\epsilon_0)C_{10}\sqrt{1+\delta}\frac{\sqrt{N}\norm{\Gamma}_{op}}{\sqrt{\Tr\left(\Gamma^2\right)}} \vee  C_{1}(\epsilon_0)C_{10}^2\frac{N\norm{\Gamma^2}_{op}}{\Tr\left(\Gamma^2\right)}\right)\\
    &\leq&4\epsilon_0^2 \left(2+  \left(\frac{\sqrt 2C_{1}(\epsilon_0)C_{10}\sqrt{N}\norm{\Gamma}_{op}}{\sqrt{\Tr\left(\Gamma^2\right)}}\right) \vee \left( C_{1}(\epsilon_0)C_{10}^2\frac{ N\norm{\Gamma}_{op}^2}{{\Tr\left(\Gamma^2\right)}}\right)\right)\\
\end{eqnarray*}
and
\begin{eqnarray*}
    \Psi^2 &\leq& \delta + C_{1}(\epsilon_0)C_{10}\sqrt{1+\delta}\frac{\sqrt{N}\norm{\Gamma}_{op}}{\sqrt{\Tr\left(\Gamma^2\right)}} \vee  C_{1}(\epsilon_0)C_{10}^2\frac{N\norm{\Gamma}_{op}^2}{\Tr\left(\Gamma^2\right)}\leq \delta + C_{1}(\epsilon_0)C_{10}^2\frac{N\norm{\Gamma}_{op}^2}{\Tr\left(\Gamma^2\right)} \vee \sqrt 2C_{1}(\epsilon_0)C_{10}\frac{\sqrt{N}\norm{\Gamma}_{op}}{\sqrt{\Tr\left(\Gamma^2\right)}}.
\end{eqnarray*}
As a result, setting $\epsilon_0=1/2$ and using that $\delta\leq1$, with probability at least $1-\gamma-N^{-2c_0}$, there exists an absolute constant $c'$ such that for all $\vlambda\in \bR^N$,
\begin{eqnarray*}
    \norm{\Gamma^{1/2}A_\eps\vlambda}_2^2 &\leq& c'\left( {\Tr(\Gamma^2)} +  {N}\norm{\Gamma}_{op}^2 +  \norm{\Gamma}_{op}{\sqrt{N\Tr(\Gamma^2)}} \right)\norm{\vlambda}_2^2\leq c''\left( \sqrt{\Tr(\Gamma^2)} +  \sqrt{N}\norm{\Gamma}_{op}\right)^2\norm{\vlambda}_2^2 
\end{eqnarray*}
This implies our goal since $A$ and $A_\eps$ have the same distribution.
{{\mbox{}\nolinebreak\hfill\rule{2mm}{2mm}\par\medbreak}}

\section{Proof of Proposition~5: an isomorphic property under weak moment assumptions}\label{sec:isomorphy_weak_moments}
To establish the benign overfitting property for the minimum $\ell_2$-norm interpolant estimator, we need to handle the behavior of the design matrix on the eigenspace $V_{1:k}$ spanned by the largest $k$ eigenvectors of $\Sigma$. In the case $k\lesssim N$, we know that it is enough to show that $\bX_{1:k}$ satisfies an isomorphic property on the entire spare $V_{1:k}$:
\begin{itemize}
    \item[3.] For all $\vv\in V_{1:k}$, $(1/2)\norm{\Sigma_{1:k}^{1/2}\vv}_2^2 \leq (1/N)\norm{\bX_{1:k}\vv}_2^2  \leq (3/2)\norm{\Sigma_{1:k}^{1/2}\vv}_2^2$.
\end{itemize} In the Gaussian case this follows from a straightforward one scale net argument, however in the heavy-tailed case, the task is much more complicated (because concentration does not balance the covering complexity of the Euclidean sphere $\cS_2^{k-1}$). We need to use a very different machinery to prove isomorphy under weak moment assumption. Fortunately,  \cite{tikhomirov_sample_2018} provides such as result under (very) weak moments assumptions. We will therefore apply this result for our statistical purpose i.e. to prove Proposition~5 from the main document. Let us first recall Proposition~5 from the main document.

\begin{Theorem}\label{theo:proposition_5_main}
Assume that the projection of $X$ on the space $V_{1:k}$ (spanned by the $k$ largest singular vectors of the covariance matrix $\Sigma$ of $X$), denoted by $P_{1:k}X$,  satisfies: 
    \begin{enumerate}
        \item[c)] there exists $q>4$, $B\geq1$ such that $\norm{\left<P_{1:k}X,\vv\right>}_{L_q} \leq B^{1/q}\norm{\left<P_{1:k}X,\vv\right>}_{L_2}$, for all $\vv\in V_{1:k}$,
        \item[d)] there are constants $c_0$ and $\gamma'<1$ such that with probability at least $1-\gamma'$, for all $i\in[N]$, $\norm{\Sigma_{1:k}^{-1/2}X_i}_2\leq c_0 \sqrt{k}$.
    \end{enumerate}
    Assume that $C k\leq  N$ where $C$ is some absolute constant depending only on $B,q$ and $c_0$ large enough, then with probability at least $1-\gamma'-1/k$, for all $\vv\in V_{1:k}$,
    \begin{equation}\label{eq:objective_rip}
        \frac{1}{2}\norm{\Sigma_{1:k}^{1/2}\vv}_2^2 \leq \frac{1}{N}\norm{\bX_{1:k}\vv}_2^2  \leq \frac{3}{2}\norm{\Sigma_{1:k}^{1/2}\vv}_2^2.
    \end{equation}
\end{Theorem}

We want to show that \eqref{eq:objective_rip} hold with large probability. By homogenity, it is equivalent to prove \eqref{eq:objective_rip} for all $\vv\in \Sigma_{1:k}^{-1/2} S_2^{k-1}$ i.e. to show that $1/2\leq (1/N)\sum_{i=1}^N \inr{Z_i,\vu}^2\leq 3/2$ for all $\vu\in S_2^{k-1}$ where $Z_i=\Sigma_{1:k}^{-1/2}X_i$ is isotropic. Hence and unlike in the last section, the anisotropic case can be handled the same way as in the isotropic case. The difference with Section~\ref{sec:extension_to_the_heavy_tailed_case} is that here we work with $\bX_{1:k}$, a random matrix with independent rows, whereas in the previous section, we worked with $\bX_{k+1:p}^\top$, a random matrix with independent columns.

We set $\hat\Sigma_N = (1/N)\sum_{i=1}^N Z_i Z_i^\top$ the empirical covariance matrix of the $Z_i$'s then \eqref{eq:objective_rip} is equivalent to show that with high probability, $\norm{\hat\Sigma_N - I_k}_{op}^2  = \underset{\vu\in S_2^{k-1}}{\mathrm{sup}}\left|\frac{1}{N}\sum_{i=1}^N\left<Z_i, \vu\right>^2-1 \right| \leq \frac{1}{2}$.
It is therefore enough to obtain an estimation result, with respect to the operator norm, on the sample covariance matrix in the isotropic case under weak moment assumptions. The following result from \cite{tikhomirov_sample_2018} is the state-of-the-art conclusion for such a problem.

\begin{Theorem}[\cite{tikhomirov_sample_2018}]\label{theo:tikhomirov_isomorphy} 
    There is a non-increasing function $\nu:(2,\infty)\to\bR_+$ with the following property. Let $q>2$, $B\geq 1$, and assume that $N\geq 2k$.  Let $Z$ be a centered isotropic $k$-dimensional random vector. Suppose we have for all $\vv\in\bR^k$,
    \begin{equation}\label{eq:norm_equivalence_rip}
        \norm{\left<Z,\vv\right>}_{L_q} \leq B^{1/q}\norm{\left<Z,\vv\right>}_{L_2}.
    \end{equation}
    Let $(Z_i)_{i\in[N]}$ be independent copies of $Z$. Then with probability at least $1-1/k$,
    \begin{equation*}
        \nu(q)^{-1}\norm{\hat\Sigma_{N} - I_k }_{op} \leq \frac{1}{N}\underset{i\leq N}{\mathrm{max}}\norm{Z_i}_2^2 + B^{2/q}\left(\frac{k}{N} \right)^{\frac{q-2}{q}}\log^4\left(\frac{N}{k}\right) + B^{2/q}\left(\frac{k}{N}\right)^{\frac{(q\wedge 4)-2}{q\wedge 4}}.
    \end{equation*}
\end{Theorem}
In particular, it follows from Theorem~\ref{theo:tikhomirov_isomorphy} that when $q>4$ and under the same assumptions as  Theorem~\ref{theo:tikhomirov_isomorphy},  with probability at least $1-1/k$,
\begin{equation*}
\norm{\hat\Sigma_{N} - I_k }_{op} \lesssim \nu(q)\left( \frac{1}{N}\underset{i\leq N}{\mathrm{max}}\norm{Z_i}_2^2 + B^{2/q}\sqrt{\frac{k}{N}} \right).
\end{equation*}
As a result, if we know that with probability $1-\gamma'$, for all $i\in[N]$ $\norm{Z_i}_2 \leq c_0 \sqrt{k}$ and $C k\leq N$ for $C> 1$ an absolute constant large enough, then we get the desired isomorphic result on $\bR^k$. Let us now state the result we will use for the projection of the input vector on $V_{1:k}$ (in that case,  $\Gamma$ is the covariance matrix $\Sigma_{1:k}$). It is a restatement of Proposition~5 in the main document where $X$ is replaced by $P_{1:k}X$ and so $\Gamma$ is $\Sigma_{1:k}$.

\begin{Theorem}\label{theo:tikhomirov_isomorphy_2} Let $X$ be a centered random vector in $\bR^k$ and let $X_1, \ldots, X_N$ be $N$ i.i.d. copies of $X$. Denote by $\Gamma$ the covariance matrix of $X$. We assume that there exists $q>4$, $B\geq1$ and $c_0$ such that
    \begin{enumerate}
        \item[H3)] $\norm{\left<X,\vv\right>}_{L_q} \leq B^{1/q}\norm{\left<X,\vv\right>}_{L_2}$, for all $\vv\in\bR^k$,
        \item[H4)] with probability at least $1-\gamma'$, for all $i\in[N]$, $\norm{\Gamma^{-1/2}X_i}_2\leq c_0 \sqrt{k}$.
    \end{enumerate}
    Assume that $C k\leq  N$ where  $C$ is an absolute constant large enough (depending only on $B,q$ and $c_0$), then with probability at least $1-\gamma'-1/k$, for all $\vv\in\bR^k$,
    \begin{equation*}
        \frac{1}{2}\norm{\Gamma^{1/2}\vv}_2^2 \leq \frac{1}{N}\norm{\bX_{1:k}\vv}_2^2  \leq \frac{3}{2}\norm{\Gamma^{1/2}\vv}_2^2.
    \end{equation*}
\end{Theorem}

\section{Some examples where Theorem~\ref{theo:stochastic_argument_weak_moments} and Theorem~\ref{theo:proposition_5_main} (i.e. Proposition~4 and 5 in the main document) hold under weak moments assumptions}\label{sec:example_weak_moments}
In this section, we provide several examples for which Theorem~\ref{theo:stochastic_argument_weak_moments} and \ref{theo:proposition_5_main}  hold. The first example is the Gaussian case (i.e. when $X$ is a Gaussian vector) which served us as an expository case to show how a geometrical viewpoint can be used to analyze the minimum $\ell_2$-norm interpolant estimator and to identify necessary and sufficient conditions for its benign overfitting property to hold. In the Gaussian setup assumptions \textit{a)} and \textit{d)} (from Theorem~\ref{theo:stochastic_argument_weak_moments} and Theorem~\ref{theo:proposition_5_main}) follow from the Borell-TIS inequality and assumptions \textit{b)} and \textit{c)} hold because the $p$-th moment of a Gaussian variable $\cN(0,1)$ is of the order of $\sqrt{p}$, for all $p\geq1$.

The second example, analyzed in the next section, is from the sub-Gaussian ensemble and relies on Hanson-Wright inequality. The third example is studied in the last sub-section. It deals with a heavy-tailed case when the marginal of $X$ satisfy a $L_{R\log{(eN)}}-L_2$ norms-equivalence condition. Such an assumption is standard in empirical process theory, compressed sensing and robust statistics \cite{mendelson_multiplier_2017,lecue_sparse_2017,chinot_robust_2021}. These two sections assume that $X = \Sigma^{1/2}Z$ for some random vector $Z$ with independent coordinates.

\subsection{Symmetric sub-Gaussian random vectors} A random vector $X$ with covariance matrix $\Gamma$ is said to be sub-Gaussian, if there is a constant $\kappa\geq 1$ such that for all vector $\vv$, $\norm{\left<X,\vv\right>}_{\Psi_2}\leq \kappa\norm{\Gamma^{1/2}\vv}_2$, where $\norm{\cdot}_{\Psi_2}$ is the $\Psi_2$ Orlicz  norm (see, for example, section~2.5 of \cite{vershynin_high-dimensional_2018}). This is equivalent to have a moment growth condition as in \textit{b)} and \textit{c)} for $\alpha=2$ but for all $q\geq2$ (see \cite[Proposition 2.5.2]{vershynin_high-dimensional_2018}). Hence, under a sub-gaussian assumption \textit{b)} and \textit{c)} of Theorem~\ref{theo:stochastic_argument_weak_moments} and \ref{theo:proposition_5_main}  hold and so it only remains to check assumptions \textit{a)} and \textit{d)} of these theorems. To that end we use an extension of Hanson-Wright inequality to the subgaussian case that we recall now.

\begin{Theorem}\label{theo:hanson_wright_ineq_sub-gauss}\cite[Section 6.3]{vershynin_high-dimensional_2018}
There are absolute constants $c_0$ and $c_1$ such that the following holds. Let $Z=(z_j)_j$ be a random vector with independent centered, variance one  coordinates and $\Gamma$ be a PSD matrix. With probability at least
\begin{equation}\label{eq:proba_theo_HS_sub-gauss}
1-c_0\exp\left(-c_1 d_*(\Gamma^{-1/2}B_2)\right),
\end{equation}
\begin{equation*}
\left| \norm{\Gamma^{1/2}Z}_2^2 - \bE \norm{\Gamma^{1/2}Z}_2^2 \right|\leq \frac{1}{4}\bE\norm{\Gamma^{1/2}Z}_2^2.
\end{equation*}
 \end{Theorem} 

It follows from Theorem~\ref{theo:hanson_wright_ineq_sub-gauss} (applied to $\Gamma = \Sigma_{k+1:p}$) and an union bound that if $\log N\leq c^\prime d_*(\Sigma_{k+1:p}^{-1/2}B_2)$ where $c^\prime$ is some absolute constant then, with the same probability as in  \eqref{eq:proba_theo_HS_sub-gauss}, for all $i\in[N]$ 
    \begin{equation*}
        \left|\frac{\norm{P_{k+1:p}X_i}_2^2}{\bE\norm{P_{k+1:p}X}_2^2}-1\right|\leq \frac{1}{4}
    \end{equation*} where $P_{k+1:p}X_1,\dots, P_{k+1:p}X_N$ are $N$ i.i.d. copies of $P_{k+1:p}X=\Sigma_{k+1:p}^{1/2}Z$ and $Z$ is a random vector with centered and variance one independent coordinates $(z_j)_j$. We therefore recover assumption \textit{a)} from Theorem~\ref{theo:stochastic_argument_weak_moments}  for $\delta=1/4$ and $\gamma$ that can be deduced from \eqref{eq:proba_theo_HS_sub-gauss} for $\Gamma = \Sigma_{k+1:p}$. Assumption \textit{d)} also follows from Theorem~\ref{theo:hanson_wright_ineq_sub-gauss}  and an union bound applied directly to the isotropic vectors $Z_i=\Sigma_{1:k}^{-1/2}X_i$'s (i.e. for $\Gamma$ which is the identity matrix $I_k$)  and, therefore, holds for $c_0=2$ and $\gamma'=\exp(-c_2k)$ as long as $\log N \lesssim k$. 
\subsection{Proof of Proposition~8 from the main document: a Hanson-Wright type inequality under weak moments assumption from \cite{latala_bounding_2022}}\label{subsec:HW_Latala}
In this section, we suppose that there is for some constant $\kappa>0$ such that for all $2\leq q\leq R\log{(e N)}$ and for every $\vv\in\bR^p$ we have $\norm{\left<X,\vv\right>}_{L_q}\leq \kappa\norm{\Sigma^{1/2}\vv}_2$. Consequently, assumptions \textit{b)} and \textit{c)} from Theorem~\ref{theo:stochastic_argument_weak_moments} and \ref{theo:proposition_5_main} are satisfied, and we must now verify assumptions \textit{a)} and \textit{d)} under this assumption.

We assume that $X=\Sigma^{1/2}Z$ where $Z$ is a random vector with independent and variance one coordinates $(z_j)_j$. As a result of $L_{R\log{(e N)}}-L_2$ equivalence above, the coordinates $(Z_j)_{j}$ satisfy $\max_j\norm{Z_j}_{L_{4+4\varepsilon}}\leq \kappa$ for some $\varepsilon>0$. We denote by $\norm{A}_{HS} = \sqrt{\Tr(A^\top A)}$ the Hilbert-Schmidt norm of a matrix $A$. In this part, we employ the following lemma from \cite{latala_bounding_2022}.

\begin{Lemma}\label{lemma:latala_diagonal}\cite[Lemma 14]{latala_bounding_2022}
    Let $Z = (z_j)_j$ be a random vector with independent mean zero and variance $1$ real-valued coordinates. We assume that for all $j$'s,  $\norm{z_j}_{L_r}\leq \kappa$ for some $\kappa>0$ and $r>4$. Then, there exists some absolute constant $C_\kappa$ (depending only on $\kappa$) such that for any  symmetric matrix $\Gamma \in\bR^{p\times p}$,
    \begin{equation*}
        \norm{\norm{\Gamma^{1/2}Z}_2^2 - \norm{\Gamma^{1/2}}_{\mathrm{HS}}^2}_{L_{r/2}} \leq C_\kappa \norm{\Gamma}_{\mathrm{HS}}.
    \end{equation*}
\end{Lemma}
By Lemma~\ref{lemma:latala_diagonal} and Markov's inequality, for all $t>0$, we have
\begin{align*}
    \bP\left(\left|\norm{\Gamma^{1/2}Z}_2^2 - \norm{\Gamma^{1/2}}_{\mathrm{HS}}^2\right|  \geq t \norm{\Gamma^{1/2}}_{\mathrm{HS}}^2\right) & \leq \frac{\norm{\norm{\Gamma^{1/2}Z}_2^2 - \norm{\Gamma^{1/2}}_{\mathrm{HS}}^2}_{L_{r/2}}^{r/2}}{\left(t \norm{\Gamma^{1/2}}_{\mathrm{HS}}^2\right)^{r/2}} \\
    &\leq \left(\frac{C_\kappa\norm{\Gamma}_{HS}}{t \norm{\Gamma^{1/2}}_{\mathrm{HS}}^2}\right)^{r/2} 
    \leq \left(\frac{C_\kappa}{t}\cdot\frac{\sqrt{\Tr(\Gamma^2)}}{\Tr\Gamma}\right)^{r/2}.
\end{align*} When the later inequality is applied to  $\Gamma = \Sigma_{k+1:p}$, we get that, if the condition $N \Tr(\Sigma_{k+1:p}^2)\lesssim \Tr(\Sigma_{k+1:p})$ - which is weaker than the necessary condition $N \Tr(\Sigma_{k+1:p}^2)=o(\Tr(\Sigma_{k+1:p}))$ for benign overfitting - is satisfied then 
\begin{equation*}
    \bP\left(\left|\norm{\Sigma_{k+1:p}^{1/2}Z}_2^2 - \norm{\Sigma_{k+1:p}^{1/2}}_{\mathrm{HS}}^2\right|\geq \frac{1}{4} \norm{\Sigma_{k+1:p}^{1/2}}_{\mathrm{HS}}^2\right) \leq \left(\frac{C_\kappa'}{N}\right)^{r/4}
\end{equation*}and so for $r = 4+4\varepsilon$ for some $\eps>0$, a union bound yields to
\begin{equation*}
    \max_{i\in[N]}\left|\norm{\Sigma_{k+1:p}^{1/2}Z_i}_2^2 - \norm{\Sigma_{k+1:p}^{1/2}}_{\mathrm{HS}}^2\right|\geq \frac{1}{4} \norm{\Sigma_{k+1:p}^{1/2}}_{\mathrm{HS}}^2
\end{equation*} with probability at least $1-(c_0/N)^\eps$ for some absolute constant $c_0$. Since $\norm{\Sigma_{k+1:p}^{1/2}}_{\mathrm{HS}}^2 = \bE \norm{\Sigma_{k+1:p}^{1/2}Z}_2^2$ we get that 
assumption \textit{a)} of Theorem~\ref{theo:stochastic_argument_weak_moments} holds with $\delta = 1/4$ and $\gamma = (c_0/N)^\varepsilon$. Assumption \textit{d)} of Theorem~\ref{theo:stochastic_argument_weak_moments} follows by a similar approach and so we omit it. We showed that the following result holds which is a restatement of Proposition~8 in the main document.

\begin{Proposition}
Let $X$ be a random vector in $\bR^p$ and denote by $X_1, \ldots, X_N$ i.i.d. copies of $X$. We assume that there exists $Z$ a random vector in $\bR^p$ with independent coordinates $z_1, \ldots, z_p$ so that $X=\Sigma^{1/2}Z$. We also assume that there exists $\kappa>0$ and $r>0$ so that for all $j\in[p]$,  $\norm{z_j}_{L_r}\leq \kappa$. Then, there exists an absolute constant $c_1$ depending only on $\kappa$ such that assumptions \textit{a)} and \textit{d)} from Theorem~\ref{theo:stochastic_argument_weak_moments} and \ref{theo:proposition_5_main} hold with $\delta=1/4$, $c_0=2$ and $\gamma = \gamma' = (c_1/N)^{(r-4)/4}$. 
\end{Proposition}


\section{Other technical tools used to handle the heavy-tailed case. Proof of Proposition~6 and Proposition~7} 
\label{sec:other_technical_tools_used_to_handle_the_heavy_tailed_case}
In this final section, we gather all the other tools we need to obtain the main result in the heavy-tailed case i.e. Theorem~11 from the main document. Indeed, in this section, we provide proofs for the last two results  from Section~6 in the main document which are  Proposition~6 and Proposition~7 in this companion document.

\subsection{Extending property~4 to the heavy-tailed case} 
\label{sub:extending_property_4_to_the_heavy_tailed_case}
The last property of the design matrix we have to extend from the Gaussian case to the heavy-tailed case is
\begin{itemize}
    \item[4.] for $D = X_{k+1:p}^\top A$, $\Tr\left(DD^\top\right)\leq \frac{c_0 N\Tr\left(\Sigma_{k+1:p}^2\right)}{\left(\Tr\left(\Sigma_{k+1:p}\right)\right)^2}$.
\end{itemize}

Property~\textit{4.} is the final property we need on $\bX$ for the proof of Proposition~4 in the main document. It requires information on the sum of i.i.d. random variables.  To deal with this type of quantities under weak moment assumptions we will use a result (see Lemma 3.2) from \cite{MR3565471} on the $L_r$-norm of a sum of i.i.d. random variables to deduce the following result.

\begin{Lemma}\label{lem:moment-sum-variid-positive}
    There exists absolute constant $c_0,c_1,c_2$ such that the following holds. Let $1\leq r<q$, set $Z\in L_q$ and put $Z_1,\cdots,Z_N$ to be independent copies of $Z$. Fix $1\leq p\leq N$, let $j_0 = \lceil (c_0 p)/\left( ((q/r)-1)\log{(4+eN/p)} \right) \rceil$ and $t>2$.
    If $j_0 = 1$ and $0<\beta<(q/r)-1$ then with probability at least $1-c_2t^{-q}N^{-\beta}$,
    \begin{align*}
        \left(\sum_{j=1}^N \left|Z_i\right|^r\right)^{1/r} \leq c_1\left(\frac{q}{q-(\beta+1)r}\right)^{1/r}t\norm{Z}_{L_q}N^{1/r}.
    \end{align*}
\end{Lemma}

Lemma~\ref{lem:moment-sum-variid-positive} shows that even under a  weak moment assumption, the $\ell_r$ norm of a random vector with i.i.d. coordinates is well-concentrated. We are now combining Lemma~\ref{lem:moment-sum-variid-positive} and Lemma~\ref{lemma:latala_diagonal}

\begin{Proposition}
  We assume that $X = \Sigma^{1/2}Z$ where $Z$ is a random vector with independent and variance one coordinates $(Z_j)_{j}$, and $X$ satisfies the $L_{R\log{eN}}-L_2$ equivalence assumption as in Appendix~\ref{subsec:HW_Latala}.
    Then with probability at least $1-c_2N^{-1/2}$ with $c_2$ borrowed from Lemma~\ref{lem:moment-sum-variid-positive}, point~\textit{4.} holds true.
\end{Proposition}
\beginproof Notice that that by point 2,
\begin{align*}
    \Tr\left(DD^\top\right)=\Tr\left(D^\top D\right)\leq \frac{c\Tr\left(X_{k+1:p}\Sigma X_{k+1:p}\right)}{\left(\Tr(\Sigma_{k+1:p})\right)^2}=\frac{c\sum_{i=1}^N \norm{\Sigma^{1/2}P_{k+1:p}X_i}_2 ^2}{\left(\Tr(\Sigma_{k+1:p})\right)^2},
\end{align*}we are left with the upper estimate for $\sum_{i=1}^N \norm{\Sigma^{1/2}P_{k+1:p}X_i}_2 ^2 = \sum_{i=1}^N \norm{\Sigma_{k+1:p}Z_i}_2^2$. We know that $\max_j \norm{Z_j}_{L_{4+4\eps}}\leq \kappa$ for some $\eps>0$. Applying Lemma~\ref{lemma:latala_diagonal} with $Z = Z$ and $\Gamma = \Sigma_{k+1:p}^2$. There exists an absolute constant $C_\kappa$ depending only on $\kappa$ such that $\norm{\norm{\Sigma_{k+1:p}Z}_2^2 - \norm{\Sigma_{k+1:p}}_{HS}^2}_{L_{2+2\eps}}\leq C_\kappa \norm{\Sigma_{k+1:p}^2}_{HS}$. By $\norm{\cdot}_4\leq\norm{\cdot}_2$, $\norm{\Sigma_{k+1:p}^2}_{HS} =  \left((\sum_{j=k+1}^p \sigma_j^4)^{1/4}\right)^2 \leq \sum_{j=k+1}^p \sigma_j^2 = \norm{\Sigma_{k+1:p}}_{HS}^2$. Therefore, by H{\"o}lder's inequality,
\begin{align*}
    \norm{\norm{\Sigma_{k+1:p}Z}_2}_{L_4}\leq \left(\bE\norm{\Sigma_{k+1:p}Z}_2^{2(2+2\eps)}\right)^{\frac{1}{2(2+2\eps)}}\leq \left((1+C_\kappa)\right)^{\frac{1}{2}}\norm{\Sigma_{k+1:p}}_{HS}.
\end{align*}Now we apply Lemma~\ref{lem:moment-sum-variid-positive} with $p=1$, $r=2$, $q=4$, $\beta = 1/2$, then $j_0 = \lceil c_0/\log{(4+eN)} \rceil = 1$ provided that $N>(2/e)(e^{c_0}-4)$. With probability at least $1 - c_2N^{-1/2}$,
\begin{align*}
    \left(\sum_{j=1}^N\norm{\Sigma_{k+1:p}Z_i}_2^2\right)^{1/2}\leq 4c_1\sqrt{N}\norm{ \norm{\Sigma_{k+1:p}Z}_2 }_{L_4}\leq 4c_1\sqrt{1+C_\kappa}\sqrt{N}\norm{\Sigma_{k+1:p}}_{HS}.
\end{align*}As a result, $\Tr\left(DD^\top\right)\leq \frac{16cc_1^2(1+C_\kappa)N\Tr\left(\Sigma_{k+1:p}^2\right)}{\left(\Tr\left(\Sigma_{k+1:p}\right)\right)^2}$.

{{\mbox{}\nolinebreak\hfill\rule{2mm}{2mm}\par\medbreak}}


\subsection{Checking property~5 on the noise under a heavy-tailed assumption} 
\label{sub:property_on_the_noise}
In this section, we deal with the only property that the noise vector $\bxi$ needs to satisfy: if $D$ is a $p\times N$ matrix then with large probability $\norm{D \bxi}_2\leq c_0\sigma_\xi\sqrt{\Tr(DD^\top)}$ -- which is Property~\textit{5.} from Section~6.1 in the main document. 

In fact, we already obtained such a result thanks to Lemma~14 from \cite{latala_bounding_2022} which is recalled in Lemma~\ref{lemma:latala_diagonal}  -- note that the result in Lemma~14 from \cite{latala_bounding_2022} hold for any square matrix (not necessarily symmetric), it is however straightforward to check that it is also true for any rectangular matrix. Following the same argument as in Section~\ref{subsec:HW_Latala}, we obtain the next result which is a restatement of Proposition~7 in the main document.

\begin{Lemma}\label{lemma:latala_diagonal_noise}
    Let $\bxi = (\xi_i)_{i=1}^N$ be a random vector with independent mean zero and variance $\sigma_\xi$ real-valued coordinates. We assume that for all $i$'s,  $\norm{\xi_i}_{L_r}\leq \kappa\sigma_\xi$ for some $\kappa>0$ and $r>4$. Then, there exists some absolute constant $C_\kappa$ (depending only on $\kappa$) such that for any matrix $D \in\bR^{p\times N}$ the following holds: if for some integer $k$ for which  $\sqrt{k} \norm{D}_{op}\leq \sqrt{\Tr(DD^\top)}$ then with probability at least $1-(C_\kappa/k)^{r/4}$,
    \begin{equation*}
        \norm{D\bxi}_2 \leq  (3/2) \sigma_\xi \sqrt{\Tr(DD^\top)}.
    \end{equation*}
\end{Lemma}

\beginproof
Without loss of generality we assume that $\sigma_\xi =1$. It follows from Lemma~14 from \cite{latala_bounding_2022} (with the direct extension to the rectangular case) that $\norm{\norm{D\bxi}_2^2 - \norm{D}_{\mathrm{HS}}^2}_{L_{r/2}} \leq C_\kappa \norm{D^\top D}_{\mathrm{HS}}$. Hence, by Markov's inequality, for all $t>0$,
\begin{align*}
\bP\left[\left|\norm{D\bxi}_2^2 - \norm{D}_{\mathrm{HS}}^2\right|\geq t\norm{D}_{\mathrm{HS}}^2 \right] \leq \left(\frac{C_\kappa \norm{D^\top D}_{\mathrm{HS}}}{t \norm{D}_{\mathrm{HS}}^2 }\right)^{r/2}\leq \left(\frac{C_\kappa}{t\sqrt{k}}\right)^{r/2}   
\end{align*}  where we used that $\norm{D^\top D}_{\mathrm{HS}}\leq \norm{D}_{op}\norm{D}_{\mathrm{HS}}$  and that $\sqrt{k} \norm{D}_{op}\leq \norm{D}_{\mathrm{HS}}$. The result follows for $t=1/2$.
{{\mbox{}\nolinebreak\hfill\rule{2mm}{2mm}\par\medbreak}}

\section{Proof of Theorem~11}\label{sec:proof_theo_11}
 The proof of Theorem~11 follows the same strategy as the one of Theorem~\ref{theo:main} for $k=k^*_b$ with $b=4/\kappa_{DM}$  in the case where $k^*_b\leq \kappa_{iso} N$ where $\kappa_{iso}$ is the absolute constant $C$ from  Theorem~\ref{theo:tikhomirov_isomorphy_2} and so $\sigma_1 N\geq \kappa_{DM}\ell_*^2(\Sigma_{k^*_b+1}B_2^p)$. In that case and using the same notation as in Theorem~\ref{theo:main}, we have $t=k^*_b$, $\sigma^2(\square, \triangle)=\square$ and on the event where points \textit{1.} to \textit{4.} (from Section~6.1 in the main document) hold, $\sqrt{\Tr(D D^\top)}\gtrsim \sqrt{k^*_b}\norm{D}_{op}$ for $D = \tilde \Sigma_1^{-1/2}X_{1:k}^\top (X_{k+1:p}X_{k+1:p}^\top)^{-1}$ or $D = \Sigma_{k+1:p}^{1/2}A$ with  $A=X_{k+1:p}^\top (X_{k+1:p} X_{k+1:p}^\top)^{-1}$. Hence, one can apply Lemma~\ref{lemma:latala_diagonal_noise} (see also Proposition~7 in the main document) with $k=k^*_b$ and get the properties we need on the noise, i.e. Property~\textit{5.} in the main document. Finally, all the properties we need on the design matrix $\bX$, i.e. properties \textit{1.} to \textit{4.} in Section~6.1 in the main document are implied by Proposition~4,5 and 6 in the same document.




\begin{footnotesize}
\bibliographystyle{plain}
\bibliography{biblio}
\end{footnotesize}

\end{document}

%% file: packages_notes.tex
\usepackage{graphicx}
\usepackage{color}
\usepackage{amsmath}
\usepackage{amssymb}
\usepackage{amscd}
\usepackage{bbm}
\usepackage{tikz}
\usetikzlibrary{patterns}
\usetikzlibrary{decorations.pathmorphing}
\usetikzlibrary{shapes.geometric, arrows}
\usepackage{hyperref}
\hypersetup{
    bookmarks=true,         
    colorlinks=true,       
    linkcolor=red,          
    citecolor=green,        
    filecolor=magenta,      
    urlcolor=cyan           
}
\usepackage{xcolor}
\usepackage{eurosym}
\usepackage{float}
\usepackage{adjustbox}

\usepackage[utf8]{inputenc}
\usepackage{amsthm}
\usepackage{bbm} 
\usepackage[linesnumbered,lined,boxed,commentsnumbered]{algorithm2e}

\usepackage{marginnote}

\usepackage{changepage}

\usepackage{enumitem}

\usepackage{multirow,bigdelim}
\usepackage{stackengine}
\usepackage{glossaries}	

\newtheorem{Theorem}{Theorem}

\newtheorem{Definition}{Definition}
\newtheorem{Lemma}{Lemma}
\newtheorem{Proposition}{Proposition}
\newtheorem{Corollary}{Corollary}
\newtheorem{Remark}{Remark}

\newenvironment{itemize*}%
  {\begin{itemize}%
    \setlength{\itemsep}{0pt}%
    \setlength{\parskip}{0pt}}%
  {\end{itemize}}


%% file: commandes_guillaume.tex
\newcommand{\inr}[1]{\bigl< #1 \bigr>}

\newcommand{\norm}[1]{\left\|#1\right\|}%
\newcommand{\vertiii}[1]{{\left\vert\kern-0.25ex\left\vert\kern-0.25ex\left\vert #1 
    \right\vert\kern-0.25ex\right\vert\kern-0.25ex\right\vert}}

\newcommand{\vertiiii}[1]{{\left\vert\kern-0.25ex\left\vert\kern-0.25ex\left\vert\kern-0.25ex\left\vert #1 
    \right\vert\kern-0.25ex\right\vert\kern-0.25ex\right\vert\kern-0.25ex\right\vert}}

\newcommand\eps{\epsilon}

\newcommand{\beginproof}{{\bf Proof. {\hspace{0.2cm}}}}

\DeclareMathOperator*{\argmin}{argmin}

\DeclareMathOperator{\diag}{diag}

\DeclareMathOperator{\conv}{conv}

\def\ds1{\textrm{1\kern-0.25emI}} 

\newcommand{\1}{\ensuremath{\mathbbm{1}}}

\newcommand \R{\mathbb{R}}

\newcommand \cB{{\cal B}}
\newcommand \cC{{\cal C}}

\newcommand \cE{{\cal E}}
\newcommand \cF{{\cal F}}

\newcommand \cL{{\cal L}}
\newcommand \cM{{\cal M}}
\newcommand \cN{{\cal N}}
\newcommand \cO{{\cal O}}

\newcommand \cQ{{\cal Q}}
\newcommand \cR{{\cal R}}
\newcommand \cS{{\cal S}}

\newcommand \cY{{\cal Y}}

\newcommand \bE{{\mathbb E}}

\newcommand \bG{{\mathbb G}}

\newcommand \bN{{\mathbb N}}

\newcommand \bP{{\mathbb P}}

\newcommand \bR{{\mathbb R}}

\newcommand \bX{{\mathbb X}}
\newcommand \bY{{\mathbb Y}}

\newcommand{\vlambda}{{\boldsymbol{\lambda}}}
\newcommand{\vv}{{\boldsymbol{v}}}

\newcommand{\vx}{{\boldsymbol{x}}}
\newcommand{\vz}{{\boldsymbol{z}}}

\newcommand{\vu}{{\boldsymbol{u}}}

\newcommand{\vy}{\mathbf{y}}

\DeclareMathOperator*{\Tr}{Tr}
\DeclareMathOperator*{\diam}{diam}

\newcommand{\bxi}{{\boldsymbol{\xi}}}
\newcommand{\bbeta}{{\boldsymbol{\beta}}}
\newcommand{\blambda}{{\boldsymbol{\lambda}}}

\newcommand{\bv}{\boldsymbol{v}}

%% file: ell_2_PTRF.bbl
\begin{thebibliography}{10}

\bibitem{bartl_random_2022}
Daniel Bartl and Shahar Mendelson.
\newblock Random embeddings with an almost {Gaussian} distortion.
\newblock {\em Advances in Mathematics}, 400:108261, May 2022.

\bibitem{MR4263288}
Peter~L. Bartlett, Philip~M. Long, G\'{a}bor Lugosi, and Alexander Tsigler.
\newblock Benign overfitting in linear regression.
\newblock {\em Proc. Natl. Acad. Sci. USA}, 117(48):30063--30070, 2020.

\bibitem{BMR}
Peter~L. Bartlett, Andreas Montanari, and Alexander Rakhlin.
\newblock Deep learning: a statistical viewpoint.
\newblock Technical report, ArXiv 2103.09177, 2021.

\bibitem{bednorz}
Wiltord Bednorz.
\newblock Concentration via chaining method and its applications.
\newblock 2014.

\bibitem{DBLP:journals/corr/abs-2105-14368}
Mikhail Belkin.
\newblock Fit without fear: remarkable mathematical phenomena of deep learning
  through the prism of interpolation.
\newblock {\em CoRR}, abs/2105.14368, 2021.

\bibitem{MR3997901}
Mikhail Belkin, Daniel Hsu, Siyuan Ma, and Soumik Mandal.
\newblock Reconciling modern machine-learning practice and the classical
  bias-variance trade-off.
\newblock {\em Proc. Natl. Acad. Sci. USA}, 116(32):15849--15854, 2019.

\bibitem{DBLP:journals/simods/BelkinHX20}
Mikhail Belkin, Daniel Hsu, and Ji~Xu.
\newblock Two models of double descent for weak features.
\newblock {\em {SIAM} J. Math. Data Sci.}, 2(4):1167--1180, 2020.

\bibitem{DBLP:conf/nips/BelkinHM18}
Mikhail Belkin, Daniel~J. Hsu, and Partha Mitra.
\newblock Overfitting or perfect fitting? risk bounds for classification and
  regression rules that interpolate.
\newblock In Samy Bengio, Hanna~M. Wallach, Hugo Larochelle, Kristen Grauman,
  Nicol{\`{o}} Cesa{-}Bianchi, and Roman Garnett, editors, {\em Advances in
  Neural Information Processing Systems 31: Annual Conference on Neural
  Information Processing Systems 2018, NeurIPS 2018, December 3-8, 2018,
  Montr{\'{e}}al, Canada}, pages 2306--2317, 2018.

\bibitem{DBLP:conf/aistats/BelkinRT19}
Mikhail Belkin, Alexander Rakhlin, and Alexandre~B. Tsybakov.
\newblock Does data interpolation contradict statistical optimality?
\newblock In Kamalika Chaudhuri and Masashi Sugiyama, editors, {\em The 22nd
  International Conference on Artificial Intelligence and Statistics, {AISTATS}
  2019, 16-18 April 2019, Naha, Okinawa, Japan}, volume~89 of {\em Proceedings
  of Machine Learning Research}, pages 1611--1619. {PMLR}, 2019.

\bibitem{CaT06}
Emmanuel~J. Candes and Terence Tao.
\newblock Near-optimal signal recovery from random projections: universal
  encoding strategies?
\newblock {\em IEEE Trans. Inform. Theory}, 52(12):5406--5425, 2006.

\bibitem{steph_BO}
Emmanuel Caron and St{\'e}phane Chr{\'e}tien.
\newblock A finite sample analysis of the benign overfitting phenomenon for
  ridge function estimation.
\newblock Technical report, 2020.

\bibitem{DBLP:journals/corr/abs-2105-02083}
Geoffrey Chinot, Felix Kuchelmeister, Matthias L{\"{o}}ffler, and Sara~A.
  van~de Geer.
\newblock Adaboost and robust one-bit compressed sensing.
\newblock {\em CoRR}, abs/2105.02083, 2021.

\bibitem{chinot_robust_2021}
Geoffrey Chinot, Guillaume Lecué, and Matthieu Lerasle.
\newblock Robust high dimensional learning for {Lipschitz} and convex losses,
  January 2021.
\newblock arXiv:1905.04281 [math, stat].

\bibitem{DBLP:journals/corr/abs-2012-00807}
Geoffrey Chinot, Matthias L{\"{o}}ffler, and Sara~A. van~de Geer.
\newblock On the robustness of minimum-norm interpolators.
\newblock {\em CoRR}, abs/2012.00807, 2020.

\bibitem{Chinot_Lerasle}
Geoffrey Chinot and Lerasle Matthieu.
\newblock On the robustness of the minimum l2 interpolator.
\newblock Technical report, ENSAE, CREST, IPParis, 2020.

\bibitem{dirksen2015tail}
Sjoerd Dirksen.
\newblock Tail bounds via generic chaining.
\newblock {\em Electronic Journal of Probability}, 20, 2015.

\bibitem{DBLP:conf/icml/DonhauserRSY22}
Konstantin Donhauser, Nicol{\`{o}} Ruggeri, Stefan Stojanovic, and Fanny Yang.
\newblock Fast rates for noisy interpolation require rethinking the effect of
  inductive bias.
\newblock In Kamalika Chaudhuri, Stefanie Jegelka, Le~Song, Csaba
  Szepesv{\'{a}}ri, Gang Niu, and Sivan Sabato, editors, {\em International
  Conference on Machine Learning, {ICML} 2022, 17-23 July 2022, Baltimore,
  Maryland, {USA}}, volume 162 of {\em Proceedings of Machine Learning
  Research}, pages 5397--5428. {PMLR}, 2022.

\bibitem{DBLP:journals/corr/abs-1903-08560}
Trevor Hastie, Andrea Montanari, Saharon Rosset, and Ryan~J. Tibshirani.
\newblock Surprises in high-dimensional ridgeless least squares interpolation.
\newblock {\em CoRR}, abs/1903.08560, 2019.

\bibitem{DBLP:conf/stoc/JacotGH21}
Arthur Jacot, Franck Gabriel, and Cl{\'{e}}ment Hongler.
\newblock Neural tangent kernel: convergence and generalization in neural
  networks (invited paper).
\newblock In Samir Khuller and Virginia~Vassilevska Williams, editors, {\em
  {STOC} '21: 53rd Annual {ACM} {SIGACT} Symposium on Theory of Computing,
  Virtual Event, Italy, June 21-25, 2021}, page~6. {ACM}, 2021.

\bibitem{DBLP:conf/nips/JacotHG18}
Arthur Jacot, Cl{\'{e}}ment Hongler, and Franck Gabriel.
\newblock Neural tangent kernel: Convergence and generalization in neural
  networks.
\newblock In Samy Bengio, Hanna~M. Wallach, Hugo Larochelle, Kristen Grauman,
  Nicol{\`{o}} Cesa{-}Bianchi, and Roman Garnett, editors, {\em Advances in
  Neural Information Processing Systems 31: Annual Conference on Neural
  Information Processing Systems 2018, NeurIPS 2018, December 3-8, 2018,
  Montr{\'{e}}al, Canada}, pages 8580--8589, 2018.

\bibitem{MR2149924}
B.~Klartag and S.~Mendelson.
\newblock Empirical processes and random projections.
\newblock {\em J. Funct. Anal.}, 225(1):229--245, 2005.

\bibitem{DBLP:journals/corr/abs-2106-09276}
Frederic Koehler, Lijia Zhou, Danica~J. Sutherland, and Nathan Srebro.
\newblock Uniform convergence of interpolators: Gaussian width, norm bounds,
  and benign overfitting.
\newblock {\em CoRR}, abs/2106.09276, 2021.

\bibitem{latala_bounding_2022}
Rafał Latała.
\newblock Bounding suprema of canonical processes via convex hull, April 2022.
\newblock arXiv:2204.09463 [math].

\bibitem{lecue2013learning}
Guillaume Lecu{\'e} and Shahar Mendelson.
\newblock Learning subgaussian classes: Upper and minimax bounds.
\newblock {\em arXiv preprint arXiv:1305.4825}, 2013.

\bibitem{MR3782379}
Guillaume Lecu\'{e} and Shahar Mendelson.
\newblock Regularization and the small-ball method {I}: {S}parse recovery.
\newblock {\em Ann. Statist.}, 46(2):611--641, 2018.

\bibitem{lecue_sparse_2017}
Guillaume Lecué and Shahar Mendelson.
\newblock Sparse recovery under weak moment assumptions.
\newblock {\em Journal of the European Mathematical Society}, 19(3):881--904,
  2017.

\bibitem{Led01}
Michel Ledoux.
\newblock {\em The concentration of measure phenomenon}, volume~89 of {\em
  Mathematical Surveys and Monographs}.
\newblock American Mathematical Society, Providence, RI, 2001.

\bibitem{MR2814399}
Michel Ledoux and Michel Talagrand.
\newblock {\em Probability in {B}anach spaces}.
\newblock Classics in Mathematics. Springer-Verlag, Berlin, 2011.
\newblock Isoperimetry and processes, Reprint of the 1991 edition.

\bibitem{MR4124325}
Tengyuan Liang and Alexander Rakhlin.
\newblock Just interpolate: kernel ``ridgeless'' regression can generalize.
\newblock {\em Ann. Statist.}, 48(3):1329--1347, 2020.

\bibitem{DBLP:journals/corr/abs-2101-11815}
Tengyuan Liang and Benjamin Recht.
\newblock Interpolating classifiers make few mistakes.
\newblock {\em CoRR}, abs/2101.11815, 2021.

\bibitem{DBLP:conf/nips/0001ZB20}
Chaoyue Liu, Libin Zhu, and Mikhail Belkin.
\newblock On the linearity of large non-linear models: when and why the tangent
  kernel is constant.
\newblock In Hugo Larochelle, Marc'Aurelio Ranzato, Raia Hadsell,
  Maria{-}Florina Balcan, and Hsuan{-}Tien Lin, editors, {\em Advances in
  Neural Information Processing Systems 33: Annual Conference on Neural
  Information Processing Systems 2020, NeurIPS 2020, December 6-12, 2020,
  virtual}, 2020.

\bibitem{Mei_Montanari}
Song Mei and Andreas Montanari.
\newblock The generalization error of random features regression: Precise
  asymptotics and double descent curve.
\newblock Technical report, Stanford University, 2020.

\bibitem{MR3565471}
Shahar Mendelson.
\newblock Upper bounds on product and multiplier empirical processes.
\newblock {\em Stochastic Process. Appl.}, 126(12):3652--3680, 2016.

\bibitem{mendelson_multiplier_2017}
Shahar Mendelson.
\newblock On {Multiplier} {Processes} {Under} {Weak} {Moment} {Assumptions}.
\newblock In Bo'az Klartag and Emanuel Milman, editors, {\em Geometric
  {Aspects} of {Functional} {Analysis}: {Israel} {Seminar} ({GAFA})
  2014–2016}, Lecture {Notes} in {Mathematics}, pages 301--318. Springer
  International Publishing, Cham, 2017.

\bibitem{MR2373017}
Shahar Mendelson, Alain Pajor, and Nicole Tomczak-Jaegermann.
\newblock Reconstruction and subgaussian operators in asymptotic geometric
  analysis.
\newblock {\em Geom. Funct. Anal.}, 17(4):1248--1282, 2007.

\bibitem{MR0293374}
V.~D. Milman.
\newblock A new proof of {A}. {D}voretzky's theorem on cross-sections of convex
  bodies.
\newblock {\em Funkcional. Anal. i Prilo\v{z}en.}, 5(4):28--37, 1971.

\bibitem{DBLP:journals/corr/abs-1906-03667}
Partha~P. Mitra.
\newblock Understanding overfitting peaks in generalization error: Analytical
  risk curves for l\({}_{\mbox{2}}\) and l\({}_{\mbox{1}}\) penalized
  interpolation.
\newblock {\em CoRR}, abs/1906.03667, 2019.

\bibitem{DBLP:journals/jsait/MuthukumarVSS20}
Vidya Muthukumar, Kailas Vodrahalli, Vignesh Subramanian, and Anant Sahai.
\newblock Harmless interpolation of noisy data in regression.
\newblock {\em {IEEE} J. Sel. Areas Inf. Theory}, 1(1):67--83, 2020.

\bibitem{MR1036275}
Gilles Pisier.
\newblock {\em The volume of convex bodies and {B}anach space geometry},
  volume~94 of {\em Cambridge Tracts in Mathematics}.
\newblock Cambridge University Press, Cambridge, 1989.

\bibitem{DBLP:conf/aistats/RichardsMR21}
Dominic Richards, Jaouad Mourtada, and Lorenzo Rosasco.
\newblock Asymptotics of ridge(less) regression under general source condition.
\newblock In Arindam Banerjee and Kenji Fukumizu, editors, {\em The 24th
  International Conference on Artificial Intelligence and Statistics, {AISTATS}
  2021, April 13-15, 2021, Virtual Event}, volume 130 of {\em Proceedings of
  Machine Learning Research}, pages 3889--3897. {PMLR}, 2021.

\bibitem{MR1673273}
Robert~E. Schapire, Yoav Freund, Peter Bartlett, and Wee~Sun Lee.
\newblock Boosting the margin: a new explanation for the effectiveness of
  voting methods.
\newblock {\em Ann. Statist.}, 26(5):1651--1686, 1998.

\bibitem{MR3184689}
Michel Talagrand.
\newblock {\em Upper and lower bounds for stochastic processes}, volume~60 of
  {\em Ergebnisse der Mathematik und ihrer Grenzgebiete. 3. Folge. A Series of
  Modern Surveys in Mathematics [Results in Mathematics and Related Areas. 3rd
  Series. A Series of Modern Surveys in Mathematics]}.
\newblock Springer, Heidelberg, 2014.
\newblock Modern methods and classical problems.

\bibitem{tikhomirov_sample_2018}
Konstantin Tikhomirov.
\newblock Sample {Covariance} {Matrices} of {Heavy}-{Tailed} {Distributions}.
\newblock {\em International Mathematics Research Notices},
  2018(20):6254--6289, October 2018.

\bibitem{TB21}
Alexander Tsigler and Peter~L. Bartlett.
\newblock Benign overfitting in ridge regression.
\newblock 2021.

\bibitem{MR3837109}
Roman Vershynin.
\newblock {\em High-dimensional probability}, volume~47 of {\em Cambridge
  Series in Statistical and Probabilistic Mathematics}.
\newblock Cambridge University Press, Cambridge, 2018.
\newblock An introduction with applications in data science, With a foreword by
  Sara van de Geer.

\bibitem{vershynin_high-dimensional_2018}
Roman Vershynin.
\newblock {\em High-{Dimensional} {Probability}: {An} {Introduction} with
  {Applications} in {Data} {Science}}.
\newblock Cambridge {Series} in {Statistical} and {Probabilistic}
  {Mathematics}. Cambridge University Press, Cambridge, 2018.

\bibitem{DBLP:journals/corr/abs-2111-05987}
Guillaume Wang, Konstantin Donhauser, and Fanny Yang.
\newblock Tight bounds for minimum l1-norm interpolation of noisy data.
\newblock {\em CoRR}, abs/2111.05987, 2021.

\bibitem{DBLP:journals/jmlr/WynerOBM17}
Abraham~J. Wyner, Matthew Olson, Justin Bleich, and David Mease.
\newblock Explaining the success of adaboost and random forests as
  interpolating classifiers.
\newblock {\em J. Mach. Learn. Res.}, 18:48:1--48:33, 2017.

\bibitem{DBLP:conf/nips/XuH19}
Ji~Xu and Daniel~J. Hsu.
\newblock On the number of variables to use in principal component regression.
\newblock In Hanna~M. Wallach, Hugo Larochelle, Alina Beygelzimer, Florence
  d'Alch{\'{e}}{-}Buc, Emily~B. Fox, and Roman Garnett, editors, {\em Advances
  in Neural Information Processing Systems 32: Annual Conference on Neural
  Information Processing Systems 2019, NeurIPS 2019, December 8-14, 2019,
  Vancouver, BC, Canada}, pages 5095--5104, 2019.

\bibitem{DBLP:conf/iclr/ZhangBHRV17}
Chiyuan Zhang, Samy Bengio, Moritz Hardt, Benjamin Recht, and Oriol Vinyals.
\newblock Understanding deep learning requires rethinking generalization.
\newblock In {\em 5th International Conference on Learning Representations,
  {ICLR} 2017, Toulon, France, April 24-26, 2017, Conference Track
  Proceedings}. OpenReview.net, 2017.

\bibitem{DBLP:journals/cacm/ZhangBHRV21}
Chiyuan Zhang, Samy Bengio, Moritz Hardt, Benjamin Recht, and Oriol Vinyals.
\newblock Understanding deep learning (still) requires rethinking
  generalization.
\newblock {\em Commun. {ACM}}, 64(3):107--115, 2021.

\bibitem{DBLP:journals/corr/abs-2112-04470}
Lijia Zhou, Frederic Koehler, Danica~J. Sutherland, and Nathan Srebro.
\newblock Optimistic rates: {A} unifying theory for interpolation learning and
  regularization in linear regression.
\newblock {\em CoRR}, abs/2112.04470, 2021.

\bibitem{DBLP:conf/colt/ZouWBGK21}
Difan Zou, Jingfeng Wu, Vladimir Braverman, Quanquan Gu, and Sham~M. Kakade.
\newblock Benign overfitting of constant-stepsize {SGD} for linear regression.
\newblock In Mikhail Belkin and Samory Kpotufe, editors, {\em Conference on
  Learning Theory, {COLT} 2021, 15-19 August 2021, Boulder, Colorado, {USA}},
  volume 134 of {\em Proceedings of Machine Learning Research}, pages
  4633--4635. {PMLR}, 2021.

\end{thebibliography}
